\documentclass[11pt]{article}
\usepackage[utf8]{inputenc}
\usepackage[T1]{fontenc}
\usepackage{amsmath,amssymb}
\usepackage{amsthm}
\usepackage{rotating,pdflscape,afterpage}
\usepackage[sc,osf]{mathpazo}

	\AtBeginDocument{
		\DeclareSymbolFont{AMSb}{U}{msb}{m}{n}
		\DeclareSymbolFontAlphabet{\mathbb}{AMSb}
	}
\usepackage{enumitem}
\usepackage{mathtools}
\usepackage[title]{appendix}
\newcommand{\abyss}[1]{}
\usepackage{stmaryrd}
\expandafter\def\csname opt@stmaryrd.sty\endcsname
{only,shortleftarrow,shortrightarrow}
\SetSymbolFont{stmry}{bold}{U}{stmry}{m}{n}

\makeatletter
\let\@wraptoccontribs\wraptoccontribs
\makeatother

\usepackage[font=small,labelfont=bf]{caption}

\mathtoolsset{mathic}

\newlist{enumarabic}{enumerate}{1}
\setlist[enumarabic]{font = \normalfont,label = (\arabic*),leftmargin = 0.3in}
\newlist{enumroman}{enumerate}{1}
\setlist[enumroman]{font = \normalfont,label = (\roman*),leftmargin = 0.3in}

\makeatletter
\g@addto@macro\bfseries{\boldmath}
\makeatother

\def\cf{\emph{cf.}}

\def\eg{\emph{e.g.}}
\def\ie{\emph{i.e.}}

\def\otimesunder#1{\mathbin{\mathop\otimes\limits_{\mkern-20mu #1\mkern-20mu}}}

\numberwithin{equation}{section}

\allowdisplaybreaks[4]

\usepackage{graphicx}
\usepackage{extpfeil}
\usepackage{mathabx} 
\usepackage{fancyvrb}
\usepackage{float}
\newcommand{\mockalph}[1]{\!}
\usepackage{enumitem}
\usepackage[nouppercase]{scrlayer-scrpage}
\usepackage{avant}
\usepackage{mathrsfs}
\usepackage{mathtools}
\usepackage{faktor,xfrac}
\usepackage[all]{xypic}\xyoption{rotate}
	\newdir{ >}{{}*!/-7pt/\dir{>}}
	\newdir{i}{{}*!/-7pt/\dir{(}	}
\usepackage{thmtools}

\usepackage[hyphens]{url}		
\usepackage{hyperref}	
\usepackage{cleveref} 	

\makeatletter
	\let\c@figure\c@table
	\let\c@equation\c@table
\makeatother

\numberwithin{table}{section}
\numberwithin{figure}{section}
\newtheorem{abcthm}[table]{Theorem}

\newtheorem{theorem}[table]{Theorem}
\newtheorem{proposition}[table]{Proposition}

\newtheorem{corollary}[table]{Corollary}

\newtheorem{lemma}[table]{Lemma}
\newtheorem{prop}[table]{Proposition}

\newtheorem{convention}[table]{Convention}

\theoremstyle{definition}

\newtheorem{definition}[table]{Definition}

\newtheorem{notation}[table]{Notation}

\newtheorem{discussion}[table]{Discussion}

\newtheorem*{layout*}{Outline}
\newtheorem{sublayout}{}[table]

\newtheorem{notmain}{Notation}[table]
\newtheorem{subdef}{}[table]

\theoremstyle{remark}

\newtheorem{example}[table]{Example}

\newtheorem{remark}[table]{Remark}

\theoremstyle{plain}
\newtheorem*{thm*}{Theorem}
\newtheorem*{theorem*}{Theorem}
\newtheorem*{prop*}{Proposition}
\newtheorem*{proposition*}{Proposition}
\newtheorem*{lemma*}{Lemma}
\newtheorem*{corollary*}{Corollary}
\newtheorem*{cor*}{Corollary}

\theoremstyle{definition}
\newtheorem*{definition*}{Definition}
\newtheorem*{defn*}{Definition}
\newtheorem*{QQ*}{Question}
\newtheorem*{obs*}{Observation}
\newtheorem*{notation*}{Notation}

\theoremstyle{remark}
\newtheorem*{rmk*}{Remark}
\newtheorem*{remark*}{Remark}

\newtheorem*{examples*}{Examples}
\newtheorem*{example*}{Example}
\newtheorem*{EG*}{Example}
\newtheorem*{EGs*}{Examples}
\newtheorem*{fact*}{Fact}
\newtheorem*{prob*}{Problem}

\usepackage{etoolbox}
\usepackage[usenames,dvipsnames,svgnames,table]{xcolor}

\usepackage[disable]{todonotes}	

\newcommand		{\defd}[1]	{\textcolor{RoyalBlue}{\textbf{\textit{#1}}}}
\newcommand		{\defm}[1]	{\textcolor{RoyalBlue}{#1}}

\tracingpatches
\makeatletter
\patchcmd{\@setref}{\bfseries ??}{\bfseries\color{red} FIX ME!}{}{}
\patchcmd{\@setcref}         {??}{\color{red} FIX ME!}{}{}
\patchcmd{\@setcref}         {??}{\color{red} FIX ME!}{}{}
\patchcmd{\@setcrefrange}    {??}{\color{red} FIX ME!}{}{}
\patchcmd{\@setcrefrange}    {??}{\color{red} FIX ME!}{}{}
\patchcmd{\@setcrefrange}    {??}{\color{red} FIX ME!}{}{}
\patchcmd{\@setcrefrange}    {??}{\color{red} FIX ME!}{}{}
\patchcmd{\@setcrefrange}    {??}{\color{red} FIX ME!}{}{}
\patchcmd{\@setcrefrange}    {??}{\color{red} FIX ME!}{}{}
\patchcmd{\@setnamecref}     {??}{\color{red} FIX ME!}{}{}
\patchcmd{\@setnamecref}     {??}{\color{red} FIX ME!}{}{}
\patchcmd{\@setcpageref}     {??}{\color{red} FIX ME!}{}{}
\patchcmd{\@setcpageref}     {??}{\color{red} FIX ME!}{}{}
\patchcmd{\@setcpagerefrange}{??}{\color{red} FIX ME!}{}{}
\patchcmd{\@setcpagerefrange}{??}{\color{red} FIX ME!}{}{}
\patchcmd{\@setcpagerefrange}{??}{\color{red} FIX ME!}{}{}
\patchcmd{\@setcpagerefrange}{??}{\color{red} FIX ME!}{}{}
\patchcmd{\@setcpagerefrange}{??}{\color{red} FIX ME!}{}{}
\patchcmd{\@cref}            {??}{\color{red} FIX ME!}{}{}
\patchcmd{\@setcite}{\bfseries ?!}{\bfseries\color{red} FIX ME!}{}{}
\patchcmd{\@citex}{\bfseries ?}{\color{red}{\@citeb}}{}{%
	\GenericWarning{}{Failed to patch \protect\@citex}}
\def\blx@citation@entry#1#2{%
	\blx@bibreq{#1}%
	\ifinlist{#1}{\blx@cites}
	{}
	{\listgadd{\blx@cites}{#1}%
		\blx@auxwrite\@mainaux{}{\string\abx@aux@cite{#1}}}%
	\ifinlistcs{#1}{blx@segm@\the\c@refsection @\the\c@refsegment}
	{}
	{\listcsgadd{blx@segm@\the\c@refsection @\the\c@refsegment}{#1}}%
	\blx@ifdata{#1}%
	{}%
	{\ifcsdef{blx@miss@\the\c@refsection}%
		{\ifinlistcs{#1}{blx@miss@\the\c@refsection}%
			{{\bfseries\color{red} cite:} }%
			{\blx@logreq@active{#2{#1}}}}%
		{\blx@logreq@active{#2{#1}}}}}
\def\blx@citeadd#1{%
	\ifcsdef{blx@keyalias@\the\c@refsection @#1}
	{\edef\blx@realkey{\csuse{blx@keyalias@\the\c@refsection @#1}}}
	{\def\blx@realkey{#1}}%
	\expandafter\blx@citation\expandafter{\blx@realkey}\blx@msg@cundefon
	\expandafter\blx@ifdata\expandafter{\blx@realkey}
	{\advance\blx@tempcnta\@ne
		\listeadd\blx@tempa{\blx@realkey}}
	{\ifnum\blx@tempcntb>\z@\multicitedelim\fi
		\expandafter\abx@missing\expandafter{\blx@realkey}%
		\advance\blx@tempcntb\@ne}}

\DeclarePairedDelimiterX{\pmodx}[1]{(}{)}{{\operator@font mod}\mkern6mu#1}
\renewcommand{\pmod}{%
  \allowbreak
  \if@display\mkern18mu\else\mkern8mu\fi
  \pmodx
}
\makeatother

\makeatletter
\newcommand{\oset}[3][0ex]{%
\raisebox{.175ex}{$%
  \mathrel{\mathop{#3}\limits^{
    \vbox to#1{\kern-2\ex@
    \hbox{$\scriptstyle#2$}\vss}}}
    $}%
    }
\makeatother

\ExplSyntaxOn
\NewDocumentEnvironment{adjunctions}{O{}}
{
	\cs_set_eq:cN {@arraycr} \farin_arraycr:
	\keys_set:nn { farin/adjunction } { #1 }
	\begin{array}
		{
			@{ \hspace { \dim_eval:n { \l_farin_left_shift_dim + \l_farin_padding_dim } } }
			r
			@{ {\farin_strut:} \l_farin_symbol_tl {} }
			l
			@{ \hspace { \dim_eval:n { \l_farin_right_shift_dim + \l_farin_padding_dim } } }
		}
	}
	{
	\end{array}
}
\keys_define:nn { farin/adjunction }
{
	leftshift       .dim_set:N = \l_farin_left_shift_dim,
	leftshift       .initial:n = 0pt,
	rightshift      .dim_set:N = \l_farin_right_shift_dim,
	rightshift      .initial:n = 0pt,
	padding         .dim_set:N = \l_farin_padding_dim,
	padding         .initial:n = 6pt,
	symbol          .tl_set:N  = \l_farin_symbol_tl,
	symbol          .initial:n = \longrightarrow,
	verticalspacing .dim_set:N  = \l_farin_vertspac_dim,
	verticalspacing .initial:n = {3pt},
}
\cs_new_protected:Npn \farin_strut:
{
	\vrule height \dim_eval:n { \ht\strutbox + 1.2\l_farin_vertspac_dim }
	depth  \dim_eval:n { \dp\strutbox + \l_farin_vertspac_dim }
	width 0pt
}
\makeatletter
\exp_args:NNo \cs_new:Npn \farin_arraycr:
{
	\@arraycr\hline
}
\makeatother
\ExplSyntaxOff

\newcommand{\myred}{BrickRed}
\hypersetup{
		colorlinks			 = true, 	
		urlcolor			 = \myred, 
		linkcolor			 = \myred, 	
		citecolor			 = PineGreen 	
}

\usepackage{xspace}

\usepackage[left = 2.54cm,right = 2.54cm,top = 2.54cm,bottom = 2.54cm]{geometry}
\usepackage{tikz}
\usetikzlibrary{matrix,fit,backgrounds,calc,arrows} 

\tikzstyle{image} = [rectangle,fill = Red!20,inner sep = -2pt]
\tikzstyle{nonzero} = [rectangle,fill = Navy!20,inner sep = 0pt]
\tikzstyle{nonzerosm} = [rectangle,fill = Navy!20,inner sep = -2pt]

\makeatletter
\newbox\xrat@below
\newbox\xrat@above
\newcommand{\xrightarrowtail}[2][]{%
  \setbox\xrat@below = \hbox{\ensuremath{\scriptstyle #1}}%
  \setbox\xrat@above = \hbox{\ensuremath{\scriptstyle #2}}%
  \pgfmathsetlengthmacro{\xrat@len}{max(\wd\xrat@below,\wd\xrat@above)+.6em}%
  \mathrel{\tikz [>->,baseline = -.55ex]
                 \draw (0,0) -- node[below = -2pt] {\box\xrat@below}
                                node[above = -2pt] {\box\xrat@above}
                       (\xrat@len,0) ;}}
\newbox\xrat@below
\newbox\xrat@above
\renewcommand{\xtwoheadrightarrow}[2][]{%
  \setbox\xrat@below = \hbox{\ensuremath{\scriptstyle #1}}%
  \setbox\xrat@above = \hbox{\ensuremath{\scriptstyle #2}}%
  \pgfmathsetlengthmacro{\xrat@len}{max(\wd\xrat@below,\wd\xrat@above)+.6em}%
  \mathrel{\tikz [->>,baseline = -.55ex]
                 \draw (0,0) -- node[below = -2pt] {\box\xrat@below}
                                node[above = -2pt] {\box\xrat@above}
                       (\xrat@len,0) ;}}
\makeatother

\newcommand{\xepi}{\xtwoheadrightarrow}
\newcommand{\epi}{\xepi{\phantom{\ \, }}}

\newcommand{\presectionskip}{-1.5\baselineskip}
\newcommand{\postsectionskip}{0.3\baselineskip}
\usepackage{titlesec}
\makeatletter
\renewcommand{\section}{\@startsection
	{chapter}{0}{0mm}
	{\presectionskip}
	{\postsectionskip}
	{\sffamily\huge}}
\renewcommand{\section}{\@startsection
	{section}{1}{0mm}
	{\presectionskip}
	{\postsectionskip}
	{\sffamily\LARGE}}
\renewcommand{\subsection}{\@startsection
	{subsection}{2}{0mm}
	{\presectionskip}
	{\postsectionskip}
	{\sffamily\Large}}
\renewcommand{\subsubsection}{\@startsection
	{subsubsection}{3}{0mm}
	{\presectionskip}
	{\postsectionskip}
	{\sffamily\normalsize}}
\renewcommand{\@seccntformat}[1]{\csname the#1\endcsname.\quad}
\newcommand\HUGE{\@setfontsize\Huge{30}{47}} 
\makeatother

\renewcommand{\SS}{\textsection}

\newcommand		{\SSS}	{Serre spectral sequence\xspace}
\newcommand		{\EMSS}	{Eilenberg--Moore spectral sequence\xspace}
\newcommand		{\exterior}	{\Lambda}
\newcommand		{\ext}		{\exterior}

\renewcommand	{\th}		{^{\mathrm{th}}}

\newcommand		{\SHC}		{\textsc{shc}\xspace}
\newcommand		{\SHCA}		{\SHC-algebra\xspace}

\newcommand		{\CGA}		{\textsc{cga}\xspace}
\newcommand		{\HGA}		{\textsc{hga}\xspace}
\newcommand		{\HGAs}		{\textsc{hga}s\xspace}

\newcommand		{\DG}		{\textsc{dg}\xspace}
\newcommand		{\DGA}		{\textsc{dga}\xspace}
\newcommand		{\DGAs}		{\textsc{dga}s\xspace}
\newcommand		{\DGC}		{\textsc{dgc}\xspace}
\newcommand		{\DGCs}		{\textsc{dgc}s\xspace}
\newcommand		{\CDGA}		{\textsc{cdga}\xspace}
\newcommand		{\CDGAs}	{\textsc{cdga}s\xspace}

\def\iter#1#2{#1^{[#2]}}

\makeatletter
\newcommand{\subalign}[1]{%
  \vcenter{%
    \Let@ \restore@math@cr \default@tag
    \baselineskip\fontdimen10 \scriptfont\tw@
    \advance\baselineskip\fontdimen12 \scriptfont\tw@
    \lineskip\thr@@\fontdimen8 \scriptfont\thr@@
    \lineskiplimit\lineskip
    \ialign{\hfil$\m@th\scriptstyle##$&$\m@th\scriptstyle{}##$\crcr
      #1\crcr
    }%
  }
}
\makeatother

\newcommand		{\bl}		{\bullet}

\newcommand		{\eqn}[1]			{\begin{align*} #1 \end{align*}}
\newcommand		{\quation}[1]		{\begin{equation} #1 \end{equation}}
\newcommand		{\case}[1]			{\begin{cases} #1 \end{cases}}

\newcommand		{\bs}				{\bigskip}

\newcommand		{\mn}				{\mspace{-2mu}}
\newcommand		{\mnn}				{\mspace{-1mu}}

\newcommand		{\nd}			{\noindent}
\newcommand		{\dsp}			{\displaystyle}
\newcommand		{\ol}			{\overline}
\newcommand		{\os}			{\overset}
\newcommand		{\us}			{\underset}

\makeatletter
	\def\smallunderbrace#1	{\mathop{\vtop{\m@th\ialign{##\crcr
		$\hfil\displaystyle{#1}\hfil$\crcr
   		\noalign{\kern3\p@\nointerlineskip}%
   		\tiny\upbracefill\crcr\noalign{\kern3\p@}}}}\limits
						}
\makeatother

\newcommand		{\ul}			{\underline}

\newcommand		{\wt}			{\widetilde}

\renewcommand	{\o}		{\circ}

\renewcommand	{\:}		{\colon}

\renewcommand	{\-}		{^{-1}}
\renewcommand	{\o}		{\circ}
\renewcommand	{\.}		{\cdot}
\newcommand		{\x}		{\times}

\DeclareMathOperator*{\otimesvariable}{%
			\mathchoice {\raisebox{.85pt}{$\displaystyle\otimes$}}
						{\raisebox{.85pt}{$\mspace{0.625mu}\otimes\mspace{0.625mu}$}}
						{\raisebox{0.7pt}{$\scriptstyle\otimes$}}
						{\raisebox{0.2pt}{$\scriptscriptstyle\otimes$}}
						}
\newcommand		{\tensor}		{\otimesvariable}
\newcommand		{\ox}			{\tensor}

\newcommand		{\Tensor}		{\bigotimes}
\newcommand		{\Direct}		{\bigoplus}

\newcommand		{\mr}			{\mathrm}

\newcommand		{\f}			{\mathfrak}
\newcommand		{\ms}			{\mathscr}

\newcommand		{\fa}		{\f a}
\newcommand		{\fb}		{\f b}

\renewcommand	{\epsilon}	{\varepsilon}
\renewcommand	{\a}		{\alpha}
\renewcommand	{\b}		{\beta}
\renewcommand	{\d}		{\delta}
\newcommand		{\h}		{\eta}

\renewcommand	{\l}		{\lambda}
\newcommand		{\vt}		{\vartheta}
\newcommand		{\vk}		{\varkappa}

\newcommand		{\e}		{\epsilon}

\newcommand		{\s}		{\sigma}

\newcommand		{\G}		{\Gamma}
\newcommand		{\D}		{\Delta}

\DeclareSymbolFont{cmletters}{OT1}{cmr}{m}{n}
\DeclareMathSymbol{\Ups}	{\mathalpha}{cmletters}{"7}
\renewcommand	{\Upsilon}	{\Ups}

\newcommand		{\F}	{\mathbb F}

\newcommand		{\Z}	{\mathbb Z}
\newcommand		{\Q}	{\mathbb Q}
\newcommand		{\R}	{\mathbb R}

\renewcommand	{\cup}		{\mspace{-1mu}\smile\mspace{-1mu}}

\DeclareRobustCommand{\lq}	{\text{\reflectbox{$/$}}}	

\newcommand		{\limit}		{\varprojlim}

\DeclareMathOperator{\diag}		{diag}

\DeclareMathOperator{\rk}		{rk }
\DeclareMathOperator{\im}		{im }

\DeclareMathOperator{\Tor}		{Tor}

\DeclareMathOperator{\Hom}		{Hom}

\newcommand		{\U}			{\mr{U}}
\newcommand		{\SU}			{\mr{SU}}

\newcommand		{\longto} 		{\longrightarrow}
\newcommand		{\lt}			{\longto}

\newcommand		{\xtoo}			{\xrightarrow} 
\newcommand		{\lmt}			{\longmapsto}

\newcommand		{\from}			{\leftarrow}
\newcommand		{\longfrom}		{\longleftarrow}

\newcommand		{\inc}			{\hookrightarrow}
\newcommand		{\xinc}			{\xhookrightarrow}
\newcommand		{\longinc}		{\xinc[]{\ \ \ \ }}

\newcommand		{\longepi}		{\xepi[]{\ \ \ \ }}

\newcommand		{\simto}		{\xrightarrow{\sim}}

\newcommand		{\longsimto}	{\os\sim\longto}
\newcommand		{\isoto}		{\longsimto}

\newcommand		{\vertsim}		{\rotatebox{90}{$\sim$}}

\newcommand		{\longhomeoto}	{\os{\smash{\homeo}}\longto}

\newcommand		{\ceq}			{\coloneqq}
\newcommand		{\eqc}			{\eqqcolon}
\newcommand		{\ideal}		{\unlhd}
\newcommand		{\idealneq}		{\lhd}

\newcommand		{\hmt}			{\simeq}
\newcommand		{\iso}			{\cong}
\newcommand		{\homeo}		{\approx}
\newcommand		{\id}			{\mr{id}}

\newcommand 		{\C}	{C^*}
\renewcommand 		{\H}	{H^*}

\newcommand 		{\HK}	{\H_K}
\newcommand 		{\HS}	{\H_S}
\newcommand			{\KGH}	{K \lq G / H}

\def\cupone{\mathbin{\cup\!\mnn_1}}
\def\cuptwo{\mathbin{\cup\!_2}}

\def\deg#1{|#1|}

\let\epsilon\varepsilon
\let\phi\varphi

\def\EEE{\mathfrak{E}}

\def\kkk{\mathfrak{k}}

\def\EE{\mathbf{E}}

\DeclareMathOperator{\Sing}{Sing}

\def\BB{\mathbf{B}}
\newcommand\B\BB
\def\OM{\boldsymbol{\Omega}}

\def\Ai{$A_{\infty}$}

\newcommand{\vpe}{\vphantom{\id_{\B A}^{\ox\,2}\big)}}
\newcommand		{\use}[2]	{\underset{\textstyle{#1\vpe}}{#2\vpe}}

\newcommand		{\clone}[1] {\use{#1}{#1}}

\usepackage{tikz-cd}
\numberwithin{equation}{section}
\usepackage{multirow}

\makeatletter
\def\eqrefhgashc#1{\textup{\tagform@{\ref*{s@#1}}}}
\makeatother

\newif\ifdebug                                                      
\debugfalse
\newcommand {\revision}[1] {\ifdebug\textcolor{Red}{#1}\else{#1}\fi}

\def\kk{k}

\newcommand{\formalitymap}{f}
\newcommand{\form}{\formalitymap}

\newcommand{\genmap}{\lambda}

\newcommand{\simpsetvar}{X}
\newcommand{\simpset}{\simpsetvar_\bl}
\newcommand{\simpgroupvar}{G}
\newcommand{\simpgrouptwovar}{H}
\newcommand{\simpgroup}{\simpgroupvar_\bl}
\newcommand{\simpgrouptwo}{\simpgrouptwovar_\bl}
\newcommand{\simpuniversalspace}{W}
\newcommand{\simpclassifyingspace}{\ol W}
\newcommand{\W}{\simpclassifyingspace}
\newcommand{\reduction}{\ol}
\newcommand{\rd}{\reduction}
\newcommand\gG{\genmap_{G}}
\newcommand\gK{\genmap_{K}}
\newcommand\gH{\genmap_{H}}
\newcommand{\onecomponent}[1]{#1{\vphantom{X^X}}^{\!\mnn+\ \mnn}{}_{\!\!\!\!\!\!(1)}}
\newcommand{\onecomponenta}[1]{#1{\vphantom{X^X}}^{+}{}_{\!\!\!\!\!\!\mn(1)}}
\newcommand\gKo{\onecomponent{(\genmap_K)}}
\newcommand\gHo{\onecomponent{(\genmap_H)}}
\renewcommand{\restriction}{\rho}
\newcommand{\rest}{\restriction}
\newcommand{\twistcohom}{t_{\mr H}}
\newcommand{\twistcochain}{t_{\mr C}}
\newcommand{\quasiisomorphism}{\Theta}
\newcommand{\quism}{\quasiisomorphism}
\newcommand{\concatenate}{\D_{\from}}

\newcommand{\Hp}{H_0}
\newcommand{\Ho}{H_1}
\newcommand{\Kp}{K_0}
\newcommand{\Ko}{K_1}
\newcommand{\Gp}{G_0}
\newcommand{\Go}{G_1}
\newcommand{\quismp}{\quism_0}
\newcommand{\quismo}{\quism_1}
\def\susp{s}
\def\desusp{s^{-1}}
\def\aa{{\vec{a}}}
\def\bb{{\vec{b}}}

\def\xx{{\vec{x}}}
\newcommand		{\tA}		{t_{\mn A}}
\newcommand		{\dA}		{d_{\mn A}}
\newcommand		{\eA}		{\e_{\mn A}}
\newcommand		{\hA}		{\h_{\mnn A}}
\newcommand		{\muA}		{\mu_{\mn A}}
\newcommand		{\dint}		{d_{\mathrm{int}}}
\newcommand		{\dext}		{d_{\mathrm{ext}}}
\newcommand		{\dbA}		{\desusp \! \bar A}
\newcommand		{\dextBAA}	{d_{\smash{\BB(A,A,A)}}^{\smash{\,\mnn\mathrm{ext}}}}
\newcommand		{\twist}	{_}
\DeclareMathOperator{\pr}{pr}
\DeclareMathOperator*{\T}{%
	\mathchoice {\raisebox{.85pt}{$\displaystyle\underline{{\otimes}}$}}
	{\mathbin{{\raisebox{.85pt}{$\underline{{\otimes}}$}}}}
	{\mathbin{{\raisebox{.7pt}{$\scriptstyle\underline{{\otimes}}$}}}}
	{\mathbin{{\raisebox{.2pt}{$\scriptscriptstyle\underline{{\otimes}}$}}}}
}
\newcommand		{\Zt}	{\mathbb{Z}[\sfrac 1{\, 2}]}

\newcommand{\Homd}{D}
\newcommand{\Homdp}{D_{\otimes}}
\newcommand{\bartimes}{*}

\begin{document}

\title
	{The cohomology of biquotients\\
	via a product on the two-sided bar construction\footnote{\ 
		The present version of this paper is the one intended for publication.
		A longer expository version, 
		with more generous exposition,
		much more historical context,
		and less abbreviated proofs,
		is available on the 
		{arXiv}
		and recommended for 
		consultation in any instances where details here 
		are found wanting~\cite{carlsonfranzlong}.
	} 
}
\author{\phantom{with an appendix by} %
		Jeffrey D. Carlson %
		\phantom{and Matthias Franz}%
		$\vphantom{X_{X_{X_{X_j}}}}\mspace{0mu}$\\
	with an appendix by Jeffrey D. Carlson
	and Matthias Franz%
		\rlap{\thanks{\ %
				M.\,F.\ was partially supported by an NSERC Discovery Grant.%
				}}
}



\maketitle

\begin{abstract}
	We compute the Borel equivariant cohomology ring
	of the left $K$-action on a homogeneous space~$G/H$, 
	where $G$ is a connected Lie group,
	$H$ and $K$ are closed, connected subgroups
	and $2$ and the torsion primes of the Lie groups 
	are units of the coefficient ring.
	As a special case, this gives the singular cohomology 
	rings of biquotients~$H \lq G / K$.
	This depends on a version of the Eilenberg--Moore~theorem 
	developed in the appendix,
	where a novel multiplicative structure on the 
	two-sided bar construction~$\BB(A',A,A'')$ 
	is defined, valid when
	$A' \from A \to A''$ is a pair of maps 
	of homotopy Gerstenhaber algebras.
\end{abstract}

Homogeneous spaces, 
which can be realized as coset spaces $G/H$ 
for $G$ a transitively acting Lie group
and $H$ the stabilizer of a point, 
are arguably the most highly symmetric, most canonical,
and most thoroughly investigated objects of study in differential geometry
after Lie groups. 
A generalization of perennial interest,
which offers many interesting examples in positive--curvature geometry, 
is the class of \emph{biquotients},
the orbit spaces $\KGH$ of $G$ under free left--right actions 
by products $K \x H$ of two closed subgroups.

The cohomology rings of a Lie group $G$ over $\Q$ 
and those finite fields $\F_p$ for which $\H(G;\Z)$ lacks $p$-torsion
have been known to be exterior algebras since the fundamental 1941 work of Hopf~\cite{hopf1941hopf},
and that of a homogeneous space $G/H$ with connected stablilizer $H$ 
has been known over $\R$ since 
work of Henri Cartan from 1950~\cite{cartan1950transgression}:
\[
	\H(G/H) 
		\iso 
	\Tor_{\H(BG)} \big(\R, \H(BH)\big)
		\mathrlap.
\]

\nd In his 1952 dissertation~\cite[\SS30]{borelthesis},
Borel used the \SSS of the Borel fibration $G \to G/H \to BH$
to show the same ring isomorphism
also holds when $H$ is of maximal rank,
with $\F_p$ coefficients if $\H(G;\Z)$ and $\H(H;\Z)$ lack $p$-torsion 
and with $\Z$ coefficients if they are torsion-free.
Beyond this, however, ring isomorphisms proved hard to find.
Starting with Baum's 1962 dissertation, 
a program began to obtain a more general result, additively at least, 
using the then-new \EMSS of the fibration $G/H \to BH \to BG$.
Cartan's result implies this spectral sequence's collapse over $\R$,
and Baum proved collapse under hypotheses covering the best-studied homogeneous spaces.
Subsequent work of authors including Baum, Larry Smith, Gugenheim, May, Munkholm, 
and Joel Wolf in the 1960s--70s proved collapse under substantially more general hypotheses%
~\cite{baum1968homogeneous,smith1967emss,gugenheimmay,munkholm1974emss,wolf1977homogeneous}
and with accordingly more involved proof techniques.
Again, however, almost all of these proofs provided only the additive structure of $\H(G/H)$,
and those which in special cases gave the ring structure 
could be factored through Borel's theory.
A collapse result guaranteeing a multiplicative isomorphism was 
provided by Franz only in 2019~\cite{franz2019homogeneous}.

Arguably the most general of the Eilenberg--Moore collapse results
actually applies to a more general case than homogeneous spaces:
Munkholm proves the collapse of the spectral sequence
corresponding to the total space of a pullback bundle
\[
	\xymatrix@C=1.25em{X \x_B E \ar[r]\ar[d]& E\ar[d]\\X \ar[r]&B}
\]
when $X$, $B$, and $E$ have polynomial cohomology on countably many generators and
the fundamental group $\pi_1(B)$ is trivial%
~\cite{munkholm1974emss}\footnote{\ 
	and additionally, in characteristic $2$,
	the cup-$1$ squares
	of polynomial generators of $\H(X)$ and $\H(E)$ vanish
	}%
.
A biquotient---and hence, in particular, a homogeneous space---%
fits into this setting when 
$X = BK$,\, $B = BG$,\, and $E = BH$ are models of classifying spaces
chosen in such a way that $BH \to BG$ is a fiber bundle,
and Singhof~\cite{singhof1993} used Munkholm's result in this case 
to determine that $\H(\KGH)$ is additively isomorphic to $\smash{\Tor^*_{\H(BG)}\big(\H(BK),\H(BH)\mnn\big)}$,
and simply to $\Tor^0$ 
in the special case where the rank of $G$ is the sum of the ranks of $K$ and $H$.
Because the \EMSS is concentrated in the $0\th$ column in this case, 
there is no extension problem, so this is in fact a ring isomorphism.
One might hope Singhof's and Franz's theorems
were instances of a more general result,
and the initial motivation of this paper was to show this hope is justified.

\begin{theorem}\label{thm:main}
	Let $G$ be a connected Lie group, 
	$H$ and $K$ closed, connected subgroups,
	and $\kk$ a principal ideal domain 
	in which $2$ and
	the torsion primes of $G$, $H$, and $K$
	are units.
	Then the Borel equivariant cohomology ring
	of the left translation action of $K$ on $G/H$,
	or, equivalently, of the two-sided action 
	of $K \x H$ on $G$ by $(x,h)\.g \ceq xgh\-$,
	is
	\[
		\HK(G/H;\kk) 
			\iso
		\H_{K \x H}(G;\kk)
			\iso
		\Tor^*_{\H(BG;\kk)} \bigl(\mspace{-.5mu}
		\H(BK;\kk), \H(BH;\kk)\mnn\big)\mathrlap.
	\]
	In particular, if the two-sided action of $K \x H$ on $G$
	is free, the cohomology ring
	of the biquotient~$\KGH$ is given by
	\[
		\H(\KGH;\kk) 
			\iso 
		\Tor^*_{\H(BG;\kk)}\bigl(\mspace{-.5mu}
		\H(BK;\kk),\H(BH;\kk)\mnn\bigr)\mathrlap.
	\]
	\end{theorem}

\nd 
\revision{See \Cref{eg:computation} for a sample computation.}
This result in fact extends 
to what could be called 
\emph{generalized biquotients}
by analogy with generalized homogeneous spaces;
see \Cref{rmk:generalized},
where we recover a known result on the cohomology of 
certain free loop spaces.

In broadest outline, our proof of \Cref{thm:main}
uses the Eilenberg--Moore theorem 
to show the cohomology of the homotopy pullback of the diagram
$BK \to BG \from BH$
is $\Tor_{\C(BG)}\bigl(\C(BK),\C(BH)\mnn\big)$
and constructs maps between $\H(B\G)$ and $\C(B\G)$ for $\G \in \{G,K,H\}$
which then induce an isomorphism with $\Tor_{\H(BG)}\bigl(\H(BK),\H(BH)\mnn\big)$.
This much it has in common with many collapse results,
but our result also shows that this map 
takes the classical product on the Tor of cochain
algebras to that on the Tor of cohomology rings.
The way this multiplicativity is 
established closely follows Franz's proof, 
but among the technical underpinnings
necessary to extend his approach
from a fibration to a general pullback 
is one substantial innovation.

A Tor of \DGAs in full generality is not endowed with a product,
and the product defined on $\Tor_{\C(BG)}\bigl(\C(BK),\C(BH)\mnn\bigr)$
is synthetic, arising from the homological external product rather 
than a multiplicative structure on the resolution itself.
As a consequence, in previous collapse results, 
one could not say whether the isomorphisms shown were multiplicative.
But there is a cochain complex $\B(A',A,A'')$,
the {two-sided bar construction},
functorial in spans~$A' \from A \to A''$ of \DGA maps,
whose cohomology is $\Tor_A(A',A'')$ under mild conditions.
In \Cref{thm:def-prod-bar},
assuming the maps in the span are actually maps 
of so-called \emph{homotopy Gerstenhaber algebras}, 
a type of \DGA with extra structure
(of which cochain algebras are the main examples),
we are able to define a 
product on $\B(A',A,A'')$
inducing a product on $\Tor_A(A',A'')$
which specializes to the known products when $A'$, $A$, $A''$
are cochain algebras or cohomology rings.
This product may be the point of greatest interest in this paper;
it is certainly the most difficult.
Then the maps between $\H(B\G)$ and $\C(B\G)$ 
alluded to in the previous paragraph, 
chosen with sufficient care,
and assuming $2$ is a unit of the coefficient ring $\kk$, 
will preserve this novel product in the transition
from $\B\bigl(\H(BK),\H(BG),\H(BH)\mnn\bigr)$ 
to $\B\bigl(\C(BK),\C (BG),\C(BH)\mnn\bigr)$.

\stepcounter{table}
\begin{layout*}
The layout of the paper is as follows. 
\end{layout*}
\begin{sublayout}
In \Cref{sec:alg} we recall algebraic conventions and the two-sided bar construction.
\end{sublayout}
\begin{sublayout}
In \Cref{sec:HGA} we discuss extended homotopy Gerstenhaber algebras and 
recall results showing normalized cochain algebras 
and cohomology rings are examples.
\end{sublayout}
\begin{sublayout}
In \Cref{sec:SHC} we recall the quasi-isomorphisms we need,
\Ai-algebra maps $\genmap$ from $\H(B\G)$ to $\C(B\G)$
for $\G \in \{K,G,H\}$ 
and \DGA maps $\formalitymap\: \C(BT) \lt \H(BT)$ 
for $T \in \{T_K,T_H\}$ maximal tori in $K$ and $H$.
These $f$ will be
defined so as to annihilate the error terms distinguishing
$\genmap$ from a genuine \DGA map,
so that the composites
$\H(BK) \to \C(BK) \to \C(BT_K) \to \H(BT_K)$ 
and~$\H(BH) \lt \H(BT_H)$ are just 
the functorially induced 
$\H\big(B(T_K \inc K)\mnn\big)$ 
and $\H\big(B(T_H \inc H)\mnn\big)$.
The \DGA quasi-isomorphisms $\formalitymap$ 
do not necessarily exist if $T$ is replaced by a more
general Lie group, and to
construct them we need a simplicial model for $BT$,
which necessitates the replacements
discussed in \Cref{sec:top}.

It is an important technical point in showing the three \Ai-maps 
$\genmap$ from $\H(B\G)$ to $\C(B\G)$ induce a map of Tors
that they are essentially functorial up to homotopy.
This fact, 
and our control over the error term annihilated by $\formalitymap\mn$,
come from the existence of a certain auxiliary structure
on extended homotopy Gerstenhaber algebras, 
a so-called \emph{strongly homotopy commutative} algebra structure $\Phi$
whose existence was proven by Franz.
The explicit formula for $\Phi$ plays a key role in \Cref{sec:prod}.
The fact that $\Phi$, so defined, is a structure of the 
type sought follows from  
an extraordinarily involved cochain-level computation~\cite{franz2019shc}.
\end{sublayout}
\begin{sublayout}
In \Cref{sec:top}, 
we apply the enhanced version of the Eilenberg--Moore
theorem proven in the appendix to the span
$BK \to BG \from BH$
whose homotopy pullback $K\lq G/H$
we are interested in,
showing $\B\big(\C(BK),\C(BG),\C(BH)\big)$
carries a weak multiplicative structure 
inducing the cup product on $\H(K\lq G /H)$.
We then show this product is preserved
when we replace $BK \to BG \from BH$
with a simplicial model
making available the \DGA formality maps 
discussed in the previous paragraph. 
\end{sublayout}
\begin{sublayout}
In \Cref{sec:final}, we weave these threads together
to prove \Cref{thm:main}. 
We first construct a quasi-isomorphism 
$\quism \: \B\big(\H(BK), \H(BG), \H(BH)\mnn\big) \lt \B\big(\C(BK), \C(BG), \C(BH)\mnn\big)$,
roughly to be thought of 
as 
${\B(\genmap_K,\genmap_G,\genmap_H)}$,
which is not necessarily multiplicative.
Letting $T_K \leq K$ and~$T_H \leq H$ be maximal tori, 
we then map the codomain into~%
${\B\big(\C(BT_K),\C(BG), \C(BT_H)\mnn\big)}$,
using the maps $\rho\: \C(B\G) \lt \C(BT_\G)$ for $\G \in \{K,H\}$
to define a strictly multiplicative map
\[\Psi = \B(\formalitymap\mn\rho,\id,\formalitymap\mn\rho)\: 
\B\big(\C(BK), \C(BG), \C(BH)\mnn\big) \lt \B\big(\H(BT_K), \C(BG), \H(BT_H)\big)
\]
inducing an injection in cohomology.
Because we have chosen $f$ and $\genmap$ to compose nicely,
the composite $\Psi\quism$ is multiplicative up to homotopy, 
and because the multiplicative map $\H(\Psi)$ is injective, 
this shows $\H(\quism)$ is multiplicative as well,
concluding the proof.
%
\end{sublayout}
\begin{sublayout}
In \Cref{sec:example}, 
we give an example computation of an equivariant cohomology ring
as Tor.
\end{sublayout}
\begin{sublayout}
Finally, in \Cref{sec:prod}, 
the technical core of the work, 
we construct a product on the two-sided bar construction
and use it to prove an Eilenberg--Moore \Cref{thm:EM-product} 
in which the product of the cohomology of the pullback of $X \to B \from E$
arises from the product on $\B\big(\C(X),\C(B),\C(E)\mnn\big)$.
\end{sublayout}

\smallskip

\nd\emph{Acknowledgments.} 
	J.D.C. would like to thank Omar Antol\'in Camarena 
	for helpful discussions of and references on matters simplicial,
	Vincent G\'elinas for fielding questions 
	of varying lucidity regarding \Ai-algebras,
	and the organizers of the Fields Institute 
	Spring 2020 Thematic Program on Toric Topology and Polyhedral Products 
	for financial and moral support.

\section{Algebras, twisting cochains, and bar constructions}\label{sec:alg}

In this background section we establish notational conventions and 
recount some foundational lemmas.
Nothing is original here save, possibly,
some of the lemmas on two-sided twisted tensor products and bar constructions,
which if not published are still likely known.
%

Fix forever a commutative base ring $\defm\kk$ with unity
over which we will consider cochain complexes, 
and with respect to which all tensor products
and $\Hom$-modules will be taken.
We will take as understood the Koszul sign convention
and the notions of 
derivation, coderivation, 
differential graded $\kk$-algebra (henceforth \textcolor{RoyalBlue}{{\scshape{dga}}}),
and differential graded $\kk$-coalgebra (\textcolor{RoyalBlue}{\DGC}).
A {commutative} \DGA is a \textcolor{RoyalBlue}{{\scshape{cdga}}}.
All \textbf{\textit{algebras}} 
we consider  
\textbf{\textit{are nonnegatively-graded, associative, and connected}} 
unless otherwise noted, and all
\textbf{\textit{coalgebras nonnegatively-graded, coassociative, and cocomplete}}.
All \DGA \textbf{\textit{maps preserve unit and augmentation}}, 
and all \DGC \textbf{\textit{maps preserve counit and coaugmentation}}.
All 
\textbf{\textit{ideals}} of \DGAs  
\textbf{\textit{will be two-sided differential}} ideals.


We also assume understood 
suspensions, desuspensions,
tensor products, 
and Hom-modules of cochain complexes,
homomorphisms of \DGAs and \DGCs,
and the classical (reduced) bar construction of a \DGA.
We nevertheless need to establish notation.

\stepcounter{table}
\begin{notmain}
	\label{def:suspension}
		Given a cochain complex $B = \Direct B_n$,
		its differential is canonically written~$\defm{d_B}$.
		For the tensor differential on a tensor product 
		$\Tensor B^{(i)}$ of complexes $B^{(i)}$ 
		we write $\defm{d_{\otimes}}$.
		We write~$\defm{\deg{x}}$ for the degree $n$ 
		of a homogeneous element $x \in B_n$. 
		The \defd{desuspension} 
		$\defm{\desusp B} = \{\defm{\desusp b} : b \in B\}$ 
		is the (re)graded $\kk$-module 
		given by $\defm{(\desusp B)_n} \ceq B_{n+1}$,
		equipped with the differential 
		$\defm{d_{\desusp B}}\: \desusp b \lmt - \desusp db$.
		The desuspension map
		$\defm\desusp\: b \lmt \desusp b$
		is a degree-$(-1)$ isomorphism of cochain complexes
		with inverse $\defm{\susp}\:\desusp B \isoto B$, 
		the \defd{suspension}.
\end{notmain}
\stepcounter{subdef}
\begin{subdef}
	\label{notation:DGA}
A \DGA $A$ is a list comprising, 
besides the underlying cochain complex $(A,d_A)$,
a 
canonically-named
	multiplication
			$\defm \mu_A\: A \ox A \lt A$,
	unit
		$\defm {\eta}_{\mn A}\: \kk \lt A$, and
	augmentation
		$\defm \e_{\mn A}\: A	\lt \kk$, 
the clarifying decorations suppressed when practicable.
The augmentation ideal is $\defm{\bar A}$
and we write $\defm{\iter\mu{n}}\: A^{\otimes n} \lt A$
for iterated multiplication.
\end{subdef}

\begin{subdef}
A \DGC comprises 
a cochain complex $(C,d_c)$,
	comultiplication
	 	$\defm \D_C\: C \lt C \ox C$,
	counit
		$\defm \e_C\: C \lt \kk$, and
	coaugmentation
		$\defm \h_C\: \kk \lt C$.
\defd{Cocompleteness} of $C$
means that every element
is annihilated by some iterate of the
reduced comultiplication 
$c \lmt \D c - 1 \ox c - c \ox 1$.
\end{subdef}
%
%
%
%
\begin{subdef}
	\label{def:iter}
	We use $\defm{\bullet}$ 
	to abbreviate indices or exponents representing 
	an indefinite number of tensor factors to be determined from context;
	for instance, if $A$ is a graded algebra,
	$\defm{a_\bl}$ denotes a pure tensor
	$a_1 \otimes \cdots \otimes a_n \in A^{\otimes n}$,
	where $n$ is to be gleaned contextually.
	A repeated $\bullet$ implies summation 
	unless explicitly stated otherwise,
	so that for example
	$d_\otimes = \id^\bl \ox d\mn_A \otimes \,\mn\id^{\bl}
	\: A^{\otimes n} \lt A^{\otimes n}$
	represents the sum of the $n$
	maps $\id^{\otimes \ell} \otimes d\mn_A \otimes \,\mn\id^{\otimes n-\ell -1}$
	applying $d\mnn_A$ to one tensor factor.
	Given a pure tensor $a_\bl \in A^{\otimes n}$,
	an expression involving 
	the string $a_\bl$ multiple times
	represents a sum 
	over order-preserving 
	subdivisions of $a_\bl$ into tensor factors.
\end{subdef}
\begin{subdef}
\label{def:permutations}
%
%
%
Our notation for tensor-factor permutations is typically
cyclic.
For example, we write
$\defm{(1 \ 2)}\: A \ox B \lt B \ox A$
for the $\kk$-linear map given by
$a \otimes b \lmt (-1)^{|a||b|}b \otimes a$.
%
We may at times also abusively use arguments to indicate permutations;
for example we may instead of $(1\ 2\ 3)$ write  
$\defm{\tau_{a\otimes b;\,c}}\: A \ox B \ox C \lt C \ox A \ox B$
for the map taking
$a \otimes b \otimes c$ to $(-1)^{|c||a| + |c||b|}\,c \otimes a \otimes b$.
Amongst all permutations, 
\emph{shuffles} will be of particular interest for us:
a $\defm{(p,q)}$\defd{-shuffle} is a permutation $\pi$ of $\{1,\ldots,p+q\}$
such that $\pi(i) < \pi(j)$ whenever $i < j$
both lie in $\{1,\ldots,p\}$ or both lie in $\{p+1,\ldots,p+q\}$.
Thus shuffles interleave two blocks, 
leaving the order within each block unaffected.
To avoid writing Koszul signs, we generally manipulate such maps directly, without evaluating on elements.
\end{subdef}

%
%
%

\begin{definition}\label{def:bar}
The graded module
underlying the 
\defd{bar construction}
$\defm{\BB A}$ 
of a \DGA $A$
is the direct sum of $\defm{\BB_n A} \ceq (\dbA)^{\otimes n}$ for $n \geq 0$.
The projection $\BB A \lt {\BB_n A}$ is denoted $\defm{\pr_n}$.
Pure tensors are written as bar-words $\defm{[a_1|\cdots|a_n]} \in \BB_n A$,
or $\defm{[a_\bl]}$ when the \defd{length} $\defm \ell\big( [a_1|\cdots|a_n]\big) = n$
is immaterial.
The empty word $\defm{[]} \in \BB_0 A$ 
is the image of 
$1 \in \kk$ under the coaugmentation;
the counit can be identified with $\pr_0$.
%
The comultiplication is the \emph{deconcatenation}
$
\defm{\D_{\BB A} }\:
[a_\bl] 
\lmt
\sum_{p = 0}^{\ell(a_\bl)}
[a_1|\cdots|a_p] 
\otimes 
[a_{p+1}|\cdots|a_{\ell(a_\bl)}]
$.
We will employ Einstein--Sweedler notation for the values of iterated deconcatenations:
$
\defm{\iter{\D_{\B A}}n} [a_\bl] 
\eqc 
\defm{[a_{(1)}]} \ox \cdots \ox  \defm{[a_{(n)}]}
$,
implicitly summing over partitions of~$[a_\bl]$
into~$n$ bar-subwords.

The sum of the tensor differentials on the $\BB_n A= (\dbA)^{\otimes n}$ 
is the \defd{internal differential} 
$\defm{\dint}$;
the \defd{external differential}
$\defm{\dext}\: \BB_n A \lt \BB_{n-1} A$
is the bar-deletion operation
$\id^\bl \ox \desusp\muA\susp^{\ox\mn 2} \ox \id^\bl$
for $n \geq 2$ and is $0$ on $\BB_{\leq 1} A$.
The total differential $\defm{d_{\B A}} \ceq \dint + \dext$
makes $\B A$ a \DGC.
This \DGC structure is functorial in \DGA maps
$f\: A \lt B$,
with the \DGC map
	$\defm {\BB\mnn f}\: \BB A \lt \BB B$
	given on $\BB_n A$ by 
	$(\desusp\mn f\mnn \susp)^{\otimes n}\:
	[a_\bl]
	\lmt [f\mnn a_1|\cdots|f\mnn a_n] \eqc \defm{[f\mnn a_\bl]}$.

Neglecting the \DGC structure,
there is a structurally irrelevant but sometimes notationally
convenient \defd{concatenation} map 
$\defm{\concatenate}\:
[a_1|\cdots|a_p] \otimes [b_1|\cdots|b_q] \lmt [a_1|\cdots|a_p|b_1|\cdots|b_q]
$
which comes from viewing
$\B A$ as the tensor algebra on $\dbA$.
\end{definition}

\begin{definition}\label{def:hom}
	Given a \DGC $C$ and \DGA $A$,
	we write $\defm{\Hom_n}(C,A)$ for 
	the $\kk$-module of $\kk$-linear maps $f$ sending each $C_j$ to $A_{j+n}$,
	and set the degree $\defm{|f|}$ to $n$ for such a map.
	We write $\defm{\Hom}(C,A)$ for the graded module 
	$\Direct_{n \in \Z} \Hom_n(C,A)$.
	This module becomes a cochain complex under the
	differential~$\defm\Homd = \defm{d_{\Hom(C,A)}}$
	given by 
	$
	\defm{\Homd\mn f} \ceq d\mn_A f - (-1)^{|f|}f d_C
	$
	and a \DGA
	with respect to the \defd{cup product} 
	$
	\phantom{,}
	\defm{f \cup g} 
	\ceq 
	\muA(f \otimes g)\D_C
	$,
	unity $\hA \e_C\: C \to \kk \to A$, 
	and augmentation  
	$f \lmt \hA f \e_C$.
	If $h \in \Hom_0(C,A)$ satisfies
	$h \eta_C = \hA$,
	then a cup-inverse
	$
	\defm{h^{\cup\,-1}}
	\ceq
	\sum_{\ell = 0}^\infty (\hA \e_C - h)^{\cup \ell}
	$
	is well-defined by cocompleteness of $C$.
%
An element 
$t \in \Hom_1(C,A)$
satisfying the three conditions
\[
\eA t = 0 = t \h_C,\qquad\qquad \Homd t = t \cup t
\]
is called a \defd{twisting cochain}.
\end{definition}

\begin{example}\label{def:tautological-twisting-cochain}
	For each connected \DGA $A$, the composition 
	\[
	\defm{\tA}\: 
	\BB A 
	\os{\smash{\pr_1}}
	\longepi
	\dbA
	\os{\susp}\lt 
	\bar A
	\longinc 
	A
	\]
	of degree $1$
	is a twisting cochain 
	called the \defd{tautological twisting cochain}.
	The assignment $t_{(-)}$
	is a natural transformation
	between the functors assigning 
	a \DGA $A$ the underlying graded modules of $\B A$ and of $A$.
	
	Any twisting cochain $C \lt A$
	extends uniquely to a \DGC
	map $C \lt \BB A$, 
	and inversely, given a \DGC map $F\: C \lt \BB A$,
	postcomposing $\tA\: \BB A \lt A$ 
	gives a twisting cochain $\tA F\: C \lt A$,
	thus prescribing a bijection
	between twisting cochains $C \lt A$ 
	and \DGC maps $C \lt \B A$.
\end{example}

\begin{example}\label{ex:shuffle-product}
	Given a \DGA $A$,
	one can check 
	the map \[
	\defm{t_\nabla} = 
	\mu_A(\tA \otimes \eA + \eA \otimes \tA)\: \BB A \ox \BB A \lt A
	\mathrlap,
	\]
	which takes $[a] \otimes 1$ and $1 \otimes [a]$ to $a$
	and annihilates all other pure tensors, is a 
	twisting cochain if and only if $A$ is a \CDGA.
	In this case, the uniquely induced 
	\DGC map 
	\[
	\defm{\mu_\nabla}\:
	\BB A \ox \BB A \lt \BB A
	\mathrlap,
	\]
	called the \defd{shuffle product},
	is a product making $\BB A$ a \DG Hopf algebra. 
	The shuffle product takes $[a_\bl] \otimes [b_\bl] \in \BB_p A \ox \BB_q A$
	to the sum of all $(p,q)$-shuffles (with Koszul sign)
	of $[a_\bl|b_\bl]$.
	
	
	More generally, let $B$ be another \DGA.
	Then the 
	\defd{shuffle map}
	\[
	\defm\nabla\: \BB A \ox \BB B \lt \BB(A \ox B)
	\]
	is the direct sum of the maps $\BB_p A \ox \BB_q B \lt \BB_{p+q}(A \ox B)$
	sending $[a_\bl] \otimes [b_\bl]$ 
	to the sum of all tensor $(p,q)$-shuffles of 
	$[a_1 \otimes 1|\cdots|a_p \otimes 1|1 \otimes b_1|\cdots|1 \otimes b_q]$.
	This is a \DGC map,
	and if $A = B$ is a \CDGA, 
	then the composition $\B\mu \o \nabla$ is the product 
	$\mu_\nabla$ of the previous paragraph.
\end{example}

The bijection between \DGC maps and twisting cochains preserves 
suitable homotopy notions.

\begin{definition}\label{def:homotopy}
A \textcolor{RoyalBlue}{{\scshape{dgc}}} \defd{homotopy}
from one \DGC map $F\: C \lt B$ to another such map $G$
is a degree-($-1$) map $H \: C \lt B$
such that 
\[
	\e_B H = 0,\qquad
	H\h_C = 0,\qquad
	\Homd H = G - F,\qquad
	\D_B H = F \ox H + H \ox G\mathrlap.
\]
We write $H\: F \,\defm{\hmt}\, G$ in this situation.
Taking $B = \BB A$ for an augmented \DGA $A$,
we may translate to an appropriate notion of homotopy for twisting cochains.
Given two twisting cochains $t,u \in \Hom_1(C,A)$,
a \defd{twisting cochain homotopy} from the former to the latter
is a map $h \in \Hom_0(C,A)$
satisfying the three conditions 
\[
	\e_C h = \eA,\qquad
	h\hA = \h_C,\qquad
	\Homd h = t \cup h - h \cup u.
\]
We again write $h\: t \,\defm{\hmt}\, u$.
\end{definition}

\begin{proposition}\label{thm:homotopy-adjunction}
	Given a \DGC homotopy $H \: C \lt \BB A$,
	the map $h = \hA \e_C + \tA H\: C \lt A$ 
	is a twisting cochain homotopy,
	and the assignment $H \lmt h$
	is a bijection from the set of \DGC homotopies $F \hmt G$
	to the set of twisting cochain homotopies $\tA F \hmt \tA G$.
\end{proposition}

Aside from providing a tractable encoding of \DGC maps $C \lt \B A$
and homotopies therebetween, 
twisting cochains $C \lt A$ 
can also be harnessed to produce new differentials
on $C \ox A$ which we will use 
to define the two-sided bar construction.

\begin{definition}\label{def:cap-product}
Let $C$ be a \DGC and $A'$ and $A''$ \DGAs,
$M$ a differential right $C$-comodule,
and $N$ a differential left $A''$-module.
One defines the \defd{cap product} with 
an element~$\phi \in \Hom(C,A'')$ by
\eqn{
	\defm{\delta^{\mr R}_{\phi}} 
		\ceq 
	(\id_{M}\otimes\mu_N)\,
	(\id_{M}\otimes \phi\otimes \id_{N})\,
	(\D_M \otimes \id_{N})\:
	 M \ox N &\lt M \ox N,\\*
	x \otimes y 
		&\lmt 
	\pm x_{(1)} \otimes \phi(x_{(2)}\mnn)y\mathrlap.
}
It is standard 
that if $\phi = t''$ is a twisting cochain,
then 
$d_{\ox} - \d^{\mr R}_{t''}$
is a differential on $M \ox N$,
making it a cochain complex we denote $\defm{M \ox_{t''} N}$
and call  a \defd{twisted tensor product}.
When $M = C$ and $N = A''$, 
this prescription makes ${C \ox_{t''} A''}$
a differential left $C$-module 
and a differential right $A''$-module.

A symmetric construction
produces a cap product endomorphism $\defm{\d^{\mr L}_\phi}$ on $P \ox Q$
from a graded linear $\phi\: C \lt A'$,
a differential \emph{left} $C$-comodule~$Q$,
and a differential \emph{right} $A'$-module~$P$,
and if $\phi = t'$ is a twisting cochain,
then $\d^{\mr L}_{t'} + d_{\ox}$ is a differential on 
$P \ox Q$; we again write $\defm{P \ox_{t'} Q}$
for the twisted tensor product.
Applying this construction to $P = A'$ and $Q =  C \ox_{t''} A''$,
we obtain a
\defd{two-sided twisted tensor product}%
~\cite[Rmks.~II.5.4]{husemollermoorestasheff1974}
\[
\defm{A' \,\ox\twist{t'}\, C \,\ox\twist{t''}\, A''} \ceq 
A' \,\ox\twist{t'}\, (C  \,\ox\twist{t''}\, A'') 
\mathrlap.
\]
This is an $(A',A'')$-bimodule, and given ideals
$\fa' \ideal A'$ and $\fa'' \ideal A''$,
we write $\defm{(\fa',\fa'')}$
for the~$(A',A'')$-sub-bimodule 
$\fa' \ox C \ox A'' + A' \ox C \ox \fa''$.
\end{definition}

%
%
%
%
%
%

%
%
%

We will need some lemmas describing when (homotopy-)commutative squares of maps
\quation{\label{eq:twisted-tensor-map}
	\begin{aligned}
	\xymatrix{
		A'_0  \ar[d]& 
		C_0 \ar[d] \ar[l]_{t'_0} \ar[r]^{t''_0}& 
		A''_0 \ar[d]\\
		A'_1& 
		C_1 \ar[l]^{t'_1}\ar[r]_{t''_1}& 
		A''_1 
	}
	\end{aligned}
}
induce a map 
$A'_0 \ox\twist{t'_0} C_0 \ox\twist{t''_0} A''_0
	\lt
A'_1 \ox\twist{t'_1} C_1 \ox\twist{t''_1} A''_1$
of twisted tensor products.
These properties seem to be well known in the one-sided case 
(see, \eg, Huebschmann~\cite[p.~360]{huebschmann1989perturbation} and
Franz \cite[\SS{7}]{franz2019homogeneous}), 
so can safely suppress 
the proofs of their easily guessed two-sided generalizations.

\begin{lemma}\label{thm:homotopy-twisted}
	Let $G\: C_0 \lt C_1$ be a \DGC map,
	$f'\: A'_0 \lt A'_1$ and $f''\: A''_0 \lt A''_1$ 
	\DGA maps,
	and~$t'_j\: C_j \lt A'_j$ and $t''_j\: C_j \lt A''_j$
	twisting cochains ($j \in \{0,1\}$).
	Then $f' \otimes G \otimes f''$ is a cochain map 
	\[
	A'_0 \,\ox\twist{t'_0}\, C_0 \,\ox\twist{t''_0}\, A''_0 
	\lt 
	A'_1 \,\ox\twist{t'_1}\, C_1 \,\ox\twist{t''_1}\, A''_1
	\]
	if $t_1' G = f' t'_0$ and $t_1'' G = f'' t''_0$.
	If $G$ is one such \DGC map and $\wt G$ another,
	and $H\: C_0 \lt C_1$ is a \DGC homotopy $G \hmt \wt G$
	such that $t'_1H$ and $t''_1 H$ are $0$,
	then $f' \ox H \ox f''$ 
	is a cochain homotopy 
	$f' \otimes G \otimes f'' \hmt f' \otimes \wt G \otimes f''$.
\end{lemma}

\begin{corollary}\label{thm:quism-DGC-map}
	Let $C_0$ and $C_1$ be \DGCs,
	$A'$ and $A''$ \DGAs,
	$t'\: C_1 \lt A'$  and $t''\: C_1 \lt A''$ twisting cochains,
	and $G\: C_0 \lt C_1$ a \DGC homomorphism.
	Then $t'\circ G\: C_0\to A'$ and $t''\circ G\: C_0\to A''$
	are twisting cochains and
	\begin{equation*}
	\id_{A'} \otimes G \otimes \id_{A''}\: 
	A' \,\ox\twist{t' G}\, 	C_0 \,\ox\twist{t'' G}\,	A''
	\lt 
	A' \,\ox\twist{t'}\, 	C_1 \,\ox\twist{t''}\, 		A''
	\end{equation*}
	a cochain map.
\end{corollary}

%
%

\begin{lemma}
	\label{thm:quism-twisting-homotopy}
	Let $C$ be a \DGC,
	$A'$ and $A''$  \DGAs,
	$t'_j\:C \lt A'$ and $t''_j\:C \lt A''$ twisting cochains 
	for~$j \in \{0,1\}$,
	and $\fa' \ideal A'$ and $\fa'' \ideal A''$ ideals,
	and suppose there exist homotopies 
	$h'\: t'_0 \hmt t'_1$ and~$h''\: t''_1 \hmt t''_0$\footnote{\ 
		Note this second homotopy goes the opposite direction to the one 
		one might expect.
		}
	with $h' \ol C \leq \fa'$ and $h'' \ol C \leq \fa''$.
	Then the composition
	\[
	(\d^{\mr L}_{h'} \otimes \id_{A''})(\id_{A'} \otimes \d^{\mr R}_{h''})\:
	A' \ox C \ox A'' \lt A' \ox C \ox A''
	\]
	is a cochain isomorphism 
	$\dsp{ A' \ox\twist{t'_0} C \ox\twist{t''_0} A'' 
		\isoto
	A' \ox\twist{t'_1} C \ox\twist{t''_1} A''}$ 
	congruent to the identity modulo $(\fa',\fa'')$. 
\end{lemma}

We now define the twisted tensor products of greatest
interest to us.

\begin{definition}\label{def:bar-two-sided}
	Let $A$, $A'$, $A''$ be augmented \DGAs
	and $F'\: \B A \lt \B A'$ and $F''\: \B A \lt \B A''$ \DGC maps.\footnote{\ 
		These are \emph{\Ai-algebra maps} from $A$ to $A'$ and $A''$ 
		in the language to be introduced in \Cref{sec:SHC}.
	}	
	Associated to $F'$ and $F''$ are twisting cochains $t' = t_{A'} F'$
	and $t'' = t_{A''} F''$,
	and we define the \defd{two-sided bar construction}
	as the two-sided twisted tensor product
	\[
		\defm{\BB(A',A,A'')} 
		\ceq 
			A' \ox\twist{t_{A'}F'} \BB A \ox\twist{t_{A''}F''} A''\mathrlap{.}
	\]
	We will write a pure tensor 
	$a' \otimes [a_\bl] \otimes a'' \in \B(A',A, A'')$
	as $\defm{a'[a_\bl]a''}$.
	The differential of $\BB(A',A,A'')$
	preserves the length filtration 
	$\defm{\BB_{\leq \ell}(A',A,A'')} \ceq A' \ox \BB_{\leq \ell} A \ox A''$.%
	\label{def:length-filtration}
\end{definition}

\begin{example}\label{def:bar-two-sided-DGA}
	This two-sided bar construction,
	when we take $F' = \B f'$ and $F'' = \B f''$ 
	for \DGA maps $f'\: A \lt A'$ and $f''\: A \lt A''$,
	coincides with the familiar classical two-sided bar construction 
	of differential
	graded algebras.
	In this case, the twisting cochains $t'$ and $t''$ defining
	the twisted tensor product 
	$\B(A',A,A'') = A' \ox\twist{t'} \B A \ox\twist{t''} A''$
	are simply $t_{A'}\B f' = f' \tA$
	and $f'' \tA$, respectively.
\end{example}

We see already from the preceding lemmas that \DGA maps from $A'$ and $A''$
and \DGC maps from $\B A$ induce maps of two-sided bar constructions,
and want to extend this technique to induce maps from commutative diagrams of \DGC maps of the form

	\begin{equation}\label{eq:bar-map}
	\begin{aligned}
	\xymatrix@R=2.5em{
		\BB A'_0 \ar[d]_{G'}& 
		\BB A_0 \ar[d]|(.45)\hole|G|(.525)\hole 
				\ar[l]_{F'_0} 
				\ar[r]^{F''_0}& 
		\BB A''_0 \ar[d]^{G''}\\
		\BB A'_1&
		\BB A_1	\ar[l]^{F'_1} 
				\ar[r]_{F''_1}&
		\BB A''_1
				\mathrlap.
	}
	\end{aligned}
	\end{equation}
We write $\defm{G'_{(1)}} \: \ol{A'_0} \lt \ol{A'_1}$ for the composite
$t_{A'_1} \o G' \o \desusp_{A'_0}$ and similarly
for  
$\defm{G''_{(1)}}$; 
later, in \Cref{def:A-infty-map}, 
we will see $\h\e + G'_{(1)}\: A'_0 \lt A'_1$ 
and $\h\e + G''_{(1)}\: A''_0 \lt A''_1$ 
can be understood as a sort of approximate \DGA maps.

\begin{proposition}[Wolf~{\cite[Thm.~7]{wolf1977homogeneous}};
	\cf~Gugenheim--Munkholm~{~\cite[Thm.~3.5$_*$]{gugenheimmunkholm1974}}%
	]\label{def:Gamma-generalized}%
	Suppose given a strictly commuting diagram of \DGC maps as in \eqref{eq:bar-map}.
	Then, in the notation of \Cref{def:bar},
	we may define a cochain map 
	 $\B(A'_0 , A_0 , A''_0)\lt \B(A'_1 , A_1 , A''_1)$
	by
\quation{\label{eq:Gamma-generalized}
	\defm{\B(G',G,G'')} 
		= 
	\big(
		\Upsilon'
			\, 	\ox \, 
		G
			\, 	\ox \, 
		\Upsilon''		
	\big)
	(
		\id_{A'_0} 
			\,\ox\, 
		\iter{\D_{\BB A_0}}3 
			\,\ox\, 
		\id_{A''_0}
	)\mathrlap,
}
	where
\[
	\begin{aligned}
	\Upsilon'|_{\kk \ox \B A_0} 
		&\,=\,
	\id_{\kk} \otimes \e_{\B A_0}\: & 
	c'[a_{(1)}] &\lmt 
	c'\e[a_{(1)}] 
	\mathrlap,\\
	\Upsilon'|_{\ol{A'_0} \ox \B A_0} 
		&\,=\,
	t_{A'_1}G'\concatenate(\desusp_{A'_0} \ox F'_0)\:\ &
		a'[a_{(1)}] &\lmt 
		\pm t_{A'_1}G'\big([a'] \ox F'_0[a_{(1)}]\big)
	\mathrlap,\\
	\Upsilon''|_{\B A_0 \ox \kk}
		&\,=\,
	\e_{\B A_0} \otimes \id_{\kk}
	\:&
	[a_{(3)}] c'' &\lmt 
	\e[a_{(3)}]c'' 
	\mathrlap,\\
	\Upsilon''|_{\B A_0 \ox \ol{A''_0}} 
		&\,=\,
	t_{A''_1}G''\concatenate(F''_0 \otimes \desusp_{A''_0})\:\ &
	[a_{(3)}]a'' &\lmt 
	\pm t_{A''_1}G''\big(F''_0[a_{(3)}] \otimes [a'']\big)
		\mathrlap.
	\end{aligned}
\]

%
	If $t_{A'_1}G'$ takes $\B_{\geq 2} A'_0$ 
	into an ideal $\f a' \ideal A'_1$
	and $t_{A''_1}G''$ takes $\B_{\geq 2} A''_0$ into $\f a'' \ideal A''_1$,
	then 
	\[
	\phantom{ \quad \big(\mathrm{mod}\ (\f a',\f a'')\big)}
	\B(G',g,G'')
		\equiv 
	(\h\e + G'_{(1)}) \ox G \ox {(G''_{(1)} + \h\e)} 
	\quad \big(\mathrm{mod}\ (\f a',\f a'')\big).%
	\footnote{\ 
	Wolf~\cite{wolf1977homogeneous} cites 
	his unpublished dissertation~\cite{wolfthesis} 
		for a proof of this result;
		however, this work was not publicly accessible at writing.
		(It is now available at the author's webpage.)
		For our proof, which turns out to be nontrivial, 
		see the \texttt{arXiv} 
		version of the present paper~\cite[Prop.~1.26]{carlsonfranzlong}.
		Proofs of the lemmas 
		and many of the standard claims in this section
		can be found there as well.
	} %
	 \]
\end{proposition}

\begin{notation}\label{def:bar-map-dga}
	In the event that
	the maps $G$, $G'$, $G''$
	in \eqref{eq:bar-map},
	are respectively
	$\B g$,~$\B g'$,~$\B g''$
	for 
	$g\: A_0 \lt A_1$,
	$g'\: A'_0 \lt A'_1$,
	$g''\: A''_0 \lt A''_1$,
	we abuse notation by writing
	the map of \Cref{def:Gamma-generalized} as
	$\defm{\B(g',g,g'')} \ceq \B(G',G,G'')$.
	Such a map is also 
	a special case of that of \Cref{thm:homotopy-twisted}.
\end{notation}

We are interested in 
two-sided bar constructions because they provide 
functorial resolutions.


\begin{proposition}[See Barthel--May--Riehl {\cite[after Prop.~10.19]{barthelmayriehl2014}}]%
		\label{thm:EM-A-infty}\label{thm:bar-computes-Tor}
		Suppose $\kk$ is a principal ideal domain
		and $A \lt A'$ a map of \DGAs 
		flat over $\kk$.
		Then 
		given another \DGA map $A \lt A''$,
		the cohomology of the classical two-sided bar construction
		$\smash{\B(A',A,A'') = 
		\B(A',A,A)\ox_A A''}$
		of \Cref{def:bar-two-sided-DGA}
		is $\Tor_{A}(A' ,A'')$.
		If $\H(A)$ and $\H(A'')$ are flat over $\kk$,
		then the $E_2$ page of the associated algebraic \EMSS
		is $\Tor_{\H(A)}\big(\H(A'),\H(A'')\big)$.
\end{proposition}

\begin{discussion}
For the proof of \Cref{thm:two-sided-ring-map} later on, 
we will need a more explicit expression for the differential on the classical two-sided bar construction
$\B(A',A,A'')$ of 
\Cref{def:bar-two-sided-DGA}.	
Note that the total differential is the sum of the
tensor differential $d_{\ox}$ on $A' \ox \BB A \ox A''$
and two cap products.
Separating the summand $\id_{A'} \otimes d_{\B A} \otimes \id_{A''}$
into its internal and external components,
$d_{\B(A',A,A'')}$ decomposes as the sum of tensor differentials
and an ``external'' differential which is the sum of the two cap products
and the bar-deletion differential.
Explicitly, 
appealing to \Cref{def:cap-product},
this external differential 
$
	\defm{d_{\BB(A',A,A'')}^{\mathrm{ext}}} 
$
is given on the summand
$\BB_\ell(A',A,A'')$ by
\begin{equation}
\begin{aligned}\label{eq:def:dext-two-sided}
&\ \phantom{{}+{}} {}\mu_{A'}(\id_{A'} \otimes f'\susp) 
&\otimes\ & \id_{\dbA}^{\otimes \ell - 1} 
&\otimes\ &  \id_{A''} \\
&\ {}+{\id_{A'}} 
&\otimes\ & \dext
&\otimes\ & \id_{A''}	\\
&\ {}-{\id_{A'}} 
&\otimes\ & \id_{\dbA}^{\otimes \ell - 1} 
&\otimes\ & \mu_{A''}(f''\susp \otimes \id_{A''})
\mathrlap.
\end{aligned}
\end{equation}
	The minus sign in the last line of \eqref{eq:def:dext-two-sided} is important
	and not consistently noted in the literature.
\end{discussion}

\section{Homotopy Gerstenhaber algebras}\label{sec:HGA}
We ultimately want to describe
the multiplicative structure on
 $\Tor_{\C(B)}\big(\C(X),\C(E)\mnn\big)$
in terms of a product on the two-sided bar construction.
The special property
of a normalized singular cochain algebra $A$
enabling us to do so
will turn out to be the existence
of a \DG Hopf algebra structure on $\B A$,
or in other words a \DGC map 
$\BB A \ox \BB A \lt \BB A$
making it a \DGA.\footnote{\
	Under our blanket assumption that $A$ is graded and connected,
	a bialgebra structure on $\B A$ admits a unique antipode
	defined by a formula analogous to the one defining the cup-inverse
	in \Cref{def:hom}~\cite[Props.~1.4.14,24]{grinbergreiner}.
	}

\begin{definition}\label{def:HGA}
Let $A$ be a \DGA such that $\B A$ admits a multiplication
making it a \DG Hopf algebra.
If the twisting cochain 
$\defm \EE \ceq \tA \mu_{\BB A}\: \B A \ox \B A \lt A$
satisfies $\defm{\EE_{j,\ell}} \ceq \EE|_{\BB_j A \,\mnn\ox\,\mnn \BB_\ell A} = 0$
for $j \geq 2$, 
we call $A$, equipped with $\mu_{\B A}$ (equivalently, with $\EE$),
a \defd{homotopy Gerstenhaber algebra}
(\textcolor{RoyalBlue}{\HGA}).	
\end{definition}

\begin{notation}\label{notation:E}
	We will on occasion write the values of the product $\mu_{\B A}$ 
	in infix notation as
	$[a_\bl] \,\defm{\bartimes}\, [b_\bl]$.
	It will be useful later to translate~$\EE_{1,\bl}$ 
	into a degree-zero 
	map on 
	$A \ox \BB A$
	by taking
	\eqn{
		\defm{\EEE}\ceq \EE_{1,\bl}
		(\desusp \otimes \id_{\BB A})\: A\ox \BB A &\lt A,\\*
		{a[b_\bl]} 
		&\lmt \EE\big([a] \otimes [b_\bl]\big)
		\mathrlap.
	}
\nd	We will also have occasion to use the operations 
	$\defm{E_\ell} = \EE_{1,\ell} \o (\desusp)^{\otimes 1 + \ell}\: 
	A \ox A^{\otimes \ell} \lt A$,
	whose values we write as
	$\defm{E_\ell(a;b_\bl)} = 
	E_\ell(a \otimes b_1 \otimes \cdots \otimes b_\ell)$.
	 An \HGA structure on $A$  
	is more commonly phrased as a list of conditions%
	~\cite[(6.2)--(6.4)]{franz2019homogeneous}
	on the $E_\ell$.
\end{notation}

\begin{remark}\label{def:vanishing-components}
We note that unitality of $\mu_{\BB A}$
implies $\EE_{0,1}$ and $\EE_{1,0}$ must 
both be $\susp\: \dbA \isoto \bar A$
and~$\EE_{j,0}$ and~$\EE_{0,\ell}$ must be $0$
for $j, \ell \geq 2$.
The ``$t\h = 0$'' clause in 
the definition of a twisting cochain in \Cref{def:hom}
implies $\EE_{0,0} = \EE \h_{\BB A} = 0$.
For notational convenience,
we will extend $\EE_{0,1}$ and $\EE_{1,0}$ to both be $\susp\: \desusp A \lt A$
(thus respectively sending $[1] \otimes []$ and $[] \otimes [1]$ to $1$)
and extend the operations $\EE_{1,\ell}$ to $(\desusp A)^{\otimes 1 + \ell}$
for $\ell \geq 1$
by setting them to $0$ on pure tensors any of whose factors 
lies in $\desusp \im \hA$.
That is, 
we will sometimes consider bar-words containing a letter $c$ in 
the coefficient ring $\kk$, 
but such words 
are to be annihilated by $\EE_{j,\ell}$ unless $(j,\ell) = (1,0)$ or $(0,1)$.
\end{remark}

\begin{example}\label{ex:HGA-shuffle}
	If $A$ is a commutative \DGA,
	then the shuffle product $\mu_\nabla$ of \Cref{ex:shuffle-product} 
	makes the bar construction~$\BB A$ a \DG Hopf algebra, 
	so that $A$ becomes an \HGA.
	The corresponding twisting cochain $t_\nabla = \EE$
	satisfies~$\EE = \EE_{0,1} + \EE_{1,0}$,
	so the operations $E_{\ell}$ vanish for $\ell \geq 1$.
\end{example}

\medskip

The Maurer--Cartan identity 
$\Homd\EE = \EE \cup \EE$ 
for the twisting cochain
$\EE\: \BB A \ox \BB A \lt A$ defining an $\HGA$ structure on $A$ 
yields the following useful formulae upon mild rearrangement:
\quation{\label{eq:EEE-Cartan}
\begin{aligned}
	\EEE(\muA \ox \id_{\BB A}) 
		&= 
	\sum_\pi \muA(\EEE \ox \EEE)\pi 
	\:
	\bar A \ox \bar A \ox \BB A \lt A\mathrlap,
	\\
	\EEE\big(a_1a_2[b_\bl]\big) &= \sum \pm \EEE\big(a_1[b_\bl]\big)\EEE\big(a_2[b_\bl]\big)
	\mathrlap,
\end{aligned}
}
where the sum is over shuffles
$\pi\: [a_1|a_2] \otimes [b_\bl] \lmt \pm [a_1|b_\bl] \otimes [a_2|b_\bl]$, and
\begin{equation}\label{eq:d-prime-EEE}
\begin{aligned}
\phantom{.}
	\Homdp\EEE 
	&	= 
	\dA \EEE - \EEE d_{\ox}
	\muA(s \otimes \EEE)(1 \ 2)
	{}\,+\,
	{\EEE(\id_A \otimes \dextBAA)} 
{}	\,-\,
	{\muA(\EEE \otimes s)}
	\mathrlap,\\
	(\Homdp\EEE)\big(a[b_\bl]\big) 
		&=
	\pm b_1 \EEE\big(a[b_2|b_\bl|b_\ell]\big)
	\pm \phantom{b_1}\EEE\big(a[b_1|b_\bl|b_p b_{p+1}|b_\bl|b_\ell]\big)
	\pm \phantom{b_1}\EEE\big(a[b_1|b_\bl|b_{\ell - 1}]\big)b_\ell
\end{aligned}
\end{equation}
for $\ell \geq 2$, with the middle term omitted if $\ell = 1$.\footnote{\ 
	Compare Franz~{\cite[Lem.~3.1]{franz2020szczarba}}
	for conditions for a map $C \lt \OM C \otimes \OM C$
	to be a twisting cochain inducing a \DGC structure
	on the cobar construction $\OM C$ of a \DGC $C$.
} 
Here $d_{\ox}$ refers to 
the ``internal'' differential  $d_{\bar A} \otimes \id_{\BB A} + \id_{\bar A} \otimes d_{\ox}$ on $\bar A \otimes \BB A$,
omitting the external differential on the $\BB A$ factor,
and $\defm{\Homdp}$ 
is the derivation on $\Hom(\bar A \ox \BB A,A)$
defined with respect to this internal differential and $\dA$.

Critically, there is a natural way to make 
a variant of the cochain algebra
into an \HGA.

\begin{definition}\label{def:cochains}
	The \defd{normalized cochain algebra} $\defm{\C(\simpset;\kk)}$
	on a simplicial set $\simpset$
	is the \DG subalgebra 
	containing all and only cochains 
	vanishing on each degenerate simplex.
	It is augmented with respect to the map
	$\C(\simpset;\kk) \epi \C(x_0;\kk) \simto \kk$
	induced by the inclusion of the simplicial subset associated to 
	any chosen basepoint $x_0 \in X_0$.
\end{definition}

The so-called \emph{interval-cut operations},
among which classical cup product and Steenrod 
cup-$i$ products are the most famous examples, 
are known~\cite[Thm.~2.15]{mccluresmith2003} 
to define 
on $\C(\simpset;\kk)$
the action of a symmetric \DG-operad~$\defm{\ms X\!\mnn}$
called the 
\defd{sequence operad}~\cite{mccluresmith2003} or
\emph{surjection operad}~\cite{bergerfresse2004operad}.
This operad is known
to be a quotient of the \DG-operad $\defm{\ms E}$ associated
to the classical Barratt--Eccles simplicial operad~\cite[Thm.~1.3.2]{bergerfresse2004operad},
an $E_\infty$-operad 
filtered by an increasing sequence
of suboperads $F_n\ms E\mnn$
whose geometric realizations
are equivalent to the little $n$-cubes operads,
and the sequence operad is accordingly filtered 
by quotients $\defm{F_n \ms X}$~\cite[Lem.~1.6.1]{bergerfresse2004operad}.
As it turns out,
an \HGA structure on a \DGA $A$
is a precisely an algebra structure over 
$F_2 \ms X\!\mnn$~{%
		\cite[Thm.~4.1]{mccluresmith2003}%
		\cite[\SS1.6.6]{bergerfresse2004operad}
	(\cf~Franz~{\cite[(3.13)]{franz2020szczarba}}
	for the signs of the operations).
Franz~{\cite[\SS3.2]{franz2019shc}}
defines an 
\defd{extended homotopy Gerstenhaber algebra} structure on a \DGA $A$
to be an algebra structure over a
certain suboperad ${F'_3 \ms X\!\mnn}$ 
of $F_3 \ms X\!\mnn$
containing $F_2 \ms X\!\mnn$.
This structure comes in the form of operations
$\defm{F_{p,q}}\: A^{\otimes p} \otimes A^{\otimes q} \lt 
	A$ for $p,q \geq 1$, in addition to the $E_\ell$,
	satisfying several compatibility conditions.
An \defd{extended}
		\textcolor{RoyalBlue}{\HGA} 
		\defd{homomorphism} 
	$f\: A \lt B$ 
	is a \DGA map
	which is simultaneously an~$F'_3 \ms X\!\mnn$-algebra homomorphism
(which is to say it distributes over the $E_\ell$
and the~$F_{p,q}$).
Thus we have the following.

\begin{corollary}\label{thm:cochain-xHGA}
		For any pointed simplicial set $X_\bl$, 
		its algebra $\C(X_\bl)$
		of normalized cochains is naturally an extended \HGA.
		Any \CDGA $A$ is naturally an extended \HGA.
\end{corollary}

\begin{convention}\label{def:xHGA-CDGA}
{\textit{\textbf{All} \CDGAs \textbf{in this work, 
				particularly cohomology rings, 
				come equipped with this trivial extended} 
				\HGA \textbf{structure.}}}
Consequently, if $A$ is an extended \HGA and $B$ is a \CDGA, 
then an extended \HGA map $f\: A \lt B$ 
annihilates the values of $F_{p,q}$
and of $E_\ell$ for $\ell \geq 1$.
\end{convention}

\section{Strong homotopy commutativity}\label{sec:SHC}

As adverted to in the introduction, we will require maps
between \DGAs that are not multiplicative,
but still preserve the multiplication up to coherent homotopy.

\begin{definition}\label{def:A-infty-map}
Given two $A_\infty$-algebras $A$ and $B$,
an \defd{$A_\infty$-map} from the former to the latter
is defined to be a \DGC homomorphism $F\: \BB A \lt \BB B$.
Such a map $F$ is determined by 
the sequence of compositions 
\[
\defm{F_{(n)}} \: A^{\otimes n} \xtoo[\sim]{(\desusp)^{\otimes n}} 
\BB_n A \longinc \BB A \os F\lt \BB B \os{t_B}\lt B
\mathrlap.
\] 
%
If $A$ and $B$ are \DGAs,
then we write
$\smash{\defm{\onecomponenta{F}}} = F_{(1)} + \h_B\e_{\mn A}\: A \lt B$.
This is a cochain map, and
if $\smash{\H(\onecomponenta{F})}\: \H(A) \lt \H(B)$ is an isomorphism,
we call $F$ an \defd{\Ai-quasi-isomorphism}.\footnote{\ 
	The map $\onecomponenta{F}$
	can be seen as a \DGA map up to homotopy,
	and an ideal $\fb \ideal B$ containing the image 
	of $t_B F(\B_{\geq 2} A)$ can be seen as a measure of its deviation 
	from multiplicativity,
	for $\onecomponenta{F}$ is a \DGA map just if $t_B F(\B_{\geq 2} A) = 0$; 
	and if so, then $F = \BB\onecomponenta{F}$.
}
\end{definition}

\nd We will see it is possible to construct an \Ai-quasi-isomorphism
from the \CDGA $\H(A)$ with trivial differential 
to $A$ itself when $A$ is a \DGA with polynomial cohomology.
To do so, we need two auxiliary notions.

\begin{proposition}[{\cite[Prop.~3.3]{munkholm1974emss}}]%
	\label{def:internal-tensor}
	There exists 
	an \defd{internal tensor product}
	of \Ai-maps,
	which, given \DGAs $A$, $A'$, $B$, $B'$
	and \DGC maps $F\: \B A \lt \B A'$
	and $G\: \B B \lt \B B'$, produces a 
	\DGC map 
	\[
	\defm{F \ {\ul\ox}\  G}\: \B(A \ox B) \lt \B(A' \ox B')
	\mathrlap.
	\] 
	The operation $(F,G) \lmt F \T G$
	\begin{itemize}
		\item 
		extends the tensor product
		of \DGA maps in the sense
		that if $f\: A \lt A'$ and $g\: B \lt B'$
		are \DGA maps,
		then $\B\mn f \T \B g = \B(f \otimes g)$,
		\item is functorial in each variable separately,
		\item is unital in the sense
		that $F \T \,\mn\id_{\kk}$ and $\id_{\kk} \T F$
		can be identified with $F$, and
		\item is associative in the sense that
		given a third pair $C$, $C'$ of \DGAs 
		and a \DGC map $H\: \B C \lt \B C'$,
		the iterated products $(F \T G) \T H$
		and 
		$F \T {(G \T H)}\: \B(A \ox B \ox C) \lt \B(A' \ox B' \ox C')$
		agree.
		
	\end{itemize}

\end{proposition}

This notion of internal tensor product is connected to the ordinary external one 
by the shuffle map $\nabla$ of \Cref{ex:shuffle-product}.

\begin{proposition}[{\cite[Lem.~4.4]{franz2019homogeneous}}]%
\label{prop:shuffle-tensor}
	Let \DGAs $A'_0$, $A'_1$, $A''_0$, $A''_1$ 
	and \DGC maps $G'\:\B A'_0 \lt \B A'_1$ and $G''\:\B A''_0 \lt \B A''_1$ be given.
	Then the following square commutes:
	\[
	\xymatrix@C=1.25em@R=4.5em{
		\B A'_0 \ox \B A''_0	\ar[r]^\nabla		\ar[d]_{G' \ox\, G''}	&	
		\B (A'_0 \ox A''_0)							\ar[d]^{G' \T\, G''}	\\
		\B A'_1 \ox \B A''_1	\ar[r]_\nabla							&	
		\B (A'_1 \ox A''_1)												\mathrlap.
}
	\]
\end{proposition}

The internal tensor product allows us to define another notion of homotopy commutativity
\emph{a~priori} unrelated to \HGAs.

\begin{definition}\label{def:SHC}
A \defd{strongly homotopy commutative algebra} 
(henceforth \textcolor{RoyalBlue}{\SHC}\defd{-algebra})
is an augmented \DGA $A$
equipped with an $A_\infty$-map from 
$A \ox A$ to $A$ 
(\ie, a \DGC map $\Phi\: \BB(A \ox A) \lt \BB A$),
satisfying the following conditions:

\begin{enumerate}
	\item\label{def:Phi-mu}
	Its $1$-component $\Phi_{(1)}\: \bar A \ox \bar A \lt \bar A$
		is the restriction $\muA|_{\bar A \ox \bar A}$ of the given product on $A$.
	\item It is strictly unital
	in the sense that
		$\Phi \o \BB(\id_A \T \hA) = \id_{\BB A} = \Phi \o \BB(\hA \T {\,\mn\id_A})$.
	\item It is homotopy-associative:
		there is a homotopy 
		between
		$\Phi(\Phi\,\mn \T {\,\mn\id_A})\: \BB(A \ox A \ox A) \lt \BB A$ 
		and $\Phi(\id_A \T \Phi)$.
	\item It is homotopy-commutative:
		there is a homotopy 
		between
		$\Phi
		$
		and~$\Phi \o \BB(1\ 2)$, 
		where $(1\ 2)$
		is the tensor-factor transposition
		${A \ox A} \isoto {A \ox A}$
		of \Cref{def:permutations}.	
\end{enumerate}
We define the iterates
$\defm{\iter\Phi{n}}\: \BB(A^{\otimes n}) \lt \BB A$
of the structure map $\Phi$
by $\iter\Phi{2} \ceq \Phi$ and $\iter\Phi{n+1} \ceq \Phi(\iter\Phi{n} \T {\mspace{1mu}\id_A})$.
\end{definition}

An \SHC-algebra structure $\Phi$ on a \DGA $A$ 
allows one to combine maps
in a 
useful way:
given sequences $(A_j)_{j = 1}^n$ of \DGAs 
and $(F_j\: \B A_j \lt \B A)_{j=1}^n$ of \DGC maps,
the composite
\quation{\label{eq:compilation}
	\B\big(\mn\bigotimes A_j\big) 
	\xtoo{\ul\bigotimes F_j} 
	\B (A^{\otimes n})
	\xtoo{\iter \Phi n}
	\B A
}
is guaranteed to be another \DGC map.
We will say this map is \defd{compiled} from the $F_j$.

Associated to each homogeneous element $a \in A$ 
is a map from the free \DGA on one generator of degree $|a|$
taking this generator to $a$.
If $a$ is a cocycle of even degree,
this map factors through the map
$\defm{\genmap_a}\: \kk[x] \lt A$ taking $x$ to $a$,
where the differential on $\kk[x]$ is trivial.
Thus, given a list $\defm{\vec a} = (a_j)$ of even-degree cocycles of $A$
and taking $A_j \ceq \kk[x_j]$ for $|x_j| = |a_j|$
and $F_j \ceq \B\genmap_{a_j}$, the compilation procedure \eqref{eq:compilation}
yields a \DGC map
\begin{equation}\label{eq:genmap}
	\defm{\genmap_{\aa}}
		\: \,
	\B\big(\mn\bigotimes_{j=1}^n k[x_j]\,\mn\big)
		\xrightarrow{
			\B(
				\raisebox{-.5pt}{$\otimes$}
				\,\genmap_{a_j}
				)
			}
	\B(A^{\otimes n}) 
		\xtoo{\iter{\Phi}{n}} 
	\B A
		\mathrlap.
\end{equation}
Then $\smash{\onecomponent{(\genmap_\aa)}}$ is easily seen to be the optimist's
candidate for a ring homomorphism,
\eqn{
	\defm{\kk[\vec x]} \ceq
		\kk[x_1,\ldots,x_n]
		&\lt 
	A
		\mathrlap,\\
	x_1^{p_1}\cdots x_n^{p_n} 
		&\lmt 
	a_1^{p_1}\cdots a_n^{p_n}
		\mathrlap.
} 
Though this map is in fact almost never multiplicative, 
it is at least a quasi-isomorphism.

\begin{prop}[Stasheff--Halperin~{\cite[Thm.~9]{halperinstasheff1970}%
								\cite[7.2]{munkholm1974emss}}]%
								\label{thm:genmap-quism}
	If $A$ is an \SHCA whose cohomology ring $\H(A)$ is polynomial
	on classes represented by even-degree elements $a_j \in A$
	and $\genmap_{\aa}$ is defined as in \eqref{eq:genmap},
	then $\smash{\H \onecomponent{(\genmap_\aa)}}\: 
	\kk[\xx] \lt \H(A)$ 
	is an isomorphism.
\end{prop}

In the target application, 
$A$ is the normalized cochain algebra
on a classifying space $BG$.
Munkholm~\cite[Prop.~4.7]{munkholm1974emss}\footnote{\ 
	and stating a bit less, 
	Gugenheim and Munkholm \cite[Prop.~4.2]{gugenheimmunkholm1974}%
	}
showed that the cochain algebra of a simplicial set
admits a natural \SHCA structure,
but to define the product on the
two-sided bar construction we will need later for our variant~\ref{thm:EM-product} 
of the Eilenberg--Moore theorem,
we will use a result of Franz
defining an \SHCA structure in terms of extended \HGA operations.

\begin{theorem}[Franz {\cite[Thm.~1.1, (4.2)]{franz2019shc}}]
	\label{thm:xHGA-SHC}
	An extended \HGA $A$ admits an \SHC-algebra structure
	whose structure map~$\Phi$
	and associativity and commutativity homotopies
	are defined in terms of the extended \HGA operations on $A$
	and hence are natural in extended \HGA maps.
	Moreover, the composite~%
	$\Phi \o \nabla 
		\: 
	\BB A \ox \BB A 
		\to 
	\BB(A \ox A) 
		\to 
	\BB A
	$
	of the shuffle map of \Cref{ex:shuffle-product} 
	and this structure map
	is the given product $\mu_{\B A}$ making $A$ an \HGA.
\end{theorem}

We will require an explicit formula for $\tA \Phi_{\mn A}$ 
in the proof of \Cref{thm:formula-prod-bar}, 
but hold off on stating it until then.
From \Cref{thm:cochain-xHGA} and \Cref{thm:xHGA-SHC},
we immediately have the following.

\begin{corollary}\label{thm:SHC-cochain}
	The normalized cochain algebra ${\C(\simpset;\kk)}$ 
	of a pointed simplicial set
	is naturally an \SHC-algebra.
\end{corollary}

\begin{corollary}\label{thm:genmap-quism-cochain}
	Given a pointed simplicial set $\simpset$ with polynomial cohomology
	and any list $\vec a$ in $\C(\simpset;\kk)$ 
	of representatives for $\kk$-algebra generators of $\H(\simpset;\kk)$,
	the \DGC map $\genmap_{\aa}\: \B\H(X_\bl;\kk) \lt \B\C(\simpset;\kk)$
	given in \eqref{eq:genmap}
	is
	such that $\H\onecomponent{(\genmap_\aa)}$ 
	can be identified with the identity map of 
	${\H(\simpset;\kk)}$.
\end{corollary}

We have observed that though it is 
a quasi-isomorphism in the cases that interest us,
the extended 
$1$-component $\smash{\onecomponent{(\genmap_\aa)}}\: 
	\H(\simpset) \lt \C(\simpset)$
is rarely a \DGA homomorphism on the nose.
However, the extended \HGA structure 
and the resulting \SHC structure on $\C(\simpset)$
guarantee that, loosely speaking, 
it is a \DGA map
\emph{up to an error term} contained in an
ideal $\kkk_{\simpset}$ of $\C(\simpset)$ functorial in $\simpset$,
independently of the choice of representatives 
$\aa$ in $\C(\simpset)$;
thinking of 
$\kkk_{\simpset}$ as a neighborhood of $0$, 
we may consider it a sort of 
uniform bound on failure to be a \DGA map.
It will be an important point in the proof of our main result
that the bounding ideal $\kkk_{\simpset}$ 
lies in the kernel
of the formality map $\formalitymap$ to be described in \Cref{thm:formality-map},
and hence $\formalitymap$ annihilates the error term.

\begin{definition}[Franz {\cite[(10.2)]{franz2019homogeneous}}]\label{def:kkk}
Given a pointed simplicial set $\simpset$,
viewing its normalized cochain algebra $\C(\simpset)$
as an 
\revision{extended \HGA via \Cref{thm:cochain-xHGA}}
we denote by 
$
\defm\kkk = \defm{\kkk_{\simpset}} \idealneq \C(\simpset) 
$
the ideal generated by the following elements,
where $a$,~$b$,~$b_{\bullet}$,~$c_{\bullet}$
range over $\C(\simpset)$:
\begin{enumarabic}[itemsep=-3pt]
	\item 
		coboundaries,
	\item
	\label{def-ax-1}
		elements of odd degree,
	\item
		elements of the form~$E_{\ell}(a; b_\bl)$ for $\ell = \ell(b_\bl) \geq 1$,
	\item
	\label{def-ax-4}
		elements of the form~$F_{p,q}(b_\bl;c_\bl)$ with~$(p,q)\ne(1,1)$,
	\item
	\label{def-ax-5}
		elements of the form~$a\cuptwo E_{\ell}(b;c_\bl)$ with~$\ell \geq 2$,
	\item
	\label{def-ax-6}
		elements of the form~$
			a\cuptwo 
			\smash{
			\Big(\mn
				\cdots
				\mnn\big((
					b_0 \cupone b_1) \cupone b_2
				\big) 
				\cupone  \cdots 
			\mnn\Big)
			}
		$
	for cocycles $a$ and $b_\bl$.
\end{enumarabic}
\end{definition}

From the naturality of the {extended \HGA} 
structure on a cochain algebra,
the following functoriality property of $\kkk$ is immediate.

\begin{lemma}[Franz~{\cite[Prop.~10.1]{franz2019homogeneous}}]\label{thm:kkk-functorial}
	If $\phi\: Y_\bl \lt \simpset$ is a map of pointed simplicial sets, 
	then the ideals of \Cref{def:kkk} 
	satisfy $\phi^{*} \kkk_{\simpset} \leq \kkk_{Y_\bl}$.
\end{lemma}

Morally speaking, then,
our maps $\genmap_\aa$ are functorial and multiplicative modulo $\kkk$:

\begin{theorem}[Franz~{\cite[Prop.~7.2]{franz2019shc}%
					   \cite[Prop.~11.5, %
							Thm.~11.6]{franz2019homogeneous}}]%
\label{thm:genmap-natural}
	Suppose $2$ is a unit of~$\kk$.
	Then the maps
	$\genmap_{\aa}\colon \kk[\xx] \lt \C(\simpset)$ 
	of \eqref{eq:genmap}
	are functorial modulo $\kkk$ in the sense that
	given a map~$\phi\: Y_\bl \lt \simpset$ of simplicial sets with polynomial
	cohomology and sequences $\aa$ in $\C(\simpset)$ and~$\bb$ in $\C(Y_\bl)$ 
	representing 
	generators of $\H(\simpset)$ 
	and~$\H(Y_\bl)$ respectively,
	the left diagram of \eqref{eq:genmap-homotopy}
	commutes up to a \DGC homotopy 
	$\defm{H_\phi}\: \B\H(\simpset) \lt \B\C(Y_\bl)$
	whose associated twisting cochain homotopy
	$\h\e + t_{\C(Y_\bl)} H_\phi$ sends $\B_{\geq 1}\H(X_\bl)$ into $\kkk_{Y_\bl}$.
	\quation{\label{eq:genmap-homotopy}
	\begin{aligned}
	\xymatrix@R=4.75em@C=3em{
		\BB\H(\simpset){\vphantom{\big(\kk[\vec x]\big)}} \ar
		[d]_{\genmap_{\aa}} \ar[r]^{\BB\H\phi} &
		\BB	\H(Y_\bl) \ar
		[d]^{\genmap_{\bb}} \\
		\BB	\C(\simpset) \ar[r]_{\BB\C\phi} & 
		\BB	\C(Y_\bl){\vphantom{\big(\big)}}
		}
			\qquad
			\quad
			\qquad
	\xymatrix@C=1em@R=4.75em{
		\mathllap\BB\big(\kk[\xx] \otimes \kk[\xx]\mspace{1mu}\big){\vphantom{\H(X_\bl)}} 
		\ar[r]^(.55){\BB\mu{\vphantom{\H}}}
		\ar[d]_{\genmap_{\aa} \ \T\, \genmap_{\aa}}&
		\BB\big(\kk[\xx]\mspace{1mu}\big) \ar[d]^{\genmap_{\aa}}\\
		\mathllap\BB\big(\mnn\C(\simpset) \ox \C(\simpset)\big) \ar[r]_(.645){\Phi{\vphantom{\C}}}&
		\BB\C(\simpset)
		}			
	\end{aligned}
	}
	Moreover $\genmap_\aa$ is multiplicative modulo $\kkk$
	in the sense that $t_{\C(\simpset)} \genmap_\aa$ 
	takes $\BB_{\geq 2}\big(\kk[\xx]\mspace{1mu}\big)$ into $\kkk_{\simpset}$
	and the right diagram of \eqref{eq:genmap-homotopy}
	commutes up to a \DGC homotopy 
	$\defm{H_\mu}\: \B\big(\kk[\xx]^{\ox\mn 2}\big) \lt \B\C(\simpset)$
	whose associated twisting cochain homotopy $\h\e + t_{\C(\simpset)} H_\mu$ sends 
	$\B_{\geq 1}\big(\kk[\xx]^{\ox\mn 2}\big)$ into $\kkk_{\simpset}$.
\end{theorem}

We remark that the mere existence of such homotopies, 
without the bound in terms of $\kkk$, 
lies at the heart of Munkholm's collapse proof~\cite{munkholm1974emss}.
To make precise our claim that there is a map
$\C(BT) \lt \H(BT)$ 
in the other direction annihilating $\kkk$,
we recall the classical construction
of the simplicial classifying space:

\begin{proposition}[See, \eg, May {\cite[\SS21]{maysimplicial}}]\label{def:W}
	Let $\simpgroup$ be a simplicial group.
	Then there exists a contractible simplicial $\simpgroup$-space 
	$\defm \simpuniversalspace \simpgroup$
	functorial in $\simpgroup$
	That is, if $\phi\:\simpgroup \lt \simpgrouptwo$
	is a homomorphism of simplicial groups,
	then $\simpuniversalspace\mn \phi\: 
	\simpuniversalspace \simpgroup \lt \simpuniversalspace \simpgrouptwo$
	is $\phi$-equivariant in the sense that
	$(\simpuniversalspace\mn\phi)(x\.g) = (\simpuniversalspace\mn\phi)(x)\.\phi(g)$
	for $x \in (\simpuniversalspace\simpgroup)_n$
	and $g \in \simpgroupvar_n$.
	The projection
	$\simpuniversalspace\simpgroup \lt 
	\defm{\simpclassifyingspace} \simpgroup \ceq \simpuniversalspace \simpgroup / 	\simpgroup$ 
	is a principal $\simpgroup$-bundle,
	and the base $\simpclassifyingspace\simpgroup$ 
	is the classifying simplicial set
	for simplicial principal $\simpgroup$-bundles.
\end{proposition}

The promised map is then provided by the following result.

\begin{theorem}[Franz~{\cite[Thm.~9.6, Prop.~9.7]{franz2019homogeneous}}]%
	\label{thm:formality-map}
	Let $T_\bl$ be a simplicial abelian group,
	pointed at $1 \in T_0$ and
	such that the cohomology $\H(T_\bl;\Z)$ 
	of the normalized cochain complex 
	$\C(T_\bl)$ 
	is an exterior algebra on finitely many
	degree-$1$ generators.
	Then there exists a quasi-isomorphism
	$\defm\formalitymap = \formalitymap_T\: 
	\C(\simpclassifyingspace  T_\bl) \lt \H(\simpclassifyingspace  T_\bl)$
	of extended \HGAs,
	called the \defd{formality map},
	which annihilates all extended \HGA operations~$F_{p,q}$ except
	for~$F_{1,1} = -{\cuptwo}$.
	If $2$ is a unit of~$\kk$,
	the formality map
	can be chosen so as to annihilate all $\cuptwo$-products of cocycles
	and such that its kernel contains 
	the ideal $\kkk_{\simpclassifyingspace T_\bl}$ 
	of \Cref{def:kkk}.
\end{theorem}

\section{Simplicial substitution}%
\label{sec:simplicial}%
\label{sec:BT-formal}%
\label{sec:top}

In this section we introduce a simplicial form 
of the homotopy pullback whose
cohomology we will ultimately compute.

\begin{discussion}\label{disc:top-simplicial-transition}
We use the functorial Milgram model~%
\cite{milgram1967bar,steenrod1968milgram,segal1968classifying}
of the universal principal bundle
$EG \to BG$
associated to a Lie group $G$.
It is well known
that $\KGH$ is the homotopy pullback
of the diagram $BK \to BG \from BH$
when $K \x H$ acts freely on $G$~\cite[\SS2]{singhof1993}.
More generally,
writing $\defm i\: K \lt G$
and $\defm j\: H \lt G$
for the inclusions, 
the homotopy pullback
can be identified as $EK \ox_K G/H$,\footnote{\ 
	\revision{This is the quotient of
		$EK \x G/H$ under the equivalence relation
		$(ek,gH) \sim (e,kgH)$ for $k \in K$,
		more commonly denoted $EK \x_K G/H$,
		the same notation as a pullback.}
} which is the pullback
of the diagram $\smash{BK \us{Bi}\to BG \from EG/H}$,
the homeomorphism
being given by 
$e \otimes gH \lmt \big(eK,(Ei)(e)gH\big)$.

It is even better known that $Ej\: EH \lt EG$
is $H$-equivariant and induces a homotopy equivalence
$BH \lt EG/H$.
The resulting map of triples induces
a map of bar constructions
$\B\big(\C(BK),\C(BG),\C(EG/H)\big) \lt 
\B\big(\C(BK),\C(BG),\C(BH)\big)$
which the algebraic \EMSS shows is a quasi-isomorphism
by \Cref{thm:EM-A-infty}.
Moreover, by \Cref{thm:EM-product},
this map is multiplicative with respect to the product
of \Cref{thm:def-prod-bar}
and hence induces a ring quasi-isomorphism.

If we write $\defm{\rd G}$
for the singular complex $\Sing G$ of a connected Lie group,
made a simplicial group by equipping each level $(\Sing G)_n$
with the valuewise multiplication of maps $\D^n \to G$,
then 
the counit $|\rd G| \lt G$
of the standard adjunction 
$|{-}| \,\dashv\,\, {\Sing}$
between the geometric realization
and singular complex functors is a homomorphism
and a homotopy equivalence,
natural in continuous homomorphisms~{\cite[Thm.~16.6(ii)]{maysimplicial}}.
For any simplicial group $\simpgroup$,
there is a 
natural $|\simpgroup|$-equivariant homeomorphism
$|\simpuniversalspace \simpgroup| \longhomeoto E|\simpgroup|$
descending to a natural homeomorphism
$|\simpclassifyingspace \simpgroup| \longhomeoto B|\simpgroup|$%
~\cite{bergerhuebschmann1998barW},
so for $\simpgroup = \rd G$
we get a composite weak homotopy equivalence 
$|\simpuniversalspace \rd G|
\to
E|\rd G|
\to
EG$
equivariant with respect to the homomorphism $|\rd G| \lt G$.
Long exact homotopy sequences then show
the map
$\simpclassifyingspace\rd G = 
|\simpuniversalspace\rd G|/|\rd G| \to (EG)/G = BG$
is a natural weak homotopy equivalence.
By naturality,
these maps yield a 
multiplicative quasi-isomorphism
$\B\big(\C(BK),\C(BG),\C(BH)\big)
\lt \B\big(\C|\simpclassifyingspace\rd K|,
			\C|\simpclassifyingspace\rd G|,
			\C|\simpclassifyingspace\rd H|\big)$.

Finally, the unit
$X_\bl \lmt \Sing |X_\bl|$ of the standard adjunction
is a natural weak equ\revision{iv}alence of simplicial sets whenever $X_\bl$ is a Kan complex~\cite[Thm.~16.6(i)]{maysimplicial},
which we know $\W\,\mn\rd K$, $\W\,\mn\rd G$,
and $\W\,\mn\rd H$ are~\cite[Lem.~21.3]{maysimplicial},
so we have a multiplicative quasi-isomorphism
to a more economical bar construction, 
$\B\big(\C|\simpclassifyingspace\rd K|,
			\C|\simpclassifyingspace\rd G|,
			\C|\simpclassifyingspace\rd H|\big)
\lt
\B\big(\C(\simpclassifyingspace\rd K),
\C(\simpclassifyingspace\rd G),
\C(\simpclassifyingspace\rd H)\big)
$, although our proof does not strictly speaking require this simplification.

\end{discussion}

All told, we have the following.

\begin{proposition}\label{thm:iso-bar-KGH}
	Let $G$ be a compact, connected Lie group and $H$ and $K$ closed subgroups,
	and suppose the coefficient ring $\kk$ is 
	a principal ideal domain.
	Then there is an isomorphism of graded algebras
	\begin{equation*}
		\H_K(G/ H) 
			\iso 
		\H\BB\big(\C(\W\,\mn\rd K),\C(\W\,\mn\rd G),\C(\W\,\mn\rd H)\mnn\big),
	\end{equation*}
	natural with respect to the diagram~$K \inc G \hookleftarrow H$.
\end{proposition}
\begin{proof}
	It remains only to check naturality,
	but this
	follows from the functoriality of $W$, $\Sing$, $|{-}|$, and $\C$
	and the naturality of the unit, the counit, 
	and the Eilenberg--Moore isomorphism.
\end{proof}

\begin{remark}\label{rmk:generalized}
May and Neumann observed~\cite{mayneumanncohomologygeneralized}
that the general result~\cite[Thm.~4.3]{gugenheimmay} 
behind the Gugenheim--May computation
of the cohomology groups 
of a homogeneous space $G/H$ 
applies equally in many cases 
to a \emph{generalized homogenous space},
meaning the homotopy fiber $F$
of a map $B_H \lt B_G$ of path-connected spaces.
Here $F$ is to be thought of as $G/H$ for 
$G = \Omega B_G$ and $H = \Omega B_H$,
and we assume
a trivial $\pi_1(B_G)$-action on $\H(F)$
and that $\H(B_G)$ be a polynomial
algebra over the principal ideal domain $\kk$,
finitely generated in each degree.
The key additional assumption
is that there be a map
$BT \lt B_H$,
where $BT$ is a $K(\Z^n,2)$,
making $\H(B_H)$ a free module of finite rank over $\H(BT)$;
then $T = \Omega BT$ is called a \emph{maximal torus}
of $H$.

Our result similarly generalizes. 
Given three path-connected spaces $B_G$, $B_K$, $B_H$
such that $K$ and $H$ admit maximal tori 
$T_K$ and $T_H$ 
and $\H(B_G)$, $\H(B_K)$, and $\H(B_H)$
are polynomial algebras of finite type, 
then replacing $\Sing BT_K$ and $\Sing BT_H$
by $\W\,\mn\rd T_K$ and $\W\,\mn\rd T_K$ respectively,
the same proof as in the following \Cref{sec:final}
shows that the cohomology 
of the homotopy pullback $Y$ of 
the span $B_K \to B_G \from B_H$
is given as a graded group by $\Tor_{\H(B_G)}\mn\big(\H(B_K),\H(B_H)\big)$,
and as a ring if $2$ is a unit of $\kk$.

We are not aware of great interest in such spaces $Y$ in general,
but some special cases have been studied.
For instance, if $B_K$ and $B_H$ are points, 
we recover that $\H(\Omega B_G)$ is an exterior algebra,
and for the span $B_K \to B_K \x B_K \from B_K$,
with both maps the diagonal, 
we find the cohomology ring of the free loop space
$ L B_K$ is the tensor product of $\H(B_K)$ and $\H(K)$.
This is an instance of a more general result of
Saneblidze~\cite{saneblidze2009bitwisted}, 
derived using Hochschild homology
and a model of the total space of the fibration $\Omega Y \to LY \to Y$
equipped with a cochain-level product 
related to the product described in \Cref{thm:formula-prod-bar}
and the \HGA structure on $\C(Y)$.
His result, however, requires of $\H(Y)$
only that the cup-one squares 
of a set of polynomial generators 
be zero, without any maximal torus assumption,
and of $\kk$ only that it be a commutative ring.
\end{remark}

\section{The quasi-isomorphisms}%
\label{sec:final}

We have finally assembled the necessary ingredients to prove \Cref{thm:main}.

\begin{notation}\label{def:proof-notation}
In the calculation that follows,
the base ring $\kk$ \textit{\textbf{is now a principal ideal domain}},
still usually suppressed in the notation.
We will not require that $2$ be a unit initially.
Let $\defm G$ be a connected Lie group 
and $\defm H$ and $\defm K$ closed, connected subgroups
such that $BG$, $BK$, and~$BH$ have polynomial cohomology over~$\kk$
(equivalently, such that the torsion primes of $G,H,K$ are units).
We will work with normalized cochains on 
the simplicial models $\W\,\mn\rd G$ of \Cref{sec:top}, 
but in the notation for cohomology identify $\H(BG)$ with $\H(\W\,\mn\rd G)$
and so on, suppressing the natural isomorphisms 
induced by the simplicial weak equivalences 
$\W\,\mn\rd G  \to \Sing |\W\,\mn\rd G| \to \Sing BG$.
\end{notation}

From \Cref{thm:iso-bar-KGH} we have a natural isomorphism 
 \begin{equation}\label{eq:map-Tor-CGK}
 \HK(G/H) \iso \Tor_{\C(\W\,\mn\rd G)}\bigl(\C(\W\,\mn\rd K),\C(\W\,\mn\rd H)\mnn\bigr),
 \end{equation} 
of graded algebras, and our goal is to use the maps $\genmap$ of \eqref{eq:genmap}
and $\formalitymap$ of \Cref{thm:formality-map}
to induce a \CGA isomorphism 
\begin{equation}\label{eq:map-Tor-HGK}
	\HK(G/H) \iso \Tor_{\H(BG)}\bigl(\H(BK),\H(BH)\mnn\bigr)
\end{equation}
natural in inclusion diagrams~$K \inc G \hookleftarrow H$.
On the level of graded modules, 
this isomorphism 
particularly means the \EMSS 
of the homotopy pullback square of $BK \to BG \from BH$ collapses,
a theorem of Munkholm~\cite{munkholm1974emss}
we will reprove as \Cref{thm:Theta-additive}.
To improve this result to a ring isomorphism,
as both Tors arise as the cohomology of a two-sided bar construction,
we will connect the argument \DGAs of the two 
via maps preserving enough structure
to guarantee our novel product on the two-sided bar construction
is preserved up to homotopy.
The structure of this section closely follows that 
of Section 12 in 
Franz's paper~\cite{franz2019homogeneous}
and specializes to it the case $K = 1$.

\begin{discussion}
We begin by constructing an additive 
cochain map between the two bar constructions. 
Selecting arbitrarily and fixing a list $\defm{\vec a}$ of cocycles 
representing irredundant generators of~$\H(BG)$,
and similarly for $\H(BK)$ and~$\H(BH)$,
we may use \eqref{eq:genmap} 
to construct \Ai-quasi-isomorphisms
$\defm{\genmap_G}\: \B\H(BG) \lt \B\C(\W\,\mn\rd G)$,
$\defm{\genmap_K}$, 
and~$\defm{\genmap_H}$. 
Recall that these maps are selected 
so that if $x_j$ denotes the cohomology class 
of $a_j \in \C(\W\,\mn\rd G)$ and $\rk G = n$,
then the 
extension~$\smash{\onecomponent{(\genmap_G)}}
= \h\e + (\genmap_G)_{(1)}$ 
of the $1$-component
is the additive quasi-isomorphism
\eqn{
	\mathllap{\kk[\xx] \iso {}}
		\H(BG) 
	&\lt 
		\C(\W\,\mn\rd G)
	\mathrlap,\\
		x_1^{p_1}\cdots x_n^{p_n} 
	&\lmt 
		a_1^{p_1}\cdots a_n^{p_n}
}
and similarly for $\genmap_K$ and $\genmap_H$. 
We write the canonical twisting cochains $\BB\H(BG) \lt \H(BG)$ 
and $\BB\C(\W\,\mn\rd G) \lt \C(\W\,\mn\rd G)$ respectively as $\defm{\twistcohom}$ and $\defm{\twistcochain}$,
and for any homomorphism $L' \lt L$ 
of topological groups will write 
$\defm{\restriction} = \defm{\restriction_{L'}^L}$ 
for the functorially induced \DGA maps
$\C(\W\,\mn\rd L) \lt \C(\W\,\mn\rd L')$ and
$\H(BL) \lt \H(BL')$.
Then the candidate quasi-isomorphism
\begin{equation}\label{eq:quism}
\defm{\quasiisomorphism}\: 
	\BB\big(\H(BK),\H(BG),\H(BH)\big) \lt \BB\big(\C(\W\,\mn\rd K),\C(\W\,\mn\rd G),\C(\W\,\mn\rd H)\big)
\end{equation}
is defined as the composition of the cochain maps
\[
\xymatrix@R=3em@C=6.75em{
	\H(BK)			\ \otimesunder{\,\twistcohom\BB\restriction\,} \
		\BB \H(BG)	\ \otimesunder{\,\twistcohom\BB\restriction\,} \ 
		\H(BH) 
		\ar@{}[d]^(.1){}="a"^(.45){\BB(\genmap_K,\id,\genmap_H)} 
		\ar "a";[d] 
&
				\ 
				\C(\W\,\mn\rd K)	\ \ \otimesunder{\,\twistcohom\BB\restriction\genmap_G\,} \ \ 
				\BB \H(BG)  \ \ \otimesunder{\,\twistcohom\BB\restriction\genmap_G\,} \ \ 
				\C(\W\,\mn\rd H)
				\ 
						\ar@{}[d]^(.1){}="a"_(.45){\B(\id,\genmap_G,\id)} 
						\ar "a";[d] 							
		\\
		\ 
	\C(\W\,\mn\rd K)	\ \ \otimesunder{\,\twistcohom\genmap_K\BB\restriction\,} \ \ 
		\BB \H(BG)	\ \ \otimesunder{\,\twistcohom\genmap_H\BB\restriction\,} \ \ 
		\C(\W\,\mn\rd H) 
		\ 
				\ar@{}[ru]^(.15){}="a"^(.8){}="b" 
				\ar"a";"b"		
				^	
				{\revision{\theta}}
	&
	\C(\W\,\mn\rd K) 		 \ \otimesunder{\,\twistcochain\BB\restriction\,} \
		\BB \C(\W\,\mn\rd G) \ \otimesunder{\,\twistcochain\BB\restriction\,} \ 
		\C(\W\,\mn\rd H) 
}
\]

\smallskip

\nd given respectively, from beginning to end,
by 
Lemmas~\ref{def:Gamma-generalized},~\ref{thm:quism-twisting-homotopy}, and~\ref{thm:quism-DGC-map},
where 
\revision{$
	\defm\theta = 
	(\delta^{\mr L}_{h^{{K}}} \otimes\, \id) \o 
	(\id\, \otimes\, \delta^{\mr R}_{h^{{H}}})
$},
the twisting cochain homotopies~%
$\defm{h^K}\: \BB \H(BG) \lt \C(\W\,\mn\rd K)$ and 
$\defm{h^H}\: \BB \H(BG) \lt \C(\W\,\mn\rd H)$ 
com\revision{ing} from \Cref{thm:genmap-natural}.
\end{discussion}

\begin{lemma}\label{thm:Theta-modulo-k}
	The map $\quism$ defined in \eqref{eq:quism} satisfies
	\begin{equation*}
	\phantom{
		\quad
		\Big(\mn\mathrm{mod}\ \big(C^{\geq 1}(BK),C^{\geq 1}(BH)\big)\mn\Big)
	} 
	\quasiisomorphism 
		\equiv 
	\gKo \otimes \gG \otimes \gHo
		\quad 
		\big(\mr{mod}\ (\kkk_{\W\,\mn\rd K}, \kkk_{\W\,\mn\rd H})\big)
	\mathrlap.
	\end{equation*}
\end{lemma}

\begin{proof}
	As $\smash{t_{\C(\W\,\mn\rd K)}\gK}$ 
	and $\smash{t_{\C(\W\,\mn\rd H)}\gH}$ 
	respectively take 
	$\BB_{\geq 2}\H(BK)$ and~$\BB_{\geq 2}\H(BH)$ 
	into~$\kkk_{\W\,\mn\rd K}$ and~$\kkk_{\W\,\mn\rd H}$
	by \Cref{thm:genmap-natural},
	the first factor $\BB(\genmap_K,\id,\genmap_H)$  
	is congruent to~$\smash{\gKo \otimes \id \otimes \gHo}$
	modulo~$(\kkk_{\W\,\mn\rd K},\kkk_{\W\,\mn\rd H})$
	by \Cref{def:Gamma-generalized}.
	Next, by \Cref{thm:genmap-natural} again, 
	the twisting cochain homotopies~$h^K$ and~$h^H$ 
	respectively send $\BB_{\geq 1}\H(BG)$ into 
	$\kkk_{\W\,\mn\rd K}$ and $\kkk_{\W\,\mn\rd H}$,
	so 
	$(\d_{h^K}^{\mr L} \otimes \id)(\id \otimes \d_{h^H}^{\mr R})$
	is congruent to the identity modulo $(\kkk_{\W\,\mn\rd K},\kkk_{\W\,\mn\rd H})$
	by \Cref{thm:quism-twisting-homotopy}.
	Finally, $\B(\id,\genmap_G,\id)$ is $\id \otimes \genmap_G \otimes \id$.
\end{proof}

\begin{proposition}\label{thm:Theta-additive}
	We retain the notations $G$, $K$, $H$ from \Cref{def:proof-notation} and 
	 $\Theta$ from \eqref{eq:quism}.
	\begin{enumroman}
		\item
		\label{thm:Theta-additive-1}
		The induced map 
		\begin{equation*}
			\H(\quasiisomorphism)\: \Tor^*_{\H(BG)}\big(\H(BK),\H(BH)\mnn\big) 
				\lt 
			\Tor^*_{\C(\W\,\mn\rd G)}\bigl(\C(\W\,\mn\rd K),\C(\W\,\mn\rd H)\mnn\bigr)
		\end{equation*}
		of graded $\kk$-modules is an isomorphism.
		\item
		\label{thm:Theta-additive-2}
		The Eilenberg--Moore spectral sequence for 
		the homotopy pullback of $BK \to BG \from BH$
		collapses at the $E_2$ page.
	\end{enumroman}
\end{proposition}

\begin{proof} (i)
	Because $\H(BG)$ and $\H(BK)$ are flat over $\kk$,
	under the length filtration of two-sided bar constructions
	discussed in \Cref{def:length-filtration},
	the $E_2$ page of the target 
	under the map of associated filtration spectral sequences 
	induced by $\quism$ 
	is again $\Tor^*_{\H(BG)}\big(\H(BK),\H(BH)\big)$
	as noted in \Cref{thm:bar-computes-Tor}.
	Since $\gG$, $\smash{\gKo}$, and $\smash{\gHo}$
	are each quasi-isomorphisms,
	by \Cref{thm:Theta-modulo-k},
	the map of $E_2$ pages is the identity map.
	These are half-plane spectral sequences
	with exiting differentials and in the
	associated exact couples $(E_1,A_1)$
	one has $A_1^{p} = 0$ for $p > 0$,
	and	particularly $\limit\mn{}_p\, A_1^{p} = 0$,
	so they are strongly convergent~\cite[Thm.~6.1(a)]{boardman1999conditionally}
	and hence $\H(\quism)$ is a graded $\kk$-linear isomorphism~\cite[Thm.~5.3]{boardman1999conditionally}.

	\bigskip
	
	\noindent	(ii)
	In the map of spectral sequences, the codomain is the 
	Eilenberg--Moore spectral sequence of the homotopy pullback
	of the diagram~$BK \from BG \to BH$.
	We have seen the spectral sequence map is 
	a $\kk$-linear graded isomorphism from $E_2$ on, 
	so it is enough to show the domain spectral sequence collapses,
	but $E_2$ of this domain sequence is already  
	isomorphic to the sequence's target~%
	$\H\BB\big(\H(BK),\H(BG),\H(BH)\big)$
	as a graded $\kk$-module.
\end{proof}

To show $\H(\Theta)$ is multiplicative and natural
will involve the formality map of \Cref{sec:BT-formal},
and particularly 
\textit{\textbf{from here on out, we will need}} 
$2$ \textit{\textbf{to be a unit in}} $\kk$.
We will soon specialize to maximal tori, but for now
let arbitrary compact tori $\defm{T_K}$ and $\defm{T_H}$
and simplicial group homomorphisms 
$\defm{\a_K}\: \rd T_K \lt \rd K$ and 
$\defm{\a_H}\: \rd T_H \lt \rd H$ be given,
and choose formality maps~$\defm{\formalitymap_{T_K}}\: \C(\W\,\mn\rd T_K)\to \H(BT_K)$ 
and~$\defm{\formalitymap_{T_H}}\: \C(\W\,\mn\rd T_H)\to \H(BT_H)$ 
as guaranteed by \Cref{thm:formality-map},
recalling that these maps respectively annihilate
the ideals
$\kkk_{\W\,\mn\rd T_K}$ and $\kkk_{\W\,\mn\rd T_H}$
defined in \Cref{def:kkk}.
As $\formalitymap$ and $\restriction$ are \HGA 
and hence \DGA maps,
\Cref{thm:homotopy-twisted} provides a cochain map
of two-sided bar constructions
\quation{\label{eq:def-Psi}
\defm{\Psi}_{\a_K,\a_H}
\: 
	\C(\W\,\mn\rd K) \ox_{\restriction t_C} \BB \C(\W\,\mn\rd G) \ox_{\restriction t_C} \C(\W\,\mn\rd H)
		\xrightarrow{	\form\mn
						\rest
						\ox \id \ox 
						\mn\form\mn
						\rest
					}
	\H(BT_K) 
		\ox_{\form\mn\rest t_C} 
	\BB \C(\W\,\mn\rd G) 
		\ox_{\form\mn\rest t_C} 
	\H(BT_H).\ \ 
}

\smallskip


\begin{lemma}\label{thm:HPsi}
	Let $G$, $K$, $H$ be as in \Cref{def:proof-notation},
	$\Theta$ as in \eqref{eq:quism}, 
	and $\Psi$ as in \eqref{eq:def-Psi}.
	\begin{enumroman}
		\item
		\label{thm:HPsi-1}
			The cochain map~$\Psi$ 
			is multiplicative with respect to the product $\wt\mu$  of \Cref{thm:def-prod-bar}.
		\item
		\label{thm:HPsi-2}
			The composite $\Psi\quism$ is equal to
			\[
\smash{
				\restriction^K_{T_K} \otimes \genmap_G \otimes \restriction^H_{T_H}\:	
				\H(BK) \ox_{\restriction\twistcohom} \BB \H(BG) \ox_{\restriction \twistcohom} \H(BH) 
				\lt
				\H(BT_K) 
				\ox_{\form\mn\rest \twistcochain}
				\BB \C(\W\,\mn\rd G) 
				\ox_{\form\mn\rest \twistcochain} 
				\H(BT_H) \mathrlap{.}
}
			\]
			
			
		\item
		\label{thm:HPsi-3}
			If 
			$\a_K$ and $\a_H$ are 
			inclusions of maximal tori, 
			the induced map~$\H(\Psi)$ in cohomology is injective.
	\end{enumroman}
\end{lemma}


\begin{proof}
	(i) Multiplicativity of $\Psi$ 
	follows from naturality of $\wt\mu$
	since $\formalitymap$ and $\restriction$ are \HGA maps.
	
	\medskip
	
\nd (ii) 
The restriction $\restriction_{T_K}^K$ sends $\kkk_{\W\,\mn\rd K}$ to 
$\kkk_{\W\,\mn\rd T_K}$ by \Cref{thm:kkk-functorial},
while $\formalitymap_{K}$ annihilates $\kkk_{\W\,\mn\rd T_K}$
by \Cref{thm:formality-map},
and similarly 
$\smash{\formalitymap_H \restriction^H_{T_H}}$ annihilates $\kkk_{\W\,\mn\rd H}$.
As $\Theta$ is congruent to 
$\smash{\gKo \otimes \gG \otimes \gHo}$ modulo 
$(\kkk_{\W\,\mn\rd K},\kkk_{\W\,\mn\rd H})$
by \Cref{thm:Theta-modulo-k},
we then have 
\[
\Psi\Theta 
	= 
(\form\mn\rest \otimes \id \otimes \mnn\form\mn\rest)\Theta 
	=
\smash{\form\mn\rest\gKo \otimes \gG \otimes \mnn\form\mn\rest\gHo}
\mathrlap.
\]
But taking cohomology of the cochain maps
\[
\vphantom{X^{X^X}}
\smash{
	\phantom{	\os{\formalitymap}	\lt}
	\H(BK) 			\xtoo{{\onecomponent{(\genmap_K)}}}
	\C(\W\,\mn\rd K) 	\os{\restriction}	\lt
	\C(\W\,\mn\rd T_K) 	\os{\formalitymap}	\lt
	\H(BT_K)
	\phantom{\xtoo{\onecomponent{(\genmap_K)}}}
}
\]
yields
\[\vphantom{X^X_{X_X}}
\H(BK) 		\os{\id}			\lt
\H(BK) 		\os{\restriction}	\lt
\C(BT_K) 	\os{\id}			\lt
\H(BT_K)
\]
and the differentials on $\H(BK)$ and $\H(BT_K)$
are zero, implying the cochain map
$f_K \o \restriction^K_{T_K}\o \onecomponent{(\genmap_K)}$ 
is itself $\restriction^K_{T_K}$;
and similarly for $\restriction^H_{T_H}$.	
	
	\medskip
	
\nd	(iii)
	We factor $\Psi$ 
	as $(\formalitymap \otimes \id \otimes \formalitymap)
	(\restriction \otimes \id \otimes \restriction)$
	and show both factors induce injections in cohomology.
	As $\formalitymap \otimes \id \otimes \formalitymap\mn\:
	\BB\bigl(\C(\W\,\mn\rd T_K),\C(\W\,\mn\rd G),\C(\W\,\mn\rd T_H)\mnn\bigr) 
	\lt
	\BB\bigl(\mspace{-.5mu}\H(BT_K),\C(\W\,\mn\rd G),\H(BT_H)\mnn\bigr)$
	induces the identity map between the $E_2$ 
	pages of the associated filtration spectral sequences,
	it is a quasi-isomorphism.
	As for $\restriction \otimes \id \otimes \restriction$,
	consider the map of Serre spectral sequences induced by 
	the map of (vertical) fibrations
	\[
		\xymatrix@R=3.5em@C=.75em{
				{}_{T_K}\mn G_{T_H} \ar[r]\ar[d]& BT_K \x BT_H \ar[d]\\
				{}_K\mn G_H 		\ar[r]		& BK \x BH\mathrlap{.}
			}
	\]
	The homotopy fiber of both 
	fibrations is $K/T_K \x H/T_H$.
	Thus the $E_2$ page of the right spectral sequence 
	is concentrated in even degree,
	implying the sequence collapses,
	and 
	the map of spectral sequences
	implies the left spectral sequence collapses as well,
	and it follows that the left fibration
	induces an injection in cohomology.
	But by the naturality clause of \Cref{thm:iso-bar-KGH}
	this injection is 
	$\H(\restriction \otimes \id \otimes \restriction)$.
\end{proof}

\begin{theorem}
	\label{thm:Theta-mult}
	The isomorphism~$\H(\quasiisomorphism)$ 
	of \Cref{thm:Theta-additive}
	is multiplicative.
\end{theorem}
\begin{proof}
	Since we know from \Cref{thm:HPsi} that
	$\H(\Psi)$ is injective and multiplicative
	it will be enough to show the map~%
	$
	\Psi\quism = \restriction \otimes \genmap_G \otimes \restriction
	$
	is multiplicative up to homotopy.
	As $\H(BK)$ and $\H(BH)$ are \CDGAs, 
	the \HGA operations $E_k$ are zero for $k \geq 0$ 
	by \Cref{def:xHGA-CDGA}, 
	so the product $\wt\mu$ of~\eqref{eq:def:two-sided-product} 
	reduces to 
	$a'[a_\bl]a'' \. b'[b_\bl]b'' 
		= 
	\pm a'b' \otimes [a_\bl] \bartimes [b_\bl] \otimes a''b''$,
	which is just
	the tensor permutation~$\Pi$ of~\eqref{eq:Pi}
	rearranging the factors in the correct order
	followed by the coordinatewise product%
	~$\mu^{\ox\,\! 3}$.
	Recalling from the proof of \Cref{thm:def-prod-bar} 
	that $\Pi$ gives a natural cochain isomorphism 	
	from the tensor-square of the two-sided bar construction
	to the two-sided twisted tensor product 
	with respect to the twisting cochains $t'$ and $t''$
	given by
	$\restriction^{\ox\mn 2}(t \otimes \h\e + \h\e \otimes t)$,
	we can transfer 
	the desired multiplicativity
	of $\Psi\quism$ up to homotopy
	to a question about maps
	\[
	\H(BK)^{\ox\mn 2} \,\mn\ox_{t'}\,\mn {}\big(\mnn\BB\H(BG)\mnn\big)^{\mn\ox\mn 2} \,\mn\ox_{t''}\,\mn
	\H(BH)^{\ox\mn 2}
		\lt
	\H(BT_K) 
		\ox_{\form\mn\rest\twistcochain}
	\BB\C(\W\,\mn\rd G)
		\ox_{\form\mn\rest\twistcochain}
	\H(BT_H)
	\mathrlap:
	\]
	we want to find a homotopy between the cochain maps
	\[
		(\restriction \otimes \genmap_G \otimes \restriction)
		(\mu \otimes \mu_{\BB\H(BG)} \otimes \mu)
		=
		\restriction\mu \otimes \genmap_G \mu_{\BB\H(BG)} \otimes \restriction\mu
	\]
	and 
	\[
		(\mu \otimes \mu_{\BB\C(\W\,\mn\rd G)} \otimes \mu) 
		(\restriction^{\ox\mn 2} \otimes \genmap_G^{\ox\mn 2} \otimes \restriction^{\ox\mn 2})
		=
		\mu\restriction^{\ox\mn 2} \otimes 
		 \mu_{\BB\C(\W\,\mn\rd G)}\genmap_G^{\ox\mn 2} \otimes
	    \mu\restriction^{\ox\mn 2}
	\mathrlap.
	\]
	On tensor factors,
	we have $\rest\mu = \mu\rest^{\ox\mn 2}$
	because $\rest$ is a ring map,
	while \Cref{thm:genmap-natural} provides a 
	coalgebra homotopy $H_\mu\:
	\genmap_G \B \mu_{\H(BG)}
		\hmt
	\Phi(\genmap_G \T \genmap_G)$
	such that $\twistcochain H_\mu$
	takes $\B_{\geq 1}\big(\mnn\H(BG)^{\ox\mn 2}\big)$
	into~$\smash{\kkk_{\W\,\mn\rd G}}$,
	so that the twisting cochains 
	$\form\mn\rest\twistcochain$
	annihilate the image of $H_\mu$.
	As the shuffle map $\nabla$ is a natural transformation 
	${\ox} \lt {\T}$
	of bifunctors,
	we may append the square of \Cref{prop:shuffle-tensor}:
	\[
		\xymatrix@C=1em@R=4.825em{
			\big(\B\H(BG)\mnn\big)^{\!\ox\mn 2}
				\ar@{}[r]^(.24){}="a"
				\ar "a";[r] 
				\ar@{}[r]_(.425){\nabla}
			\ar@{}[d]^(.05){}="b"
				\ar "b";[d]
			\ar@{}[d]_{\genmap_G \,\ox\, \genmap_G}&
			\BB\big(\H(BG)^{\ox\mn 2}\big)
			{\vphantom{\H(X_\bl)}} 
				\ar@{}[r]^(.385){}="a"
				\ar "a";[r]  
				\ar@{}[r]_(.5){\BB\mu{\vphantom{\H}}}
				\ar@{}[d]^(.05){}="b"
				\ar "b";[d]
					|(.44)\hole
					|(.4725){\genmap_G \,\mn\T\,\mn \genmap_G}
					|(.5)\hole&
			\BB\big(\H(BG)\big) 	\ar@{}[d]^(.05){}="b"
				\ar "b";[d]
			\ar@{}[d]^{\genmap_G}\\
			\big(\B\C(\W\,\mn\rd G)\mnn\big)^{\smash{\!\ox\mn 2}}
				\ar@{}[r]^(.24){}="a"
				\ar "a";[r] 
				\ar@{}[r]_(.425){\nabla}
			&
			\BB\big(\C(\W\,\mn\rd G)^{\ox\mn 2}\big)
				\ar@{}[r]^(.385){}="a"
				\ar "a";[r]  
				\ar@{}[r]_(.5){\Phi{\vphantom{\C}}}
			&
			\BB\C(\W\,\mn\rd G)\mathrlap.	
		}
	\]
	Recalling from \Cref{thm:xHGA-SHC} 
	that $\Phi \o \nabla = \smash{\mu_{\B\C(\W\,\mn\rd G)}}$
	and from \Cref{ex:shuffle-product}
	that $\B\mu \o \nabla = \mu_{\B\H(BG)}$,
	we then see $H_\mu \nabla$ is a \DGC homotopy 
	from $\smash{\genmap_G \mu_{\B\H(BG)}}$
	to $\smash{\mu_{\B\C(\W\,\mn\rd G)}\genmap_G^{\ox\mn 2}}$.
	As $\twistcochain H_\mu \nabla$ takes the coaugmentation coideal
		of $\big(\mspace{-.5mu}\B\H(BG)\mspace{-.75mu}\big){}^{\mnn\ox\mn 2}$ 
	into $\smash{\kkk_{\W\,\mn\rd G}}$,
	we may apply \Cref{thm:homotopy-twisted} 
	to see that $\form\mn\rest \otimes H_\mu \mn\nabla \otimes \form\mn\rest$
	is a homotopy doing what we wanted.
\end{proof}

This establishes the ring isomorphism of \Cref{thm:main}.
It remains to show this isomorphism is natural in inclusion diagrams.

\begin{theorem}
	\label{thm:iso-natural}
	Given a diagram
	\begin{equation}\label{eq:iso-natural}
	\begin{aligned}
		\xymatrix{
			\Ko \,\ar@{^{(}->}[r]\ar[d]&\Go\ar[d]	&\,\Ho \ar@{_{(}->}[l]\ar[d]\\
			\Kp \,\ar@{^{(}->}[r]		 &\Gp		&\,\Hp \ar@{_{(}->}[l]
		}
	\end{aligned}
	\end{equation}
	of continuous group homomorphisms
	such that the cohomology of each of 
	the groups' 
	classifying spaces is a polynomial ring over $\kk$
	and given 
	cocycle representatives for an irredundant set
	of polynomial generators of each cohomology ring,
	so as to define quasi-isomorphisms $\quismo$
	and $\quismp$ as in \eqref{eq:quism},
	the following induced square commutes:
	\[
	\xymatrix@C=4.25em@R=3em{
	\Tor^*_{\H(B\Gp)}\bigl(\H(B\Kp),\H(B\Hp)\mnn\bigr)
		\ar@{}[d]^(.075){}="b"
				\ar "b";[d]
			\ar@{}[d]_{\H(\quismp)}^\vertsim \ar[r]^{\Tor^*_\restriction(\restriction,\restriction)} & 
	\Tor^*_{\H(BG)}\bigl(\H(B\Ko),\H(B\Ho)\mnn\bigr) 
		\ar@{}[d]^(.075){}="b"
				\ar "b";[d]
			\ar@{}[d]^{\H(\quismo)}_\vertsim \\
	\H_{\Kp}(\Gp/\Hp) \ar[r] & \H_{\Ko}(\Go/\Ho).
	}
	\]
\end{theorem}

\begin{proof}
	Let $\defm{i_{\Ko}} \: T_{\Ko} \longinc \Ko$ and $\defm{i_{\Ho}}\:T_{\Ho} \longinc \Ho$ 
	be inclusions of maximal tori.
	We expand the description of $\Psi\quism$
	in \Cref{thm:HPsi} to the diagram in Figure \ref{fig:naturality-diagram}.
\begin{figure}
\xymatrix@C=4.125em@R=4em{
\dsp
	\H(B\Kp) \,\mn \ox_{\smash{\rest\twistcohom}} \,\mn \BB \H(B\Gp) \,\mn \ox_{\smash{\rest\twistcohom}} \,\mn \H(B\Hp) 
			\ar@{}[d]^(.085){}="b"^(.86){}="c"
	\ar "b";"c"
			\ar@{}[d]_{\quismp}
		\ar[r]^{\restriction \,\otimes\, \BB\restriction \,\otimes\, \restriction} 
&
\dsp 
	\H(B\Ko) \,\mn \ox_{\smash{\rest\twistcohom}} \,\mn \BB \H(B\Go) \,\mn \ox_{\smash{\rest\twistcohom}} \,\mn \H(B\Ho) 
	\ar@{}[d]^(.085){}="b"^(.86){}="c"
	\ar "b";"c"
			\ar@{}[d]^{\quismo} 
\\
\dsp
	\C(\W\,\mn\rd \Kp) \,\mn \ox_{\smash{\rest\twistcochain}} \,\mn \BB \C(\W\,\mn\rd \Gp) \,\mn \ox_{\smash{\rest\twistcochain}} \,\mn \C(\W\,\mn\rd \Hp)
	\ar@{}[d]^(.085){}="b"^(.86){}="c"
	\ar "b";"c"
			\ar@{}[d]_{\restriction
				\,\otimes\, \id \,\otimes\, 	
				\restriction
				} 
		\ar[r]^{\restriction \,\otimes\, \BB\restriction \,\otimes\, \restriction} 
&
\dsp 
	\C(\W\,\mn\rd \Ko) \,\mn \ox_{\smash{\rest\twistcochain}} \,\mn \BB \C(\W\,\mn\rd \Go) \,\mn \ox_{\smash{\rest\twistcochain}} \,\mn \C(\W\,\mn\rd \Ho)
	\ar@{}[d]^(.085){}="b"^(.86){}="c"
	\ar "b";"c"
			\ar@{}[d]^{\restriction \,\otimes\, \id \,\otimes\, \restriction} 
\\
\dsp
	\C(\W\,\mn\rd T_{\Ko}) \,\mn \ox_{\smash{\rest\twistcochain}} \,\mn \BB \C(\W\,\mn\rd \Gp) \,\mn \ox_{\smash{\rest\twistcochain}} \,\mn \C(\W\,\mn\rd T_{\Ho}) 	
	\ar@{}[d]^(.085){}="b"^(.86){}="c"
	\ar "b";"c"
			\ar@{}[d]_{\formalitymap \,\otimes\, \id \,\otimes\, \formalitymap} 
		\ar[r]_{\id \,\otimes\, \BB\restriction \,\otimes\, \id} 
&
\dsp 
	\C(\W\,\mn\rd T_{\Ko}) \,\mn \ox_{\smash{\rest\twistcochain}} \,\mn \BB \C(\W\,\mn\rd \Go) \,\mn \ox_{\smash{\rest\twistcochain}} \,\mn \C(\W\,\mn\rd T_{\Ho}) 
	\ar@{}[d]^(.085){}="b"^(.86){}="c"
	\ar "b";"c"
			\ar@{}[d]^{\formalitymap \,\otimes\, \id \,\otimes\, \formalitymap } 
\\
\dsp
	\H(BT_{\Ko}) \,\mn \ox_{\smash{\form\mn\rest\twistcochain}} \,\mn \BB \C(\W\,\mn\rd \Gp) \,\mn \ox_{\smash{\form\mn\rest\twistcochain}} \,\mn \H(BT_{\Ho}) 
		\ar[r]_{\id \,\otimes\, \BB\restriction \,\otimes\, \id} 
&
\dsp 
	\H(BT_{\Ko}) \,\mn \ox_{\smash{\form\mn\rest\twistcochain}} \,\mn \BB \C(\W\,\mn\rd \Go) \,\mn \ox_{\smash{\form\mn\rest\twistcochain}} \,\mn \H(BT_{\Ho})
}
\caption{A diagram of two-sided bar constructions.}%
\label{fig:naturality-diagram}
\end{figure}
	Note carefully the asymmetry in the bottom square:
	on the left we have $\Gp$ rather than $\Go$ but 
	the tori $T_{\Ko}$ and $T_{\Ho}$ are subtori of $\Go$ rather than $\Gp$.
	The composite $\form\mn\rest \otimes \id \otimes \form\mn\rest$ 
	after $\quism$ in the right column 
	is the map~$\Psi_{i_{\Ko},i_{\Ho}}$ of \eqref{eq:def-Psi}
	and the composite after $\quism$ the first column 
	is~$\Psi_{j_K,j_H}$,
	where $\defm{j_K}$ and $\defm{j_H}$ are
	respectively the compositions $T_{\Ko} \inc \Ko \to \Kp$
	and $T_{\Ho} \inc \Ho \to \Hp$.

	Our goal is to show the top square commutes in cohomology, 
	but as $\H(\Psi_{i_{\Ko},i_{\Ho}})$ is injective by \Cref{thm:HPsi}\ref{thm:HPsi-3},
	it will suffice to postcompose $\Psi_{i_{\Ko},i_{\Ho}}$ 
	to the two paths around the square 
	and show the resulting maps induce the same map in cohomology.
	One confirms tensor factor by tensor factor 
	that the bottom two squares commute,
	so it will be enough to find a cochain homotopy between
	the two paths along the outer rectangle.
	By \Cref{thm:HPsi}\ref{thm:HPsi-2}, 
	the composite along the upper right is
	\begin{equation}\label{eq:comp1}
		(\restriction^{\Ko}_{T_{\Ko}} \otimes \genmap_{\Go} \otimes \restriction^{\Ho}_{T_{\Ho}})
		(\restriction^{\Kp}_{\Ko} \otimes \BB \restriction^{\Gp}_{\Go} \otimes \restriction^{\Hp}_{\Ho})
			=
		\restriction^{\Kp}_{T_{\Ko}} \otimes \genmap_{\Go} \, \BB \restriction^{\Gp}_{\Go} \otimes \restriction^{\Kp}_{T_{\Ko}}
	\end{equation}
	and the composite along the lower left is 
	\begin{equation}\label{eq:comp2}
		(\id \otimes 
		\BB \restriction^{\Gp}_{\Go} \otimes 
		\id)
		(\restriction^{\Kp}_{T_{\Ko}} \otimes 
		\genmap_{\Gp} \otimes 
		\restriction^{\Hp}_{T_{\Ho}})
			=
		\restriction^{\Kp}_{T_{\Ko}} \otimes 
		\BB \restriction^{\Gp}_{\Go}\, \genmap_{\Gp} \otimes 
		\restriction^{\Hp}_{T_{\Ho}}\mathrlap{.}
	\end{equation}
	Now \Cref{thm:genmap-natural} provides a \DGC homotopy $H_\rho$
	between $\restriction^{\Gp}_{\Go}\genmap_{\Gp}$ and~$\gG\restriction^{\Gp}_{\Go}$
	such that $\twistcochain H_\rho$ 
	takes~$\BB_{\geq 1} \H(B\Gp)$
	into $\kkk_{\W\,\mn\rd \Go}$, 
	so that the two twisting cochains $\form\mn\rest\twistcochain$
	defining the common codomain of 
	\eqref{eq:comp1} and \eqref{eq:comp2}
	annihilate the image of $H_\rho$.
	Applying
	\Cref{thm:homotopy-twisted},
	we see $\form\mn\rest \otimes H_\rho \otimes \form\mn\rest$
	is the sought-after homotopy between 
	\eqref{eq:comp1} and \eqref{eq:comp2}.
\end{proof}

\begin{corollary}
	The isomorphism~$\H(\quism)$ 
	of \Cref{thm:main} 
	does not depend on 
	the choice of representatives defining~$\quism$ in~\eqref{eq:quism}.
\end{corollary}

\begin{proof}
	One can take the vertical maps to each be the identity in \eqref{eq:iso-natural},
	defining $\quism_0$ in terms of one set of representatives
	and $\quism_1$ in terms of another.
\end{proof}
\section{Example}\label{sec:example}
\revision{In this section we provide an example where \Cref{thm:main}
provides information not previously available from the standard techniques.
}
\revision{%
This example is 
in a sense minimal in that in similar examples
where the differentials form a regular sequence,
the fact an exterior algebra is a free \CGA allows one to resolve the multiplicative
extension problem,
and in examples with $\Tor^{\leq -2} = 0$,
one can resolve the extension problem by noting
elements of $\Tor^{-1}$ annihilate one other and lie in odd degree,
whereas $\Tor^0$ is concentrated in even degree.
}
\begin{example}\label{eg:computation}\revision{
	Let $S$ be the image of the embedding 
	$z \lmt \diag(z^{-4},z,z,z)$ of 
	$\U(1)$ in $\SU(4)$. 
	\Cref{thm:main} computes the ring $\H_S\big(\SU(4)/S;\Zt\big)$
	as 
	$\Tor^*_{\H(B\SU(4);\Zt)}\big(\H(BS;\Zt),\H(BS;\Zt)\big)$
	(previously, the Sullivan model noted by Kapovich
	showed the same with $\Zt$ replaced by $\Q$).
	One finds explicitly that the induced map 
	$\Zt[c_2,c_3,c_4] \simto \H B\SU(4) \to \H BS \simto \Zt[t_2]$
	takes the universal Chern classes $c_2, c_3, c_4$ respectively 
	to $-6t^2,-8t^3,-3t^4$,
	so one can compute the Tor, using the associated Koszul complex,
	as the cohomology of the \DGA
	$\Zt[s_2,t_2] \ox \ext[v_3,w_5,x_7]$
	with differential annihilating $s$ and $t$ and with
	\eqn{
		dv_3 &= 6(s^2 - t^2),\\
		dw_5 &= 8(s^3 - t^3),\\
		dx_7 &= 3(s^4 - t^4).
	}
	Setting $p_{2j} = (s^j-t^j)/(s-t)$ for $2 \leq j \leq 4$,
	one finds $\Tor^{0} \iso \Zt[s,t] / (s-t)(6p_2,8p_4,3p_6)$;
	that $\Tor^{-1}$ is the quotient of the $\Zt[s,t]$-module
	generated by the cocycles
	\eqn{
	a &= 3p_2 w_5 - 4p_4 v_3 = \frac 1{2(s-t)}d(v_3 w_5),\\
		b &= 2p_2 x_7 -  p_6 v_3= \frac 1{3(s-t)}d(v_3 x_7),\\
		c &= 8p_4 x_7 - 3p_6 w_5 = \frac 1{s-t}d(w_5 x_7),
	}
	which satisfy $p_6 a - 4p_4 b + p_2 c = 0$,
	by the submodule generated by
		$d(v_3 w_5),
		d(v_3 x_7),
		d(w_5 x_7)$;
	and that $\Tor^{-2}$ is the quotient
	of the free $\Zt[s,t]$-module
	generated by the class of 
	\[
		e =
		3p_6 v_3 w_5 - 8p_4 v_3 x_7 + 6p_2 w_5 x_7 = 
		\frac 1{(s-t)}{d(v_3 w_5 x_7)}\mathrlap.
	\]
	by the submodule generated by $d(v_3 w_5 x_7)$.
	All told, one finds 
	$\HS\big(\SU(4)/S;\Zt\big)$
	is the quotient of 
	$\Zt[s_2,t_2,e_{14}] \ox \ext[a_7,b_9,c_{11}]$
	by the ideal generated by
	\eqn{
		(s^3+s^2t+st^2+t)a-4(s^2+st+t^2)b+(s+t)c,\\
		2(s-t)a,\quad
		3(s-t)b,\quad
		(s-t)c,\quad
		(s-t)e,\\
		ab-(s+t)e,\quad
		ac-4(s^2+st+t^2)e,\quad
		bc-(s^3+s^2t+st^2+t^3)e,\\
		ae,\quad
		be,\quad
		ce,\quad
		e^2
		\mathrlap.
	}
This computation resolves, for example, the multiplicative extension problem 
for products of the classes corresponding to $a,b,c,e$
on the $E_\infty$ page of the Serre spectral sequence 
beginning with $E_2 =
\H(BS \x BS) \otimes \H \SU(4)$ and
converging to $\HS\big(\SU(4)/S;\Z\big)$.
}

\end{example}

\bs

\bs

\newpage
\nd{\textsf{{\huge Appendix} (with Matthias Franz)
	}}

\vspace{-1.5em}

\appendix
\section{A product on the two-sided bar construction}\label{sec:prod}

Given that an \HGA structure defines a \DGA structure
on the bar construction of a \DGA, 
one might hope 
\HGA homomorphisms $A' \from A \to A''$ 
similarly induce a \DGA structure on the
two-sided bar construction.
We cannot assert this,
but we can at least define a non-associative product
sufficient to prove the variant Eilenberg--Moore \cref{thm:EM-product} 
needed for the proof of \Cref{thm:main}.

%
%
%

\begin{theorem}\label{thm:def-prod-bar}
	Let $A' \os{f'}\from A \os{f''}\to A''$ be \HGA
	homomorphisms.
	Then there exists a cochain map 
	\[
	\defm{\wt\mu}\: \BB(A',A,A'') \ox \BB(A',A,A'') \lt \BB(A',A,A'')\mathrlap,
	\]
	which we think of as a product,
	natural in the sense that given a commutative diagram
	\[
	\xymatrix{
		A' \ar[d]_{g'}& 
		A \ar[d]|(.425)\hole|{g}|(.55)\hole \ar[l]\ar[r]& 
		A'' \ar[d]^{g''}\\
		B'& 
		B \ar[l]\ar[r]& 
		B'' 
	}
	\]
	of \HGA homomorphisms, 
	the induced map $\B(g',g,g'')\: \B(A',A,A'') \lt \B(B',B,B'')$ 
	of \Cref{def:bar-map-dga}
	is multiplicative.
\end{theorem}
\begin{remark}
	It requires the extended \HGA operations $F_{p,q}$
	to define the natural homotopy $h^{\mr c}$
	witnessing the homotopy-commutativity axiom 
	for Franz's natural \SHCA structure map~$\Phi_A$ 
	described in \Cref{thm:xHGA-SHC},
	but to define $\Phi_A$ itself 
	and show that $\Phi_A \o \nabla = \mu_{\B A}$,
	which is all we require here, 
	one needs only that $A$ be an \HGA.
\end{remark}
\begin{proof}[Proof of \Cref{thm:def-prod-bar}]
	By \Cref{thm:xHGA-SHC},
	there exists a \DGC map $\Phi_{\mn A}\: \BB(A \ox A) \lt \BB A$
	giving an \SHC-algebra structure on $A$,
	and similarly for $A'$ and $A''$.
	We set 
	\begin{equation}\label{eq:mu-tilde}
	\defm{\wt\mu} 
	\ceq
	\underbrace{
		\BB(\Phi_{\mn A'},\Phi_{\mn A},\Phi_{\mn A''}) 
		\o
		(\id_{A' \ox A'} \otimes \nabla \otimes \id_{A'' \ox A''})	
	}_{\BB(\Phi_{\mn A'},\,\mu_{\BB A},\,\Phi_{\mn A''})} 
	{\mn} \o {\Pi}
	\mathrlap,
	\end{equation}
	where the $\BB(-,-,-)$ maps 
	are as defined in \Cref{def:Gamma-generalized}
	and the initial map
	\begin{equation}\label{eq:Pi}
	\begin{aligned}
	\defm \Pi\:
	\BB(A',A,A'') \ox \BB(A',A,A'') 
	&\lt 
	(A' \ox A') \ox_{t'} {}
	(\BB A \ox \BB A) \ox_{t''} {}
	(A'' \ox A'')
	\\
	a'[a_\bl]a'' \otimes b'[b_\bl]b''
	&\lmt 
	\pm a' \otimes b' \otimes [a_\bl] \otimes [b_\bl] \otimes a'' \otimes b''
	\end{aligned}
	\end{equation}
	is the tensor permutation~$(2 \ 3 \ 5\ 4)$.
	Here the 
	two-sided twisted tensor product
	is determined by the twisting cochains~%
	\eqn{\defm{t'} = 
		(f')^{\ox\mn 2}t_{A^{\ox\mn 2}}\nabla =
		(f')^{\ox\mn 2}
		(t_{A} \otimes \hA\e_{\BB A} + \hA\e_{\BB A} \otimes t_{A})
		&\: (\BB A)^{\ox\mn 2} \lt (A')^{\ox\mn 2} \mathrlap,\\
		\defm{t''} = 
		(f'')^{\ox\mn 2}t_{A^{\ox\mn 2}}\nabla =
		(f'')^{\ox\mn 2} (t_{A} \otimes \hA\e_{\BB A} + \hA\e_{\BB A} \otimes t_{A})
		&\: (\BB A)^{\ox\mn 2} \lt (A'')^{\ox\mn 2}
	}
	for $\nabla$ the shuffle map of 
	\Cref{ex:shuffle-product},
	and it is not hard to check 
	$\Pi$ is a cochain map.
	To see 
	\[
	\id \ox \nabla \ox \id\: 
	(A' \ox A') \ox_{t'} {(\BB A \ox \BB A)} \ox_{t''} {(A'' \ox A'')}
	\lt
	(A' \ox A') 
	\ox_{t'} {}
	\BB (A \ox A) 
	\ox_{t''} 
	{(A'' \ox A'')}
	\] 
	is a cochain map,
	it is enough by \Cref{thm:homotopy-twisted}
	to observe 
	$t' = (f' \otimes f')t_{A \ox A} \nabla = t_{A' \ox A'}\BB(f' \otimes f') \nabla $
	and similarly $t'' = t_{A'' \ox A''}\BB(f'' \otimes f'') \nabla$.
	To see 
	\[
	\BB(\Phi_{\mn A'},\Phi_{\mn A},\Phi_{\mn A''})\: 
	\BB(A' \ox A', A \ox A, A'' \ox A'') 
	\lt 
	\BB(A',A,A'')
	\]
	is a well-defined cochain map, 
	by \Cref{def:Gamma-generalized}
	it is enough to	note the diagram
	\[
	\xymatrix@C=2.25em@R=5em{
		\BB (A' \ox A')  \ar[d]_{\Phi_{\mn A'}}&
		\BB (A 	\ox A) 	 \ar[d]|(.45)\hole|{\Phi_{\mn A}}|(.55)\hole
		\ar[l]_(.475){\BB(f')^{\ox\mn 2}} 
		\ar[r]^(.475){\BB(f'')^{\ox\mn 2}}&
		\BB (A'' \ox A'') \ar[d]^{\Phi_{\mn A''}} \\
		\BB A' &
		\BB A  \ar[l]^{\BB f'}\ar[r]_{\BB f''}&
		\BB A''
		\mathrlap.
	}
	\]
	commutes by the naturality of $\Phi$ in \HGAs.
	
	Naturality follows because 
	$\Pi$ is natural in sextuples of cochain complexes,
	$\nabla$ is natural in pairs of cochain complexes,
	$\Phi$ is natural in extended \HGA maps,
	and $\BB(-,-,-)$ is functorial in triples of \DGA maps. 
\end{proof}

To make real use of this product,
we will require a more explicit formula.

\begin{theorem}\label{thm:formula-prod-bar}
	The cochain map $\wt\mu$ of \Cref{thm:def-prod-bar}
	is given  
	in terms of the notations set
	in \Cref{def:HGA} and \Cref{notation:E}
	by
	\begin{equation}
	\begin{aligned}\label{eq:def:two-sided-product}
	a'[a_\bl]a'' \ \otimes \ b'[b_\bl]b''
	\quad\lmt\quad&
	\pm a'b'						
	&\otimes&\ 
	\ \ \  [a_\bl] \bartimes [b_{(2)}]
	\ \otimes\  
	\EEE\big(a''[f''b_{(3)}]\big)b''
	\\
	&\pm a' \EEE\big(f'a_1[b'|f'b_{(1)}]\big) \!\!\!
	&\otimes &
	\ \ [a_{(2)}] \bartimes [b_{(2)}]
	\ \otimes\ 
	\EEE\big(a''[f''b_{(3)}]\big)b''
	\mathrlap,
	\end{aligned}
	\end{equation}
	where the first term is the sum over all decompositions $[b_\bl] = \defm{[b_{(2)}]} \otimes \defm{[b_{(3)}]}$ 
	into two tensor factors
	and the second is the double sum over decompositions 
	$[a_\bl] = [a_1] \otimes \defm{[a_{(2)}]}$
	where the first tensor factor is of length~$1$
	and decompositions 
	$[b_\bl] = \defm{[b_{(1)}]} \otimes \defm{[b_{(2)}]} \otimes \defm{[b_{(3)}]}$ 
	into three tensor factors,
	and where 
	we recall from \Cref{def:vanishing-components}
	our convention
	that $\EEE(f' a_1 [b'|f'b_{(1)}]) = 0$
	when $b' = 1$.
	The signs are those imposed by the Koszul convention;
	explicitly, 
	the first sum is 
	\quation{\label{eq:two-sided-first-term}
		\overbrace{
			\Big(\mu_{A'} \,\otimes\, 
			\mu_{\BB A} \,\otimes\, 
			\underbrace{
				\mu_{A''}\big(\EEE(\id_{A''} \ox \BB\mn f'') \ox \,\mn\id_{A''}\big)
				\tau_{[b_{(3)}\mnn];\, a''}
			}_{\defm{T_{0,3}}}
			\Big)
		}^{\defm{T_0}}
		\o\,
		(\id_{A' \ox A'} \otimes \,\mn\id_{\BB A} \otimes \D_{\BB A} \otimes \,\mn\id_{A'' \ox A''})
		\Pi
		\mathrlap,\quad
	}
	for $\defm{\tau_{[b_{(3)}\mnn];\, a''}}$ the tensor shuffle 
	taking the factors $[b_{(3)}] \otimes a'$
	to $(-1)^{|[b_{(3)}]|\.|a'|} \, a' \otimes [b_{(3)}]$
	and $\Pi$ the permutation $(2\ 3\ 5\ 4)$ from \eqref{eq:Pi},
	and the second is
	\begin{multline}
	\label{eq:two-sided-second-term}
	\overbrace{
		\Big(
		\mn
		\overbrace{
			\mu_{A'}\big(\id_{A'} \otimes \EEE(f' \susp \otimes \desusp \otimes \BB\mn f')\big)
			\tau_{b';\, [a_1]}
		}^{\defm{T_{1,1}}}
		\,\otimes\,
		\mu_{\BB A}
		\,\otimes\,
		\overbrace{
			\mu_{A''}\big(
			\EEE(\id_{A''} \otimes \BB\mn f'') \otimes \id_{A''}
			\big)
			\tau_{[b_{(3)}\mnn];\, a''}\mn
		}^{\defm{T_{1,3}}}
		\Big)
	}^{\defm{T_1}}
	\qquad\qquad
	\\
	\o
	\Big(
	\id_{A'} \otimes (\id_{A'}-\e_{A'})
	\otimes 
	\big(\mnn
	({\pr_1} \otimes {\id_{\BB A}}) 
	\ox 
	{(\id_{\BB A} \otimes \id_{\BB A})}
	\ox 
	{(\e_{\BB A} \otimes \id_{\BB A})}
	\mnn\big)
	\iter{\D_{\BB A \ox \BB A}}3 
	\ox 
	\id_{A'' \ox A''}
	\Big)
	\Pi
	\mathrlap,\quad
	\end{multline}
	where $\defm{\tau_{b';\, [a_1]}}$ is again a tensor transposition.
\end{theorem}

Concatenations $\concatenate$ (as defined in \Cref{def:bar}) 
have been omitted in 
the statement and proof to aid legibility (so far as possible).

\begin{proof}
	We first establish the formula modulo $2$. 
	The composition of the first two maps in \eqref{eq:mu-tilde}
	takes a pure tensor $a'[a_\bl]a'' \ox b'[b_\bl]b''$
	to $(a' \ox b') \ox \nabla\big([a_\bl] \ox {[b_\bl]}\big) \ox {(a'' \ox b'')}$.
	Recall 
	$\nabla\big([a_\bl] \ox {[b_\bl]}\big)$
	is the signed sum of terms 
	$\defm{[\a \otimes \b]_\bl} \ceq 
	[\a_1 \otimes \b_1| \cdots |\a_\ell \otimes \b_{\defm\ell}]$ 
	where each $\a_j \otimes \b_j$ is either 
	$a_m \ox 1$ or~$1 \ox b_m$ for some $m$
	and the indices of the $a$-arguments appear in increasing order
	as one encounters them from left to right,
	as do the indices of the $b$-arguments.
	The map $\BB(\Phi_{\mn A'},\Phi_{\mn A},\Phi_{\mn A''}) $ first breaks each term 
	$[\a \otimes \b]_\bl$ into three factors 
	$[\a \otimes \b]_{(1)} = [\a_1 \otimes \b_1|\cdots|\a_p \otimes \b_p]$, 
	$[\a \otimes \b]_{(2)}$,
	$[\a \otimes \b]_{(3)}$  
	via $\smash{\iter{\D_{\BB(A^{\ox\mn 2})}}3}$,
	and applies $f'$ to each $\a$ and $\b$ in the first block and $f''$ to each in the third block.
	The verification for each of the three tensor factors will run in parallel.
	
	Now there is a case distinction to be made.
	Write $\defm{\a'_j} = f'\a_j$ and so on. 
	If $a' \ox b' \in \ol{A' \ox A'}$ 
	and~$a'' \ox b'' \in \ol{A'' \ox A''}$, 
	then by \eqref{eq:Gamma-generalized},
	$\BB(\Phi_{\mn A'},\Phi_{\mn A},\Phi_{\mn A''})$ takes 
	$(a' \ox b') \ox \nabla\big([a_\bl] \ox {[b_\bl]}\big) \ox {(a'' \ox b'')}$
	to 
	\quation{\label{eq:second-two-maps-mod-2}
		\pm
		t_{A'}\Phi_{\mn A'}[a' \otimes b'|\a'_1 \otimes \b'_1|\cdots|\a'_p \otimes \b'_{\defm p}]
		\ \ \otimes\ \ 
		[\a \otimes \b]_{(2)}
		\ \ \otimes\ \ 
		t_{A''}\Phi_{\mn A''}[\a''_{1} \otimes \b''_{1}|\cdots|\a''_{\defm q} \otimes \b''_{q}|a'' \otimes b'']
		\mathrlap,\ \ 
	}
	where $p + \ell\big([\a \otimes \b]_{(2)}\big) + q = \ell(a_\bl) + \ell(b_\bl)$.
	If instead $a' \otimes b' \in \im \h_{A' \ox A'} = \kk$,
	then the associated clause $\Upsilon'= \id_\kk \otimes \e$ from \eqref{eq:Gamma-generalized}
	applies instead,
	so the first tensor factor in \eqref{eq:second-two-maps-mod-2}
	is replaced 
	with~$a'b'$, 
	and~$p = 0$;
	similarly, if $a'' \otimes b'' \in \im \h_{A'' \ox A''}$,
	then instead the third tensor factor in \eqref{eq:second-two-maps-mod-2} 
	is $a''b''$, and $q = 0$.
	In either case, these simpler factors agree with those given in
	\eqref{eq:def:two-sided-product},
	\eqref{eq:two-sided-first-term}, and
	\eqref{eq:two-sided-second-term}
	by our observations and conventions from \Cref{def:vanishing-components}
	on the vanishing of $\EE$,
	so for the remainder of the verification we may assume 
	$a' \otimes b' \in \ol{A' \ox A'}$ and $a'' \ox b'' \in \ol{A'' \ox A''}$.
	
	At this point it becomes necessary to recall the unsigned formula~\cite[(4.2)]{franz2019shc}
	for $\tA \Phi_{\mn A}$: 
	\quation{\label{eq:Phi-mod-2}
		\tA \Phi_{\mn A} [\a \otimes \b]_\bl 
		= 
		\smash{\sum_{\smash{\s}} \pm \a_1 \EE[\a_2|\b_\bl] \cdots \EE[\a_n|\b_\bl] \b_n\mathrlap,}
	} 
	where the sum is over those permutations $\defm \s$
	separately preserving the orders of the indices of the $\a$-arguments 
	and those of the $\b$-arguments,
	and such that, moreover, at every point,
	reading left to right,
	one has encountered the symbol $\a$ more often than $\b$.
	The key point in simplifying this formula in the present context
	is that 
	$\EE[\a_j|\b_\bl]$ vanishes whenever one of the $\b$-letters
	other than $\b_n$ is $1$
	by the convention set in \Cref{def:vanishing-components}.
	
	Thus in order for 
	$t_{A'}\Phi_{\mn A'}\big([a' \otimes b'] \otimes {[\a' \otimes \b']_\bl}\big)$
	to be nonzero, one of two things must happen.
	The first possibility is that
	$p = \ell[\a' \otimes \b']_\bl = 0$,
	so that $t_{A'}\Phi_{\mn A'}[a' \otimes b'] = \pm a'b'$.
	The other possibility is that $b',\b'_1,\ldots,\b'{}\!_{p -1}$
	are all non-$1$.
	Particularly, if $b' = 1$, then $\pm a'b'$ is again the only term,
	agreeing with \eqref{eq:def:two-sided-product},
	\eqref{eq:two-sided-first-term}, and
	\eqref{eq:two-sided-second-term}, 
	so from now on we will assume $b' \in \ol{A'}$.
	Since for each $j$ one of $\a'_j$ and $\b'_j$ is $1$,
	it follows that $\a'_1 = \cdots = \a'_{p-1} = 1$.
	As $\smash{\EE_{1,0}\big([1] \otimes {[]}\big) = 1}$ 
	and otherwise~$\smash{\EE_{1,\bl}\big([1] \otimes {[\b_\bl]}\big) = 0}$
	by the conventions of \Cref{def:vanishing-components}, 
	for a term in \eqref{eq:Phi-mod-2} to be nonzero
	one needs $\a'_{p} \neq 1$, so that $\b'_{p} = 1$.
	Thus the only relevant terms have 
	\eqn{[\a' \otimes \b']_\bl &= 
		[1 \otimes \mnn f'b_1|\cdots|1 \otimes \mnn f'b_{p-1}|f'a_1 \otimes 1]\mathrlap,\\
		t_{A'}\Phi_{\mn A'}\big([a' \otimes b'] \otimes {[\a'  \otimes \b']_\bl}\big)
		&= \pm a' \EE_{1,p}\big([f'a_1] \otimes {[b'|f'b_{(1)}]}\big)
		\mathrlap,
	}
	where $[f' b_{(1)}] = [f'b_1|\cdots|f'b_{p-1}]$.
	
	Similarly, in order that 
	$t_{A''}\Phi_{\mn A''}[\a''_1 \otimes \b''_1|\cdots|\a''_q \otimes \b''_q|a'' \otimes b'']$ 
	be nonzero, 
	the arguments $\b''_j$ must be non-$1$,
	so each $\a''_j$ is $1$.
	Thus the relevant tensor factor is
	\[
	t_{A''}\Phi_{\mn A''}\big([1 \ox f''b]_{(3)} \otimes [a'' \ox b'']\big) 
	= \pm \EE_{1,q}\big([a''] \otimes [f''b_{(3)}]\big)b''
	\mathrlap,
	\]
	which vanishes by convention if $a'' = 1$.
	As for the middle tensor factor, we know from \Cref{thm:xHGA-SHC}
	that $\Phi \o \nabla = \mu_{\BB A}$.
	Combining our descriptions of the three tensor factors, 
	we have established \eqref{eq:def:two-sided-product}.
	
	\medskip
	
	It remains to determine the signs. 
	The sign for the permutation $\Pi$
	of \eqref{eq:Pi}
	is by definition the Koszul permutation sign.
	For the composition of the second two maps,
	recall from 
	\eqref{eq:Gamma-generalized}
	that before involving $\Phi_{\mn A'}$ and $\Phi_{\mn A''}$,
	the map $\BB(\Phi_{\mn A'},\,\mu_{\BB A},\,\Phi_{\mn A''})$ first breaks 
	$[\a \otimes \b]_\bl = \nabla\big([a_\bl] \otimes [b_\bl]\big)$
	into three chunks via $\iter{\D_{\BB(A \ox A)}}3$.
	As $\nabla$ is a \DGC map $\BB A \ox \BB A \lt \BB(A \ox A)$
	by \Cref{ex:shuffle-product},
	we equally well have $\iter{\D_{\BB(A \ox A)}}3\nabla = 
	\nabla^{\otimes 3}\iter{\D_{\BB A \ox \BB A}}3$.
	Applying the definition \eqref{eq:Gamma-generalized}, 
	suppressing concatenations~$\concatenate$, 
	gives
	\begin{equation}
	\begin{multlined}\label{eq:last-three-maps}
	\BB(\Phi_{\mn A'},\Phi_{\mn A},\Phi_{\mn A''})
	(\id_{A' \ox A'} \,\ox\, \nabla \,\ox\, \id_{A'' \ox A''})
	\\=
	\underbrace{
		\Big(\mn
		\overbrace{\mn
			t_{A'}\Phi_{\mn A'}\big(\desusp_{A' \ox A'} \,\otimes\, \BB(f' \otimes f')\nabla\big)
		}^{{\defm{Q_1}}} 
		\,\ox\, 
		\Phi_{\mn A}\nabla
		\,\ox\, 
		\overbrace{
			t_{A''}\Phi_{\mn A''}\big(\BB(f''\otimes f'')\nabla \,\otimes\, \desusp_{A''\ox A''}\big)
			\mn}^{\smash{\defm{Q_3}}}
		\mn\Big)
	}_{\smash{\defm Q}}
	\\
	\o
	(\id_{A'\ox A'} \,\ox\, \iter{\D_{\BB A \ox \BB A}}{3} \,\ox\, \id_{A'' \ox A''})\mathrlap.
	\end{multlined}
	\end{equation}
	We will evaluate $Q$
	on terms
	\[
	\big(
	a' \otimes b' \otimes [a_{(1)}] \otimes [b_{(1)}]
	\big) 
	\otimes 
	\big(
	[a_{(2)}] \otimes [b_{(2)}]
	\big) 
	\otimes 
	\big(
	[a_{(3)}] \otimes [b_{(3)}] \otimes a'' \otimes b''
	\big)
	\]
	and compare with the evaluations of $T_0$ and $T_1$
	from \eqref{eq:two-sided-first-term} and \eqref{eq:two-sided-second-term},
	recalling
	from the discussion of nonzero values of $t_{A'}\Phi_{\mn A'}$
	and $t_{A''}\Phi_{\mn A''}$ in the previous paragraph
	that the only nonzero terms in \eqref{eq:last-three-maps}
	arise when $ [a_{(3)}] = []$ is of length $0$ and either
	we have 
	$[a_{(1)}] = [b_{(1)}] = []$ of length $0$ as well or
	else $[a_{(1)}] = [a_1]$ is of length $1$.
	
	As each of the three tensor factors of $Q$, of $T_0$,
	and of $T_1$ is of degree zero,
	it will be enough to compare these factors separately.
	The middle tensor factor $\Phi_{\mn A}\nabla$ of $Q$ is $\mu_{\BB A}$
	by \Cref{thm:xHGA-SHC},
	as in $T_0$ and $T_1$. 
	For the first and third factors, we will
	need to be explicit about the sign for $\Phi$.
	The signed version of formula~\eqref{eq:Phi-mod-2} 
	for 
	$(\Phi_{\mn A})_{(n)} = \tA \o (\Phi_{\mn A})|_{\BB_n A} \o (\desusp_{A \ox A})^{\otimes n}$
	is
	\quation{\label{eq:Phi-sign}
		(\Phi_{\mn A})_{(n)} 
		\ceq 
		(-1)^{n-1} \sum_\s \iter{\mu_A}\bl 
		(\id_A \otimes \Tensor_m E_{\ell_m} \otimes \id_A)\s
		\mathrlap,
	}
	where the sum is over shuffles $\s$ as in \eqref{eq:Phi-mod-2},
	so that in particular the sequence of $(\ell_m)$ gives a partition of $n-1$, 
	apportioning
	subsequences of $\beta_1,\ldots,\beta_{n-1}$ 
	as arguments of
	$\EE[\alpha_2| {-} ]$ through $\EE[\alpha_n | {-} ]$,
	and $E_\ell$ is $\EE_{1,\ell}(\desusp)^{\otimes 1+\ell}$ as in \Cref{notation:E}; 
	and similarly for $\Phi_{\mn A'}$ and  $\Phi_{\mn A''}$.
	
	To interpret $Q_1$,
	note that $\ell[a_{(1)}]$ is either $0$ or $1$.
	In the case $[a_{(1)}] = []$,
	the only term  we need to worry about is
	$t_{A'}\Phi_{\mn A'}\desusp_{A' \ox A'} = (\Phi_{\mn A'})_{(1)}$,
	which is $\mu_{A'}$ by \Cref{def:SHC}.\ref{def:Phi-mu}.
	This is the first factor of $T_0$ from \eqref{eq:two-sided-first-term}.
	When $[a_{(1)}] = [a_1]$,
	we have seen from the discussion above 
	that the only terms on which
	$t_{A'}\Phi_{\mn A'}$
	will potentially not vanish  
	are those of the form
	$[a' \otimes b'|1 \otimes f'b]_{(1)} \otimes [f'a_1 \otimes 1]$.
	Here the notation $[a' \otimes b'|1 \otimes f'b]_{(1)}$
	means some initial subword
	$[a' \otimes b'|1 \otimes f'b_1|\cdots|1 \otimes f'b_\ell]$,
	which is $[a' \otimes b']$ if $\ell = 0$.
	
	A term of $\nabla\big([a_1] \otimes [b_{(1)}])$
	leading to a term of this form 
	is created by the operation 
	\[
	\defm\upsilon\: 
	[a_1] \otimes [b_{(1)}] \lmt [a_1 \otimes 1|1 \otimes b]_{(1)}
	\]
	followed by the shuffle
	$\defm \tau = \tau_{[a_1\ox 1];\, [1\ox b]_{(1)}}$
	moving $[a_1 \otimes 1]$ past each $[1 \otimes b_j]$,
	resulting in 
	\[
	\defm x 
		\ceq 
	\pm 
	[a' \otimes b'] \otimes [1 \otimes f'b]_{(1)} \otimes [f'a_1 \otimes 1]\mathrlap.
	\]	
	To determine the sign explicitly, \eqref{eq:Phi-sign}
	tells us we should rewrite things in terms of the operations~$E_\ell$.
	When we apply $\Phi_{\mn A'}$
	to $x$,
	we have seen above that the only nonzero term
	is $\pm a'E_{p}(a_1;\, b',b_{(1)})$,
	which comes from the summand in \eqref{eq:Phi-sign} 
	corresponding to the shuffle
	\[
	\defm{\s'}\: 
	(a' \otimes b') \otimes\,\mnn \Tensor_{j=1}^{p-1} \,\mn (1 \otimes f'b_j) \otimes (f'a_1 \otimes 1)
	\lmt (-1)^{|a_1	|\.(|b'|+|b_{(1)}|)}
	a' \otimes 1^{\otimes p-1} \otimes f' a_1 \otimes b' \otimes f'b_{(1)}
	\otimes 1
	\mathrlap.
	\]
	If we write $\defm \vk$ for the map omitting tensor factors $1$,
	then this nonzero term can be written as
	\[
	\mu_A'\smash{(\id_{A'} \otimes E_{p})\tau_{b' \ox b_{(1)};\, a_1} }\vk (\desusp_{A' \ox A'})^{\otimes p+1}x
	\mathrlap,
	\]
	where 
	$\defm{\tau_{b' \ox b_{(1)};\, a_1}}$ is the tensor shuffle
	moving $a_1$ past $b' \otimes b_{(1)}$.
	Since $(\desusp)^{\otimes n}\susp^{\otimes n}$ is multiplication by 
	$(-1)^{\binom{n}{2}}
	$,
	we conclude from
	\eqref{eq:Phi-sign} that on 
	${(A' \ox A') \ox \BB_1 A \ox \BB_p A}$, we have
	\eqn{
		Q_1 
		&= 
		t_{A'}\Phi_{\mn A'} 
		\o 
		\big(\desusp_{A' \ox A'} \otimes \BB(f')^{\ox\mn 2}\tau\upsilon\big)\\
		&=
		(-1)^{\binom{p+1}{2}} (\Phi_{\mn A'})_{(p+1)}
		\susp_{A'\ox A'}^{\otimes 1+p}
		\o
		\big(\desusp_{A' \ox A'} \otimes 
		\BB(f')^{\ox\mn 2} \tau\upsilon\big)
		\\
		&=
		(-1)^{\binom{p+1}{2}} 
		\Big(
		(-1)^{p}
		\mu_{A'}\big(\id_{A'} \otimes \EE_{1,p}(\desusp_{A'})^{\otimes 1+p}\big)\tau_{b' \ox b_{(1)};\, a_1} \vk 
		\Big)
		(-1)^{p}\big(\id_{A' \ox A'} \otimes 
		\susp_{A'\ox A'}^{\otimes p}\BB(f')^{\ox\mn 2}\tau\upsilon\big)
		\\
		&= 
		(-1)^{\binom{p+1}{2}} 
		\mu_{A'}\Big(\id_{A'} \otimes \EEE
		\big(\id_{A'} \otimes (\desusp_{A'})^{\otimes p}\big)\Big)\tau_{b' \ox b_{(1)};\, a_1} \vk
		\big(
		\id_{A' \ox A'} \otimes 
		\susp_{A' \ox A'}^{\otimes p}\BB(f')^{\ox\mn 2}\tau\upsilon
		\big)
		\mathrlap,
	}
	where the accounting of signs is as follows:
	the $(-1)^{\binom{p+1}{2}}$ comes from 
	$(\Phi_{\mn A'})_{(1+p)} \ceq t_{A'}\Phi_{\mn A'}\susp_{A'\ox A'}^{\otimes 1+p}$
	in \eqref{eq:Phi-sign},
	the first $(-1)^{p}$ 
	is the leftmost term
	of the right-hand side of \eqref{eq:Phi-sign},
	and the second~$(-1)^{p}$ comes 
	from moving $\susp_{A' \ox A'}^{\otimes p}$ past~$\desusp_{A'\ox A'}$.
	
	We must compare the sign 
	of the value of $Q_1$ on $\big([a' \otimes b'] \otimes [a_1|b_{(1)}]\big)$
	with that of
	the first factor~$T_{1,1} = 
	\mu_{A'}\big(\id_{A'} \otimes \EEE(f' \susp \otimes \desusp \otimes \BB f')\big) \tau_{b';\, [a_1]}$
	from \eqref{eq:two-sided-second-term}.
	This is done symbolically by determining which tensor factors are moved past which others in the computation.
	For these purposes, the factors
	$\mu_{A'}$, $\id_{A'}$, $\EEE$, $\vk$, 
	$f'$, $\BB(f' \otimes f')$, and $\upsilon$
	of degree $0$ are invisible,
	as are the tensor factors $1$ of degree $0$
	produced by $\upsilon$ and deleted by $\vk$.
	In the below diagram  
	$Q_1$ is in the left column and~$T_{1,1}$ in the right,
	and in each row of each column,
	operations are listed on the left
	and arguments on the right,
	the result of applying the operations 
	appearing on the right in the following row. 
	As usual, moving two symbols $y$ and $z$ past each other incurs the sign $(-1)^{|y|\.|z|}$,
	and the rule goes equally for $y = \susp^{\pm 1}$.
	
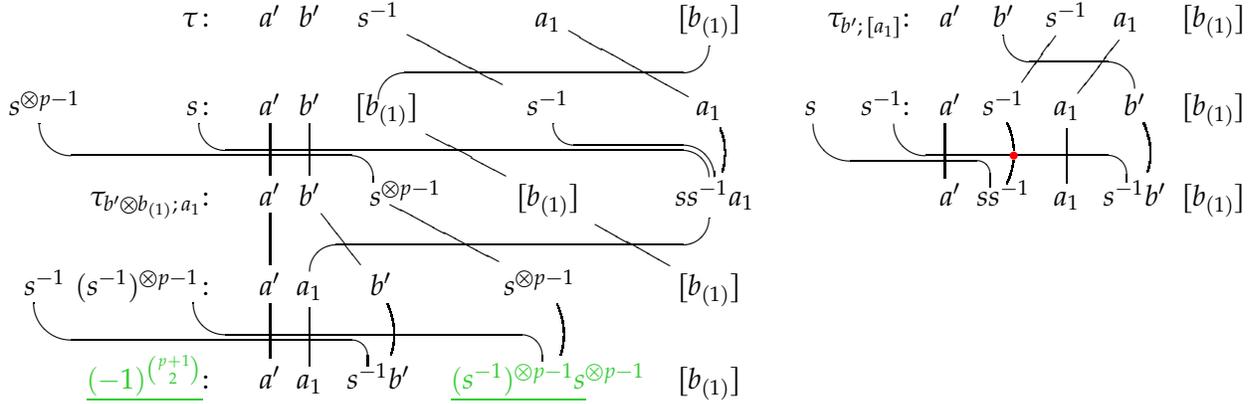
\[
	\xymatrix@R=1.25em@C=0em{
		&&\mathllap\tau\mathrlap{\:}	&\quad
		&a'	&b'	&\desusp\ar@{-}[dr]&a_1\ar@{-}[dr]&[b_{(1)}]\ar@{-} `d[l] `l [dll]&	
		&\ \  
		&&&\mathllap{\tau_{b';\, [a_1]}\mn}\mathrlap{\:}	&\	
		&a'&b'	\ar@{-} `d[r] `r [drr]&\desusp\ar@{-}[dl]&a_1	\ar@{-}[dl]\ \ &[b_{(1)}]	
		\\
		\susp^{\otimes p-1}\ar@{-}@<-2pt> `d[r] `r [drrrrrr]&&
		\mathllap\susp\mathrlap{\:}\ar@{-}@<0pt> `d[r] `r [drrrrrr]
		&\quad&a'\ar@{-}[d]&b'\ar@{-}[d]
		&\ \ [b_{(1)}]{\vphantom{X}}\ar@{-}[dr]
		&\ \ \desusp\ar@{-}@<2pt> `d[r] `r [dr]\ \ &a_1\ar@{-}@/^/[d]& 
		&\ \  
		&\susp\ \ \  \ar@{-}@<-5pt> 
		`d[r] `r [drrrrr]
		&\ \ \ 
		&\mathllap{\desusp\mn}\mathrlap{\:}\ar@{-}@<-3pt> 
		`d[r] `r [drrrrr]
		&\ &a'\ar@{-}@<-2pt>[d]
		&\desusp\ar@{-}@/^/@<-2pt>[d]|(.53){\textcolor{Red}{\bullet}}
		&a_1\ar@{-}[d]
		&b'\ar@{-}@/^/[d]
		&[b_{(1)}]
		\\
		&&\mathllap{\tau_{b' \otimes b_{(1)};\, a_1}\mn}\mathrlap{\:}
		&\quad
		&a'\ar@{-}[d]
		&{b'}\ar@{-}[dr]
		&\!\susp\mathrlap{{}^{\otimes p-1}}\ar@{-}[dr]
		&\ \ [b_{(1)}]\ar@{-}[dr]\ \ 
		&\!\!\!\susp\desusp \mathrlap{a_1}\ar@{-} `d[l] `l [dlll]& 
		&\ \  
		&&&&\	&\smash{a'}{\vphantom{t}}&\susp\smash{\desusp}&a_{\smash{1}}&\desusp \smash{b'}&[b_{(1)}]
		\\
		\desusp \ar@{-}@<-4pt> `d[r] `r [drrrrrr]
		& \qquad \quad
		&{\phantom{}}\smash{\mathllap{(\desusp)\vphantom{s}{}^{\otimes p-1}}}\mathrlap{\:}	\ar@{-}@<-2pt> `d[r] `r [drrrrr]
		&\quad
		&a'\ar@{-}[d]
		&a_1\ar@{-}[d]
		&\, b'\ar@{-}@/^/[d]&\susp^{\otimes p-1}\ar@{-}@/^/[d]\ \ &[b_{(1)}] 
		\\
		&&\mathllap{\textcolor{LimeGreen}{\ul{(-1)^{\binom{p+1} 2}}\mn}}
		\mathrlap{\:}
		&\quad&a'&a_1&\desusp b'&
		\ \ \textcolor{LimeGreen}{\ul{(\desusp)^{\otimes p-1}\susp\mathrlap{{}^{\otimes p-1}}}} \ \ \ \ \ \ \ \ \ &[b_{(1)}] 
	}
	\]
	In determining whether the signs agree, we can remove cancelling
	signs within one diagram or pairs of matching signs in both 
	diagrams, referring to this as ``cancellation'' either way.
	Then two crossings cancel if both involve transposing the same symbols
	$y$ and $z$ or if one crossing transposes $\susp$ and $z$
	while the other transposes $\desusp$ and $z$.
	After cancellation,
	the only signs remaining 
	are $(-1)^{|\desusp|\.|\desusp|} = -1$
	arising from the red crossing on the right and
	the underlined green
	\[
	(-1)^{\binom{p+1}{2}}\.|(\desusp)^{\otimes p-1}\susp^{\otimes p-1}| = (-1)^{\binom{p+1}{2}+\binom{p-1}{2}} = (-1)^{p^2 -p + 1	} = -1
	\] 
	at the bottom on the left,
	so the operations in question are indeed the same.

	As for $Q_3$ and $T_{0,3} = T_{1,3}$
	(from \eqref{eq:last-three-maps}, 
	\eqref{eq:two-sided-first-term},
	and \eqref{eq:two-sided-second-term}),
	we have seen above the terms on which the $\Phi_{\mn A''}$ 
	in $Q_3$ is potentially nonzero
	are those of the form 
	$[1 \otimes f''b]_{(3)} \otimes [a'' \otimes b'']$.
	In particular, $[a_{(3)}] = []$.
	Recall that we write~$\defm q = \ell{[b_{(3)}]}$.
	We see the many $1$ factors
	arising because~$\nabla\big([] \otimes [b_{(3)}]\big) = 
				[1 \otimes b]_{(3)}$
	contribute nothing to the iterated $\mu_{A''}$-product,
	so we may omit them with the map $\vk$.
	When they are omitted, the shuffle involved in the nonvanishing summand
	of \eqref{eq:Phi-sign}
	simplifies to 
	the shuffle $\tau_{f''b_{(3)};\,  a''}$  
	switching the $A''$ tensor factors containing~$a''$ and $f''b_{(3)}$
	in $f''b_{(3)} \otimes a'' \otimes b''$.
	Thus on $\BB_q A \ox A'' \ox A''$,
	we can expand $Q_3$ as
	\eqn{
		Q_3
		&= 
		t_{A''}\Phi_{\mn A''}
		\o
		\big(\BB(f'')^{\ox\mn 2} \otimes \desusp_{A''\ox A''}\big) 
		\\
		&=
		(-1)^{\binom{q+1}{2}} 
		(\Phi_{\mn A''})_{(1+q)}
		\o
		\susp_{A''\ox A''}^{\otimes 1+q}
		\o
		\big(\BB(f'')^{\ox\mn 2}\nabla\otimes 
		\desusp_{A'' \ox A''}\big)
		\\
		&=
		(-1)^{\binom{q+1}{2}} 
		(-1)^q
		\mu_{A''}
		\big(
		\EE_{1,q}(\desusp_{A''})^{\otimes 1+q} \otimes \id_{A''} 
		\big)
		\tau_{f''b_{(3)};\, a''}\vk
		\o
		\big(
		\susp_{A''\ox A''}^{\otimes q}
		\BB(f'')^{\ox\mn 2}\nabla \otimes 
		\id_{A'' \ox A''}
		\big)
		\\
		&=
		(-1)^{{\binom{q+1}{2}}-q} 
		\mu_{A''}
		\Big( 
		\EEE
		\big(
		\id_{A''} \otimes (\desusp_{A''})^{\otimes q}
		\big) 
		\otimes \id_{A''}
		\Big)
		\tau_{f''b_{(3)};\, a''}\vk
		\big(
		\susp_{A''\ox A''}^{\otimes q}
		\BB(f'')^{\ox\mn 2}\nabla \otimes 
		\id_{A'' \ox A''}
		\big)
		\mathrlap,
	}
	where as in the analysis of $Q_1$, 
	the $(-1)^{\binom{q+1}{2}}$
	comes from $(\Phi_{\mn A''})_{q+1} = 
	t_{A''}\Phi_{\mn A''}\susp_{A'' \ox A''}^{\otimes q+1}$
	and the $(-1)^q$ is
	the leftmost factor of the right-hand side of \eqref{eq:Phi-sign}.
	We may simplify the exponent of $-1$ by observing
	that 
	\[
	{\binom{q+1}{2}} - q = 
	{\binom{q+1}{2}} - {\binom{q}{1}}= 
	{\binom q 2}\mathrlap.
	\]
	In the diagram below the crossings for $Q_3$
	appear on the left and that for 
	\[
	T_{0,3} = T_{1,3} = 
	\mu_{A''}
	\big(\EEE(\id_{A''} \ox \BB f'') \ox \id_{A''}\big)
	\tau_{[b_{(3)}\mnn];\, a''}
	\]
	appears on the right. 
	\[
	\xymatrix@C=0em@R=1em{
		{\phantom{()^{-1}}}\susp^{\otimes q}\mathrlap{\:}\ar@{-}[drr]			&\quad	&&[b_{(3)}] 	& a'' 										& b''
		&\qquad\qquad&
		\tau_{[b_{(3)}\mnn];\, a''}\mathrlap{\:}	&		&[b_{(3)}] \ar@{-}[dr]	& a''		\ar@{-}[dl] 				& b''
		\\
		\tau_{f''b_{(3)};\, a''}\mathrlap{\:}	&	&	\susp^{\otimes q} \ar@{-}[dr]&[b_{(3)}] \ar@{-}[dr]& a'' 	
		\ar@{-} `d[l] `l [dll]									& b''
		&\qquad\qquad&
		&\quad	&a''					& [b_{(3)}] 							& b''
		\\
		(\desusp)^{\otimes q}\mathrlap{\:}		\ar@{-}@<-1.5pt>
		`d[r] `r [drrr]		&&		a''\ar@{-}[d]						&\susp^{\otimes q}	\ar@{-}@/^/[d]&[b_{(3)}] 					& b''
		\\
		\textcolor{LimeGreen}{\ul{(-1)^{\binom q 2}}}\mathrlap{\:}	&&		a''				
		&
		\textcolor{LimeGreen}{\ul{(\desusp)^{\otimes q}\susp^{\otimes q}}}&[b_{(3)}] & b''
	}
	\]
	When matching crossings are cancelled,
	the only signs remaining are
	the underlined
	$(-1)^{\binom q  2}$ and~
	$|(\desusp)^{\otimes q} \susp^{\otimes q}| = (-1)^{\binom q 2}$
	on the left,
	so the operations are equal, concluding the proof.
\end{proof}

With this more explicit formulation of $\wt\mu$, 
we are able to relate it to the product on an \HGA.

\begin{theorem}\label{thm:two-sided-ring-map}
	Given a commutative diagram 
	\[
	\xymatrix{
		&A\ar[dl]_{f'}\ar[dr]^{f''}
		\ar[dd]|(.45)\hole|{g}|(.55)\hole\\
		A'\ar[dr]_{g'}	&							&A''\ar[dl]^{g''}\\
		&\wt A
	}
	\]
	of \HGA homomorphisms, 
	the natural cochain map
\begin{equation}
\begin{aligned}\label{eq:def:two-sided-cochain-map}
	\defm\xi\ceq \iter{\mu_{\wt A}}{3} 
	(g' \otimes \hA\e_{\BB A} \otimes g'')\:
	\BB(A',A,A'') &\,\lt\, \wt A\mathrlap{,}\\
	a'[a_\bl]a'' &\,\lmt\, g'(a')
	\h_{\wt A}\e_{\BB A}[a_\bl] g''(a'')\mathrlap{.}
\end{aligned}
\end{equation}
	induces a map in cohomology
	multiplicative with respect to the 
	product on $\H \B(A',A,A'')$
	induced by the product defined by \Cref{thm:def-prod-bar}
	and the expected product on $\H(\wt A)$.
\end{theorem}
\begin{proof}
	By \Cref{thm:def-prod-bar},
	$\B(g',g,g'')$
	is a multiplicative cochain map $\BB(A',A,A'') \lt \BB(\wt A,\wt A,\wt A)$,
	so it will be enough to 
	show $\xi\: \BB(\wt A,\wt A,\wt A) \lt \wt A$
	is a cochain map inducing a multiplicative map in cohomology.
	In other words, we may as well start by assuming $A' = A = A'' = \wt A$
	and~$f' = f'' = g' = g'' = \id_A$.
	Let us then agree to write $\defm \mu = \muA$
	and $\defm \id = \id_A$
	and $\defm \e = \hA\e_{\BB A}$,
	so we can restate
	\eqref{eq:def:two-sided-cochain-map}
	as
\begin{equation}\label{eq:def:two-sided-cochain-map-v2}
\begin{aligned}
	\xi = &\ \iter{\mu}{3} (\id \otimes \e \otimes \id )\:\\
	c'[c_\bl]c'' &\lmt c' \e[c_\bl] c''
	\mathrlap.
\end{aligned}
\end{equation}
	It is trivial to check this is a cochain map.

	To show the induced map $\H(\xi)$
	in cohomology is multiplicative,
	we show that
	\begin{equation}\label{eq:homotopy-verification}
	\Homd{h} = \xi\wt\mu - \mu \xi^{\ox\mn 2}\: 
	\BB(A,A,A)^{\ox\mn 2} \lt A
	\end{equation}
	for a certain homotopy $h\:\B(A,A,A)^{\ox 2} \lt A$
	to be produced momentarily.
	It will help to adapt  
	\eqref{eq:def:two-sided-product},
	\eqref{eq:two-sided-first-term},
	\eqref{eq:two-sided-second-term}
	to the present set-up:
\begin{equation}\label{eq:def:two-sided-product-v2-sign}
\begin{aligned}
	\wt\mu
	&=
	&	\Big(\ 	\use{\mathclap{a'b'}}{\mu} 
\ \ 
		\,\use\otimes\otimes\, 
	&		\phantom{{}_{(2)}}
		\use{[a_\bl] \bartimes [b_{(2)}]}{\mu_{\BB A}} 
\ \,
		\,\use\otimes\otimes\, 	
\ 	&		\mu\smash{\big(}
		\clone{\EEE\mn\smash{\big(}}
		\use{a''}\id \ox \use{[b_{(3)}]\mnn\smash{\big)}}
		{\id_{\BB A} \mn ) }
		\ox \,\mn\use{b''}\id\smash{\big)}
			\mn\Big)
		\,\o\,
		\Pi'\phantom{'}	\\
\phantom{\wt\mu}&+
&\Big(\mnn\mu\smash{\big(}
\use{a'}\id 
\clone\otimes 
\clone{\EEE\smash{\big(}} 
\use{a_1}{\!\susp} \otimes 
\use{\!\! [b'}{\desusp} 
\use{|}\otimes 
\use{b_{(1)} ] \mnn\smash{\big)}}{\id_{\BB A}}
\mn )\smash{\big)}
\ \ \,\clone\otimes\,
&
\phantom{{}_\bl}
\use{ [a_{(2)}] \bartimes [b_{(2)}] }{\mu_{\BB A}}
\ \,
\,\clone\otimes\,
\ &			\mu\smash{\big(}
		\clone{ \EEE\mn\smash{\big(} }
		\use{a''}{\id} 
		\otimes 
		\use{ [b_{(3)} ] \mnn\smash{\big)} }
			{\id_{\BB A}
		\mn ) }
		\otimes 
		\use{b''}{\id}
			\smash{\big)}
			\mn
		\Big)
	\o\,\mn
	\Pi''		\mathrlap,
\end{aligned}
\end{equation}
where $\Pi'$ and $\Pi''$ are preparatory tensor permutations.
The permutations $\tau$ in particular have been absorbed into these
without sign change because the operators between $\tau$
and $\Pi$ in \eqref{eq:two-sided-first-term} 
and \eqref{eq:two-sided-second-term}
are of degree zero.
	Here the value, up to sign (determined by the Koszul convention,
	but not written out), 
	on the 
	standard pure tensor $a'[a_\bl]a'' \otimes b'[b_\bl]b''$ 
	is displayed below the function, to make the argument easier to follow,
	though in principle these values are optional.
	We will continue with these notations throughout the proof.
Recall that if $b'$ is in the image of the unit map $\h_A$,
then the $\EEE$ factor in the second term is defined to vanish
on bar-words containing it by \Cref{def:vanishing-components}.

It is evident from
\eqref{eq:def:two-sided-cochain-map-v2}
that
$-\mu\xi^{\ox 2}$ vanishes unless $\defm \ell = \ell(a_\bl)$
and $\defm r = \ell(b_\bl)$ are both zero,
in which case
it is
\[
-\iter \mu 6 
	(
		\use{a'}\id \otimes \h_A \otimes \use{a''}\id 
			\otimes 
		\use{b'}\id \otimes \h_A \otimes \use{b''}\id
	)
=
-\use{a'a''b'b''}{\iter \mu 4}
\mathrlap.
\]

As $\mu_{\BB A}$ is of degree $0$ 
and $\e_{\B A}$ annihilates $\B_{\geq 1} A$,
the first term of $\xi\wt\mu$ vanishes,
according to \eqref{eq:def:two-sided-cochain-map-v2}
and \eqref{eq:def:two-sided-product-v2-sign}
unless $\ell(a_\bl) = \ell(b_{(2)}) = 0$, in which case
it 
reduces to 
\[
\iter\mu 3 \Big(\ 
			\use{\mathclap{a'b'}}\mu 
				\otimes 
			\h_A 
				\otimes  
			\mu\smash{\big(}
		\clone{ \EEE\mn\smash{\big(} }
			\use{a''}\id 
			\ox 
			\use{ [b_{(3)}] \mnn\smash{\big)} }
				{\id_{\B A} \mn ) }
			\ox \,\mn
			\use{b''}{\id}
			\smash{\big)}
		\mn\Big)
		\o\,	\tau_{a'';b'}
\mathrlap.
\]
If additionally $\ell(b_{(3)}) = 0 = r$,
this further reduces to 
\[
\iter\mu 3 \Big(\ 
	\use{\mathclap{a'b'}}\mu 
		\otimes 
	\h_A 
		\otimes  
		\mu\smash{\big(}
	\use{a''}\id \
		\ox \,\mn\use{b''}\id\smash{\big)}
		\mn\Big)
		\o\,	
	\tau_{a'';b'} 
	= 
	\iter \mu 4 \,\o\, \tau_{a'';b'}
\mathrlap.
\]
The second term 
	of $\xi\wt\mu$ 
	vanishes unless $[a_\bl] = [\defm a]$ is of length $1$ 
	and $[b_{(2)}]$ is of length $0$,
	in which case it contributes
\[
	\iter \mu 4
	\big(\use{a\mathrlap{'}}\id \otimes 
		\use{\EEE\mn \smash{\big(} a_1\!\!\phantom{\susp}}	
			{\EEE(\susp\phantom{a_1}\!\!} 
		\!\!
		\otimes 
		\use{[b'}\desusp 
		\use|\otimes 
		\use{b_{(1)} ]\mnn\smash{\big)}}{\id_{\B A}}) 
		\otimes 
		\use{\EEE\mn\smash{\big(}}{\EEE(} 
		\use{a''}\id \otimes 
		\use{ [b_{(3)}] \mnn\smash{\big)} }
			{\id_{\B A}) }
		\otimes
		\use{b''}{\id}
		\big)
		\o{}
		\pi''
\]
for $\pi''$ running over tensor permutations
	$a'[a]a'' \otimes b'[b_\bl]b''
	\lmt 
	a' \otimes [a]b'[b_{(1)}] \otimes a''[b_{(3)}] \otimes b''$.
These $\pi''$ are the specialization of
$\Pi''$ from \eqref{eq:def:two-sided-product-v2-sign}
to the case $\ell(a_{(2)}) = 0$.
	
All told, on $\B_{ \ell}(A,A,A) \ox \B_{ r}(A,A,A)$
one has
\begin{equation}\label{eq:xi-mu-mu-xi}
\begin{aligned}
	\xi\wt\mu - \mu_A \xi^{\ox\mn 2} 
	= \case{
		\use{a'b'a''b'}{\iter \mu 4 \tau_{a'';\, b'}}
		\
		 - 
		\!\!\!\use{a'a''b'b''}{\iter \mu 4}
		& \mbox{if }	\ell = 0 = r\mathrlap,
		\\
		\noalign{\vskip7.5pt}
		\iter \mu 4(
		\use{a'}\id \otimes 
		\use{b'}\id \otimes 
		\use{\smash{\EEE\big(a''[b_\bl]\big)}}\EEE \otimes 
		\use{b''}\id)\tau_{a'';\, b'}
		& \mbox{if }	\ell = 0 < r\mathrlap,\\
		\noalign{\vskip7.5pt}
		\iter\mu 4
				\big(\use{a'}\id \otimes 
		\use{\EEE\mn \smash{\big(} a_1\!\!\phantom{\susp}}	
			{\EEE(\susp\phantom{a_1}\!\!} 
		\!\!
		\otimes 
		\use{[b'}\desusp 
		\use|\otimes 
		\use{b_{(1)} ]\mnn\smash{\big)}}{\id_{\B A}}) 
		\otimes 
		\use{\EEE\mn\smash{\big(}}{\EEE(} 
		\use{a''}\id \otimes 
		\use{ [b_{(3)}] \mnn\smash{\big)} }
			{\id_{\B A}) }
		\otimes
		\use{b''}{\id}
		\big)
		\o
		\pi''
		& \mbox{if }	\ell = 1\mathrlap,\\
		\noalign{\vskip7.5pt}
		0
		& \mbox{if }	\ell \geq 2\mathrlap.
	}
\end{aligned}
\end{equation}
%
%
The promised homotopy is given as
	\begin{equation}
	\begin{aligned}\label{eq:def:homotopy}
	\defm h 
	&\ceq 
	\iter{\mu}{4}
	\,\mn\smash{\big(}
	\use{a'}{\id} \otimes
	\use{\e[a_\bl]}{\e} \otimes 
	\use{\EEE}{\EEE}
	\use{\smash{\big(}\mn a''}{(\id} \otimes 
	\use{[b'}{\desusp} 
	\use|{\otimes} 
	\use{b_\bl]\mnn\smash{\big)}}{\id_{\BB A}}) \otimes 
	\use{b''}{\id\vpe}
	\smash{\big)}
	\end{aligned}
	\end{equation}
	To show $\Homd h$ agrees with $\xi\wt\mu - \mu_A \xi^{\ox\mn 2}$, 
	we follow the same case distinctions.
	
	If $\ell = \ell(a_\bl) \geq 2$, then $h$ and hence $dh$
	vanish since $\e[a_\bl] = 0$, 
	and since the differential on~$\B(A,A,A)$
	takes $\B_{\geq 2}(A,A,A)$ into $\B_{\geq 1}(A,A,A)$,
	so does $hd_{\B(A,A,A)^{\ox 2}}$.
	
	If $\ell = 1$, then $h$ and hence $dh$ vanishes
	as before,
	but $hd$ need not, 
	by the formula for the ``external'' differential
	$\dextBAA = d_{\B(A,A,A)} - d_{\ox}$ given in \eqref{eq:def:dext-two-sided}, 
	since the outer two operators reduce the length of
	the bar-word in $a'[a]a''$ to zero, taking it respectively
	to $\pm a'a[]a''$ and $\pm a'[]aa''$.
	Thus~$hd_{\otimes}$ and $h(\id_{\B(A,A,A)} \otimes d_{\B(A,A,A)})$
		vanish,
	but for $h(d_{\B(A,A,A)} \ox \id_{\B(A,A,A)})$,
	plugging \eqref{eq:def:dext-two-sided} into \eqref{eq:def:homotopy},
	we get
	\begin{equation}
	\begin{aligned}\label{eq:ell-equals-one}
&
	\phantom{-{}}
	\iter{\mu}{4}
	\,\mn\smash{\big(}
	{\id} \otimes
	{\e} \otimes 
	{\EEE}
	{(\id} \otimes 
	{\desusp} 
	{\otimes} 
	{\id_{\BB A}}) \otimes 
	{\id}
	\smash{\big)}
	&\o\ &
	\Big( \mn
		\big(\mu(\use{a'}\id \ox \use{a}\susp) \ox \use{[]}{\h_{\B A}} \ox \use{a''}\id\big)\ox \use{b'[b_\bl]b''}{\id_{\B(A,A,A)}}
\\
&&
&-		\use{a'}\id \ox \use{[]}{\h_{\B A}} \ox \mu(\use{a}\susp \ox \use{a''}\id)\big) \ox \use{b'[b_\bl]b''}{\id_{\B(A,A,A)}}\Big)
	\\
	=&-\,\iter {\mu} 4
	\big(
	\use{a'}\id \otimes
	\use{a}\susp \otimes
	\use{\EEE\mn\smash{\big(a''}}
		{ \EEE(\id } 
	\otimes 
	\use{[b'}{\desusp} 
	\use|\otimes 
	\use{b_{(1)}]\mnn\smash{\big)}}{\id_{\BB A})} \otimes
	\use{b''}\id
	\big)
	&&
	\\&
	+ \,\iter {\mu} 3
	\Big(
	\use{a'}\id \otimes
	\use{\EEE\big(}
		{\EEE\big(}
	\mu(\use{a}\susp \otimes 
	\use{a''}\id) \otimes
	\use{[b'}{\desusp} 
	\use|\otimes 
	\use{b_{(1)}]\mnn\smash{\big)}}{\id_{\BB A})} \otimes
	\use{b''}\id
	\Big)
	\mathrlap,&&
	\end{aligned}
	\end{equation}
	the change in signs coming in both terms from the 
	factor $\susp$ of degree $1$ 
	coming from $\dextBAA$
	being moved past the factor
	$\desusp$
	of degree $-1$ in the argument of $\EEE$.
	By the Cartanesque formula~\eqref{eq:EEE-Cartan},
	the second term of \eqref{eq:ell-equals-one}
	can be replaced by
	\[
	\iter \mu 4
	\big(
	\use{a'}\id 
	\otimes
	\,
	\use{(\EEE\otimes\EEE)\pi\big(}{(\EEE\otimes\EEE)\pi(}
	\use{a}\susp \otimes 
	\use{a''}{\id}
	\otimes 
	\use{[b'}\desusp 
	\use|\otimes 
	\use{b_\bl]\mn\big)}{\id_{\BB A})}
	\otimes 
	\use{b''}\id 
	\big)
	\mathrlap,
	\]
	where $\defm{\pi}$ runs over shuffles
	$
	a \otimes a'' \otimes [b'|b_\bl]
	\lmt 
	a [b'|b_{\bl}]_{(1)} \otimes a''[b'|b_{\bl}]_{(2)}
	$
	and $[b'|b_{\bl}]_{(1)} \otimes [b'|b_{\bl}]_{(2)}$
	is a deconcatenation of $[b'|b_\bl]$.
	There is a case distinction here:
	the deconcatenation yielding $[] \otimes [b'|b_\bl]$ 
	recovers
	\[
	\iter {\mu} 4
	\big(
	\use{a'}\id \otimes
	\use{a}\susp \otimes
	\use{\EEE\mn\smash{\big(a''}}
		{ \EEE(\id } 
	\otimes 
	\use{[b'}{\desusp} 
	\use|\otimes 
	\use{b_\bl]\mnn\smash{\big)}}{\id_{\BB A})} \otimes
	\use{b''}\id
	\big)
	\mathrlap,
	\]
	cancelling the first term of \eqref{eq:ell-equals-one}.
	Otherwise, $b'$ is assigned to the first word
	in the deconcatenation,
	so the value is $[b'|b_{(1)}] \otimes [b_{(3)}]$,
	and one gets operations
\[
			\big(\use{a'}\id \otimes 
		\use{\EEE\mn \smash{\big(} a_1\!\!\phantom{\susp}}	
			{\EEE(\susp\phantom{a_1}\!\!} 
		\!\!
		\otimes 
		\use{[b'}\desusp 
		\use|\otimes 
		\use{b_{(1)} ]\mnn\smash{\big)}}{\id_{\B A}}) 
		\otimes 
		\use{\EEE\mn\smash{\big(}}{\EEE(} 
		\use{a''}\id \otimes 
		\use{ [b_{(3)}] \mnn\smash{\big)} }
			{\id_{\B A}) }
		\otimes
		\use{b''}{\id}
		\big)
		\o
		\pi''
%
%
\]
	agreeing with the case $\ell = 1$
	of \eqref{eq:xi-mu-mu-xi}.
	Because $\pi$ preserves the relative positions of 
	$\susp$ and $\desusp$,
	it does not incur a sign change in moving past them
	to become $\pi''$.

%
	For the $\ell = 0$ cases,
	we write the differential on $\BB(A,A,A)^{\ox\mn 2}$
	as the sum of the ``internal'' tensor differential $d_{\ox}$
	and the ``external''
	differentials $\id_{\BB(A,A,A)} \otimes \dextBAA$
	and $\dextBAA \otimes \id_{\BB(A,A,A)}$.
	We may omit the second external differential
	in consideration of $hd_{\B(A,A,A)}$
	since~$d_{\B(A,A,A)}$ vanishes on $\B_0(A,A,A)$,
	and when $r = 0$, we may omit the first external differential
	for the same reason.
	Write $\defm{\Homdp}$
	for the differential	
	on $\Hom\mn\big(\mnn(A \ox \BB A \ox A)^{\ox\mn 2},A\big)$,
	where the domain is $\B(A,A,A)^{\ox \mn 2}$ equipped with the tensor
	differential $d_{\ox}$.
	Thus $\Homdp h = \dA h + h d_{\otimes}$.
	Because composition is a cochain map,
	as are $\mu_{A}$, $\id_A$, $\e_{\BB A}$, $\id_{\BB A}$, $\desusp$,
	we have
	\begin{equation}\label{eq:d-prime-h-brief}
	\Homdp h
	= 
	\iter {\mu} 4
	\big(
	\use{a'}\id \otimes 
	\h_{\B A} \otimes 
	\use{(\Homdp{\EEE})\mn\big(}
		{(\Homdp{\EEE})(}
	\use{a''}\id \otimes 
	\use{[b'}\desusp 
	\use|\otimes 
	\use{b_\bl]\mn\big)}{\id_{\BB A})} 
	\otimes
	\use{b''}\id
	\big)
			\mathrlap.
	\end{equation}
	Equation \eqref{eq:d-prime-EEE} shows
	that for fixed $\defm r$,
	the function
	$\Homdp \EEE\:\bar A \otimes \BB_{1+r} A \lt A$
	is given by
	\begin{equation}\label{eq:d-prime-EEE-values}
\begin{aligned}
\phantom{.}
	\Homdp\EEE 
	&	\,=\, 
	\mu(
	\use{b'}\susp 
	\otimes 
	\use	{\EEE	 \mathrlap{ \mn\big(a'' [b_\bl] \mnn\big) } }
		{\EEE	 \mathrlap{	)(1 \ 2) } }
	\phantom{ )(1 \ 2) }
	{\quad}\,+\,
	\use{\EEE\big(a'' \,\dext[b'|b_\bl]\mnn\big)}
	{\EEE(\id \otimes \dext)} 
{}	\,-\,\quad
	\use{
	\mathrlap{
		\case{
			\EEE\big(a''[]\big)b'& r = 0,
			\\
			\EEE\big(a''[b'|b_1|b_\bl|b_{r - 1}]\big)b_r	& r > 0
			}
	}
	}
	{\mathrlap{\mu(\EEE \otimes s).}
	}
\phantom{	\EEE\big(a''[b'|b_1|b_\bl|b_{r - 1}]\big)b_r	\quad r > 0	}
\end{aligned}
\end{equation}
%
%
%
%
%
%
%
%
%
%
The second term also vanishes when $r = 0$,
and then plugging
\eqref{eq:d-prime-EEE-values}
into
\eqref{eq:d-prime-h-brief}
yields
\begin{equation}\label{eq:d-prime-h-r-0}
	\begin{aligned}
	\Homdp h
	=&
	\phantom{-{}}
	\iter \mu 4\big(
	\use{a'}\id \otimes 
	(\use{b'}\susp \otimes 
	\use{\EEE\mn(a''[]\big)}{\EEE)}\tau_{a'';\, [b']}
	(\id \otimes \desusp \otimes \h_{\B A}) \otimes 
	\use{b''}\id
	\big)
	&=
	\use	{a'b'a''b''}
	{	\iter \mu 4	\o \tau_{a'';\,b'}}
	\\
	\noalign{\vskip5pt}
		&
	-
	\iter \mu 4\big(
	\use{a'}\id \otimes 
	\use{\EEE\mn(a''[]\big)}{(\EEE\phantom{)}} \otimes 
	\use{b'}\susp)
	(\id \otimes \desusp \otimes \h_{\B A})
	\otimes 
	\use{b''}\id
	\big)
	&
	\use	{a'a''b'b''}
		{\mathllap{-}	\iter \mu 4	\mathrlap.}
	\end{aligned}
	\end{equation}
Here we use that $\EEE = \id$ on $A \ox \B_0 A \iso A$
and the further simplification in the first term comes from 
the calculation
	\[(\susp \otimes \id)\tau_{a'';\, [b']}(\id \otimes \desusp)
	= \tau_{a'';\, b'}(\id \otimes \susp)(\id \otimes \desusp) = \tau_{a'';\, b'}\mathrlap.
	\]
The $\ell = 0 = r$ case of \eqref{eq:xi-mu-mu-xi}
agrees with \eqref{eq:d-prime-h-r-0},
concluding that case.

For $\ell = 0 < r$, substituting \eqref{eq:d-prime-EEE-values}
into
\eqref{eq:d-prime-h-brief} now gives instead
	\begin{equation}\label{eq:d-prime-h-r-nonzero}
	\begin{aligned}
	\Homdp h
	&\,=\,\mn 
	\iter \mu 4\big(
	\use{a'}\id \otimes (
	\use{b'}\susp \otimes 	
	\use	{\EEE	 \mathrlap{ \mn\big(a'' [b_\bl] \mnn\big) } }
		{\EEE	}	)\tau_{a'';\, [b']}  
	(\id \otimes \desusp \otimes \id_{\BB_r A}) 
	\otimes 
	\use{b''}\id\big)
	\\	\noalign{\vskip5pt}
	&\,+\,
	\iter \mu 3\big(
	\use{a'}\id \otimes 
	\use{
		\EEE
		\mathrlap{
			\mn\big(a''[	b'b_1|b_2|\cdots|b_r]\mnn\big)
		}
	}{\EEE}
	(\id \otimes \desusp\mu\susp^{\ox\mn 2} \otimes \id_{\BB_{r-1} A})
	(\id \otimes \desusp \otimes \id_{\dbA} \otimes \id_{\BB_{r-1} A})
	\otimes \use{b''}\id
	\big)\\	\noalign{\vskip5pt}
	&\,+\,
	\iter \mu 3\big(
	\use{a'}\id \otimes 
	\use{\EEE\big(a'' }
	{\EEE(\id} \otimes 
	\use{[b']}\id 
	\clone{\otimes} 
	\use{\dext[b_\bl]\mnn\big)}
	{\dext|_{\BB_r A})}
	(\id \otimes \desusp \otimes \id_{\BB_r A})
	\otimes \use{b''}\id
	\big)\\	\noalign{\vskip5pt}
	&\,-\,
	\iter \mu 4\big(
	\use{a'}\id \otimes 
	\use
		{\EEE\big(a''[b'|b_1|\cdots|b_{r - 1}]\big)}
		{(\EEE}
	 \otimes \use{b_r}{\susp)}
	(\id \otimes \desusp \otimes \id_{\BB_{r-1} A} \otimes \id_{\dbA})
	\otimes \use{b''}\id
	\big)
	\mathrlap.
	\end{aligned}
	\end{equation} 
	As for $h(\id_{\BB(A,A,A)} \otimes \dextBAA)$,	
	composing its definition
	from \eqref{eq:def:dext-two-sided} 
	with the expression for $h$ from \eqref{eq:def:homotopy}
	yields 
	\begin{equation}\label{eq:h-id-dext}
	\begin{aligned}
	\iter {\mu} 3 
	&
	\Big(
	\use{a'}\id \otimes
	\clone{\EEE\big(}
	\use{a''}\id \otimes 
	\use{[b'b_1}\desusp\mu(\id \otimes \susp) 
	\use|\otimes 
	\use{\cdots\,|b_r]\mnn\big)}{\id_{\BB_{r-1} A}\big)}
	\otimes \use{b''}\id
	\Big)\\	\noalign{\vskip5pt}
	+\,\iter {\mu} 3 
	&
	\Big(
	\use{a'}\id \otimes
	\clone{\EEE\big(}
	\use{a''}\id \otimes 
	\use{[b']}\desusp 
	\clone\otimes 
	\use{\dext[b_\bl]\mnn\big)}{\dext|_{\BB_r A}\big)}
	\otimes \use{b''}\id
	\Big)\\	\noalign{\vskip5pt}
	-\,\iter {\mu} 3 
	&
	\Big(
	\use{a'}\id \otimes
	\clone{\EEE\big(}
	\use{a''}\id \otimes \use{[b'}\desusp 
	\use|\otimes 
	\use{b_1| \cdots | b_{r-1}]\big)}{\id_{\BB_{r-1} A}}\big)
	\otimes \mu(
	\use{b_r}\susp \otimes \use{b''}\id)
	\Big)
	\mathrlap.
	\end{aligned}
	\end{equation}
	These three terms are the respective opposites of the second through fourth
	terms of \eqref{eq:d-prime-h-r-nonzero},
	as transitioning from one 
	ordering of symbols to the other involves moving $\desusp$ past 
	$\susp$ in the second and fourth terms and past $\dext$
	in the third, incurring each time a sign change of $(-1)^{-1 \cdot 1}$.
	Thus these three pairs of terms cancel.
	It remains only to check the uncancelled first term of
	\eqref{eq:d-prime-h-r-nonzero}
	agrees with the $\ell = 0 < r$ clause in \eqref{eq:xi-mu-mu-xi}.
	For this, note that $\id$ has degree $0$,
	so we have~$\tau_{a'';\,[b']}(\id \ox \desusp) = 
				(\desusp \ox \id)\tau_{a'';\,b'}$,
	and $\EEE$ has degree zero as well, so 
	
\eqn{
\iter \mu 4\big(
	\id \otimes (
	\susp \otimes 	
		{\EEE	}	)\tau_{a'';\, [b']}  
	(\id \otimes \desusp \otimes \id_{\BB_r A}) 
	\otimes 
	\id\big)
&=
\iter \mu 4\big(
	\id \otimes (
	\susp \otimes 	
		{\EEE	})
	(\desusp \otimes \id  \otimes \id_{\BB_r A}) 
	\otimes 
	\id\big)
	\tau_{a'';\, b'}  
\\&=
\iter \mu 4(
		\id \otimes 
		\id \otimes 
		\EEE \otimes \id
		)
		\tau_{a'';\, b'}
\mathrlap.\tag*{\qedhere}
}
\end{proof}

\begin{remark}
	We expect that 
	if $A'$, $A$, and $A''$ are extended \HGAs,
	the product $\wt\mu$
	of \Cref{thm:def-prod-bar}
	is the $2$-component of a differential on $\B\B(A',A,A'')$
	making $\B(A',A,A'')$ an \Ai-algebra
	and making the map $\xi$
	of \Cref{thm:two-sided-ring-map} 
	the extended $1$-component $\smash{\onecomponent{\Xi}}$ of an
	\Ai-map $\Xi$ from $\B(A',A,A'')$ to $\smash{\wt A}$,
	but we will leave the exploration of this possibility
	for another occasion.
\end{remark}

We now can use the new product to establish the particular version of the Eilenberg--Moore
theorem we will need.


\begin{theorem}[Eilenberg--Moore with induced product]%
	\label{thm:EM-product}\label{thm:same-product}
	Suppose the coefficient ring $\kk$ is a principal ideal domain
	and suppose given a pullback square 
	\[
	\xymatrix{
		Y \ar[r]^{\defm\b}\ar[d]_{\defm\a}			& E\ar[d]\\
		X \ar[r]									& B
	}
	\]
	of pointed topological spaces in which $E \to B$ is a Serre fibration, 
	$B$ and $X$ are path-connected,
	the action of $\pi_1(B)$ on the cohomology of
	the fiber $F$ is trivial, 
	and (1)
	each $H^n(F)$ is a finitely generated $\kk$-module
	or (2) each $H^n(B)$ and each $H^n(X)$ is
	a finitely generated $\kk$-module.
	Then there is a natural quasi-isomorphism
	\eqn{
		\defm\xi\:
		\BB\big(\C(X),\C(B),\C(E)\mnn\big) 
		&\lt \C(Y)\mathrlap,\\
		x[b_\bl]e
		&\lmt \a^*(x) \cup 
		\h_{\C(Y)}\e_{\BB \C(E)}[b_\bl] \cup 
		\b^*(e)\mathrlap,
	}
	inducing a ring isomorphism
	\[
	\Tor^*_{\C(B)}\mn\big(\mnn\C(X),\C(E)\mnn\big) 
	\isoto 
	\H(Y)
	\]
	with respect to the product on the domain
	induced	
	by the product given in \Cref{thm:def-prod-bar}.
\end{theorem}

We rely on the presentation of Gugenheim and May~\cite[Thm.~3.3]{gugenheimmay}
with some modification of hypotheses.

\begin{proof}
	We assume $B$ and $X$ are path-connected to guarantee
	that the homotopy type of the fiber $F$ be well defined 
	and we only have to discuss 
	one group $\pi_1(B)$;
	one otherwise needs a separate argument for each path-component.
	This assumption then also implies
	$\pi_1(X)$ acts trivially on $\H(F)$.
	Given any proper projective resolution $P^\bl$
	of $\C(X)$ as a $\C(B)$-module 
	(or resolution in the more general sense of Gugenheim--May%
	~\cite[Defs.~1.1]{gugenheimmay}), then,
	Gugenheim--May show
	the expected composite filtered map~$\defm\vt\: 
	P^\bl \ox_{\C(B)} \C(E) \lt \C(X) \ox_{\C(B)} \C(E) \to \C(Y)$
	is a quasi-isomorphism~\cite[Thm.~3.3, Cor.~3.5]{gugenheimmay}
	under slightly different finiteness hypotheses, 
	namely that $\kk$ is Noetherian and the groups
	$H_n(X;\Z)$ and $H_n(B;\Z)$ are all finitely generated.
	
	These hypotheses 
	are used only to see the canonical maps
	$\C(B) \ox \H(F) \lt \C\big(B;\H(F)\big)$
	and
	$\C(X) \ox \H(F) \lt \C\big(X;\H(F)\big)$
	are quasi-isomorphisms,
	using a lemma~\cite[Lem.~3.2]{gugenheimmay}
	asserting that if $\kk$ is a commutative Noetherian ring, 
	$G$ a $\kk$-module,
	and $C$ a chain complex over $\Z$
	with each $H_n(C)$ finite,
	then $\Hom_\Z(C,\kk) \otimes_\kk G \lt \Hom_\Z(C,G)$
	is a quasi-isomorphism.
	The proof of this lemma
	uses only that $\Z$ is 
	a principal ideal domain,
	and hence the same argument shows
	that if $\kk$ is a principal ideal domain,
	$C$ now a differential graded $\kk$-module,
	and $G$ a $\kk$-module, 
	then $\Hom_\kk(C,\kk) \otimes_\kk G \lt \Hom_\kk(C,G)$
	is a quasi-isomorphism.
	Thus, assuming $\kk$ is a principal ideal domain, 
	we can replace Gugenheim--May's hypothesis
	that each integral homology group $H_n(B;\Z)$ and $H_n(X;\Z)$ 
	is finitely generated 
	with the weaker hypothesis (2) that the (co)homology
	groups with coefficients in $\kk$ are.
	If, alternatively, we assume (1)
	that the $H^n(F)$ is finitely generated over $\kk$,
	then 
	again assuming $\kk$ is a principal ideal domain, 
	the decomposition of each $H^n(F)$ 
	as a finite product of cyclic $\kk$-modules
	shows 
	$\C(B) \ox H^n(F) \lt \C\big(B;H^n(F)\big)$
	and~$\C(X) \ox H^n(F) \lt \C\big(X;H^n(F)\big)$
	are isomorphisms of differential graded $\kk$-modules.	

	To see the suppressed details 
	in the proof of the multiplicativity of $\H(\vt)$
	with respect to the classical product on Tor%
	~\cite[p.~26]{gugenheimmay},
	subdivide the vertical cohomological Eilenberg--Zilber map
	featuring along 
	the upper-right of McCleary's diagram gesturing at such a proof%
	~\cite[p.~255]{mcclearyspectral} as
	\[
	\H(Y) \ox \H(Y)
	\xtoo{\H(\Phi)}
	\H\big( C^*(Y) \ox C^*(Y) \big)
	\os{i}\longto
	\H\Big(\mn\big( C_*(Y) \ox C_*(Y) \big)\mn{\vphantom{)}}^*\mn\Big) 
	\os{\ \mr{EZ}^*\!}\longfrom
	\H(Y \x Y)\mathrlap. 
	\]
	Then there are evident horizontal maps subdividing the region 
	into three rectangles that can be seen 
	to commute on choosing maps between resolutions%
	~\cite[Thm.~1.7]{gugenheimmay} 
	and expanding out the definition of the external product.
	
	As $\kk$ is a principal ideal domain,
	\Cref{thm:bar-computes-Tor} 
	implies that
	$\defm{\wt C} = \B\big(\C(X),\C(B),\C(B)\big)$
	gives a resolution $\wt C \lt \C(X)$ of $\C(X)$ as a differential $\C(B)$-module
	and $\smash{\H(\wt C)}$ is the desired Tor.
	The map $\vt$ then specializes to our $\xi$.
	As $\C$ is a functor valued in  \HGAs,
	by \Cref{thm:def-prod-bar,thm:two-sided-ring-map}
	there is a natural cochain map
	$\wt\mu\: \wt C^{\ox\mn 2} \lt \wt C$
	such that 
	the quasi-isomorphism $\xi$ is 
	multiplicative up to homotopy
	in the sense that if $\mu$ is the cup product on $\C(Y)$,
	then there exists $h\: \wt C \lt \C(Y)$ with
	$\Homd h = \xi\wt\mu - \mu \xi^{\ox\mn 2}$.
	Thus the $\kk$-module isomorphism~$\H(\xi)\: \H(\wt C) \lt \H(Y)$ 
	takes the induced product
	\[
	\H(\wt C) \ox \H(\wt C) \os{\x}\lt \H(\wt C \ox \wt C) \xtoo{\H(\wt\mu)} \H(\wt C)
	\]  
	on Tor to the cup product on $\H(Y)$.
	Moreover, since the classical proof shows $\H(\xi)$ is multiplicative
	with respect the standard product on Tor
	and the cup product on $\H(Y)$,
	we conclude that $\H(\wt\mu) \o {\x}$ is the standard product.
\end{proof}


\begin{remark}\label{rmk:EM-original}
	The classical Eilenberg--Moore theorem
	obtains the ring structure
	on the domain 
	through the external product
	and the Eilenberg--Zilber theorem,
	without reference to a cochain-level product on 
	any complex computing Tor.
	We went to all this trouble 
	because in \Cref{sec:final},
	we need the fact 
	that the product on Tor 
	is induced by such a product 
	in order to show certain maps of Tors
	are multiplicative and finally obtain \Cref{thm:main}.
\end{remark}
%
%

\smallskip


{\footnotesize\bibliography{bibshort} }

\begin{thebibliography}{BaumS67}

\bibitem[BaMR14]{barthelmayriehl2014}
Tobias Barthel, J.~Peter May, and Emily Riehl.
\newblock Six model structures for {DG}-modules over {DGAs}: model category
theory in homological action.
\newblock {\em New York J. Math}, 20:1077--1159, 2014.
\newblock 
\href
{http://nyjm.albany.edu/j/2014/20-53.html}
{\path{nyjm.albany.edu/j/2014/20-53.html}},
\href
{http://arxiv.org/abs/1310.1159} {\path{arXiv:1310.1159}}.

\bibitem[Baum68]{baum1968homogeneous}
Paul~F. Baum.
\newblock {On the cohomology of homogeneous spaces}.
\newblock {\em Topology}, 7(1):15--38, 1968.
\newblock \href {http://dx.doi.org/10.1016/0040-9383(86)90012-1}
{\path{doi:10.1016/0040-9383(86)90012-1}}.

\bibitem[BF04]{bergerfresse2004operad}
Clemens Berger and Benoit Fresse.
\newblock Combinatorial operad actions on cochains.
\newblock {\em Math. Proc. Cambridge Philos. Soc.}, 137(1):135--174, 2004.
\newblock \href {http://arxiv.org/abs/0109158} {\path{arXiv:0109158}}, 
\href {http://dx.doi.org/10.1017/S0305004103007138}
{\path{doi:10.1017/S0305004103007138}}.

\bibitem[BH98]{bergerhuebschmann1998barW}
Clemens Berger and Johannes Huebschmann.
\newblock Comparison of the geometric bar and $w$-constructions.
\newblock {\em J. Pure Appl. Algebra}, 131(2):109--123, 1998.
\newblock \href {http://arxiv.org/abs/1303.4018} {\path{arXiv:1303.4018}},
\href {http://dx.doi.org/10.1016/S0022-4049(98)00020-6}
{\path{doi:10.1016/S0022-4049(98)00020-6}}.

\bibitem[Bo99]{boardman1999conditionally}
J.~Michael Boardman.
\newblock Conditionally convergent spectral sequences.
\newblock {\em Contemp. Math.}, 239:49--84, 1999.
\newblock  \url{http://hopf.math.purdue.edu/Boardman/ccspseq.pdf}, \href
{http://dx.doi.org/10.1090/conm/239/03597}
{\path{doi:10.1090/conm/239/03597}}.

\bibitem[Bor53]{borelthesis}
Armand Borel.
\newblock {Sur la cohomologie des espaces fibr{\'e}s principaux et des espaces
homog{\`e}nes de groupes de {L}ie compacts}.
\newblock {\em Ann. of Math. (2)}, 57(1):115--207, 1953.
\newblock 
\href
{http://web.math.rochester.edu/people/faculty/doug/otherpapers/Borel-Sur.pdf}
{\path{web.math.rochester.edu/people/faculty/doug/otherpapers/Borel-Sur.pdf}},
\href {http://dx.doi.org/10.2307/1969728} {\path{doi:10.2307/1969728}}.

\bibitem[CaF21]{carlsonfranzlong}
Jeffrey D.~Carlson (appendix joint~with Matthias~Franz).
\newblock The cohomology of biquotients via a product on the two-sided bar
construction (expository version).
\newblock 2021.
\newblock \href {http://arxiv.org/abs/2106.02986v1}
{\path{arXiv:2106.02986v1}}.

\bibitem[Car51]{cartan1950transgression}
Henri Cartan.
\newblock {La transgression dans un groupe de {Lie} et dans un espace fibr{\'e}
principal}.
\newblock In {\em {Colloque de topologie (espace fibr{\'e}s), {Bruxelles}
1950}}, pages 57--71, Li{\`e}ge/Paris, 1951. Centre belge de recherches
math{\'e}matiques, Georges Thone/Masson et companie.
\newblock  \url{http://eudml.org/doc/112227}.

\bibitem[Fr19]{franz2019homogeneous}
Matthias Franz.
\newblock The cohomology rings of homogeneous spaces.
\newblock 2019.
\newblock \href {http://arxiv.org/abs/1907.04777} {\path{arXiv:1907.04777}}.

\bibitem[Fr20a]{franz2019shc}
Matthias Franz.
\newblock Homotopy {Gerstenhaber} algebras are strongly homotopy commutative.
\newblock {\em J. Homotopy Relat. Struct.}, 15(3):557--595, 2020.
\newblock \href {http://arxiv.org/abs/1907.04778} {\path{arXiv:1907.04778}},
\href {http://dx.doi.org/10.1007/s40062-020-00268-y}
{\path{doi:10.1007/s40062-020-00268-y}}.

\bibitem[Fr20b]{franz2020szczarba}
Matthias Franz.
\newblock {Szczarba's} twisting cochain is comultiplicative.
\newblock 2020.
\newblock \href {http://arxiv.org/abs/2008.08943} {\path{arXiv:2008.08943}}.

\bibitem[GrR14]{grinbergreiner}
Darij Grinberg and Victor Reiner.
\newblock Hopf algebras in combinatorics.
\newblock 2014.
\newblock 
\href{http://www.cip.ifi.lmu.de/~grinberg/algebra/HopfComb-sols.pdf}
{\path{cip.ifi.lmu.de/~grinberg/algebra/HopfComb-sols.pdf}}, 
\href
{http://arxiv.org/abs/1409.8356} {\path{arXiv:1409.8356}}.

\bibitem[Gu60]{gugenheim1960brown}
Victor~K.A.M. Gugenheim.
\newblock On a theorem of {E.} {H.} {Brown}.
\newblock {\em Illinois J. Math.}, 4(2):292--311, 1960.
\newblock \href {http://dx.doi.org/10.1215/ijm/1255455870}
{\path{doi:10.1215/ijm/1255455870}}.

\bibitem[GuM]{gugenheimmay}
Victor~K.A.M. Gugenheim and J.~Peter May.
\newblock {\em On the Theory and Applications of Differential Torsion
Products}, volume 142 of {\em Mem. Amer. Math. Soc.}
\newblock Amer. Math. Soc., 1974.

\bibitem[GuMu74]{gugenheimmunkholm1974}
Victor~K.A.M. Gugenheim and Hans~J. Munkholm.
\newblock On the extended functoriality of {Tor} and {Cotor}.
\newblock {\em J. Pure Appl. Algebra}, 4(1):9--29, 1974.
\newblock \href {http://dx.doi.org/10.1016/0022-4049(74)90026-7}
{\path{doi:10.1016/0022-4049(74)90026-7}}.

\bibitem[Hopf41]{hopf1941hopf}
Heinz Hopf.
\newblock {{\"U}ber eie {Topologie} der {Gruppen-Mannigfaltigkeiten} und ihre
{Verallgemeinerungen}}.
\newblock {\em Ann. of Math. (2)}, 42(1):22--52, Jan 1941.
\newblock \href {http://dx.doi.org/10.2307/1968985}
{\path{doi:10.2307/1968985}}.

\bibitem[Hu89]{huebschmann1989perturbation}
Johannes Huebschmann.
\newblock Perturbation theory and free resolutions for nilpotent groups of
class 2.
\newblock {\em J. Algebra}, 126(2):348--399, 1989.
\newblock \href {http://dx.doi.org/10.1016/0021-8693(89)90310-4}
{\path{doi:10.1016/0021-8693(89)90310-4}}.

\bibitem[HuMS74]{husemollermoorestasheff1974}
Dale Husemoller, John~C. Moore, and James Stasheff.
\newblock Differential homological algebra and homogeneous spaces.
\newblock {\em J. Pure Appl. Algebra}, 5(2):113--185, 1974.
\newblock 
\href{http://math.mit.edu/~hrm/18.917/husemoller-moore-stasheff.pdf}
{\path{math.mit.edu/~hrm/18.917/husemoller-moore-stasheff.pdf}}.

\bibitem[May67]{maysimplicial}
J.~Peter May.
\newblock {\em Simplicial objects in algebraic topology}, volume~11 of {\em
Chicago Lect. Math.}
\newblock Univ. Chicago Press, 1992 (1967).
\newblock \url{https://www.maths.ed.ac.uk/~v1ranick/papers/maybook.pdf}.

\bibitem[May68]{may1968principal}
J.~Peter May.
\newblock The cohomology of principal bundles, homogeneous spaces, and
two-stage {Postnikov} systems.
\newblock {\em Bull. Amer. Math. Soc.}, 74(2):334--339, 1968.
\newblock  \href{http://www.math.uchicago.edu/~may/PAPERS/6.pdf}
{\path{math.uchicago.edu/~may/PAPERS/6.pdf}}.

\bibitem[May72]{mayloop}
J.~Peter May.
\newblock {\em The geometry of iterated loop spaces}, volume 271 of {\em
Lecture Notes in Math.}
\newblock Springer, 2006 (1972).
\newblock 
\href{http://guests.mpim-bonn.mpg.de/mpenney/operads/May_Iterated-loop-spaces.pdf}
{\path{guests.mpim-bonn.mpg.de/mpenney/operads/May_Iterated-loop-spaces.pdf}}.

\bibitem[MayN02]{mayneumanncohomologygeneralized}
J.~Peter May and Frank Neumann.
\newblock On the cohomology of generalized homogeneous spaces.
\newblock {\em Proc. Amer. Math. Soc.}, 130(1):267--270, 2002.
\newblock \href {https://doi.org/10.1090/S0002-9939-01-06372-9}
{\path{doi:10.1090/S0002-9939-01-06372-9}}.

\bibitem[Mc]{mcclearyspectral}
John McCleary.
\newblock {\em {A User's Guide to Spectral Sequences}}, volume~58 of {\em
{Cambridge Stud. Adv. Math.}}
\newblock Cambridge Univ. Press, Cambridge, 2001.

\bibitem[McS02]{mccluresmith2002deligne}
James~E. McClure and Jeffrey~H. Smith.
\newblock A solution of {Deligne’s} {Hochschild} cohomology conjecture.
\newblock In {\em Recent progress in homotopy theory ({Baltimore}, {MD},
2000)}, volume 293 of {\em Contemp. Math.}, pages 153--193, 2002.
\newblock \href {http://arxiv.org/abs/9910126} {\path{arXiv:9910126}}.

\bibitem[McS03]{mccluresmith2003}
James McClure and Jeffrey Smith.
\newblock Multivariable cochain operations and little $n$-cubes.
\newblock {\em J. Amer. Math. Soc.}, 16(3):681--704, 2003.
\newblock \href {http://arxiv.org/abs/0106024} {\path{arXiv:0106024}}, \href
{http://dx.doi.org/10.1090/S0894-0347-03-00419-3}
{\path{doi:10.1090/S0894-0347-03-00419-3}}.

\bibitem[Mi67]{milgram1967bar}
R.~James Milgram.
\newblock The bar construction and abelian {$H$}-spaces.
\newblock {\em Illinois J. Math.}, 11(2):242--250, 1967.
\newblock \href {http://dx.doi.org/10.1215/ijm/1256054662}
{\path{doi:10.1215/ijm/1256054662}}.

\bibitem[Mu74]{munkholm1974emss}
Hans~J. Munkholm.
\newblock {The {Eilenberg--Moore} spectral sequence and strongly homotopy
multiplicative maps}.
\newblock {\em J. Pure Appl. Algebra}, 5(1):1--50, 1974.
\newblock \href {http://dx.doi.org/10.1016/0022-4049(74)90002-4}
{\path{doi:10.1016/0022-4049(74)90002-4}}.

\bibitem[Sa09]{saneblidze2009bitwisted}
Samson Saneblidze.
\newblock The bitwisted cartesian model for the free loop fibration.
\newblock {\em Topology Appl.}, 156(5):897--910, 2009.
\newblock \href {http://arxiv.org/abs/0707.0614} {\path{arXiv:0707.0614}},
\href {https://doi.org/10.1016/j.topol.2008.11.002}
{\path{doi:10.1016/j.topol.2008.11.002}}.

\bibitem[Se68]{segal1968classifying}
Graeme Segal.
\newblock Classifying spaces and spectral sequences.
\newblock {\em Publ. Math. Inst. Hautes {\'E}tudes Sci.}, 34:105--112, 1968.
\newblock \url{http://numdam.org/item/PMIHES_1968__34__105_0/}.

\bibitem[Si93]{singhof1993}
Wilhelm Singhof.
\newblock {On the topology of double coset manifolds}.
\newblock {\em Math. Ann.}, 297(1):133--146, 1993.
\newblock \href {http://dx.doi.org/10.1007/BF01459492}
{\path{doi:10.1007/BF01459492}}.

\bibitem[Sm67]{smith1967emss}
Larry Smith.
\newblock Homological algebra and the {Eilenberg}--{Moore} spectral sequence.
\newblock {\em Trans. Amer. Math. Soc.}, 129:58--93, 1967.
\newblock \href {https://doi.org/10.2307/1994364} {\path{doi:10.2307/1994364}}.

\bibitem[StaH70]{halperinstasheff1970}
James Stasheff and Steve Halperin.
\newblock Differential algebra in its own rite [{\em sic}].
\newblock In {\em Proc. Adv. Study Inst. Alg. Top. (Aarhus 1970)}, volume~3,
pages 567--577, 1970.
\newblock 
\url {https://dropbox.com/s/zueujfbpt6g8ohw/Stasheff--Halperin_differential_algebra_own_rite.pdf?dl=0}.


\bibitem[Ste68]{steenrod1968milgram}
Norman~Earl Steenrod.
\newblock Milgram's classifying space of a topological group.
\newblock {\em Topology}, 7(4):349--368, 1968.
\newblock \href {http://dx.doi.org/10.1016/0040-9383(68)90012-8}
{\path{doi:10.1016/0040-9383(68)90012-8}}.

\bibitem[Wol73]{wolfthesis}
Joel~L. Wolf.
\newblock {\em The cohomology of homogeneous spaces and related topics}.
\newblock PhD thesis, Brown University, 1973.
\newblock 
  \url{https://jdkcarlson.github.io/cohomology_homogeneous_spaces(Wolf,DISSERTATION).pdf}.

\bibitem[Wolf77]{wolf1977homogeneous}
Joel~L. Wolf.
\newblock {The cohomology of homogeneous spaces}.
\newblock {\em Amer. J. Math.}, pages 312--340, 1977.
\newblock \href {http://dx.doi.org/10.2307/2373822}
{\path{doi:10.2307/2373822}}.

\end{thebibliography}

\nd\footnotesize{%
	\url{jeffrey.carlson@tufts.edu}
}

\medskip

\nd\footnotesize{%
	\textsc{Department of Mathematics\\
		Western University
	}\\
	\url{mfranz@uwo.ca}
}
\end{document}